\documentclass[11pt]{amsart}
\usepackage[top=1.15in, bottom=1.15in, left=1.17in, right=1.17in]{geometry}
\usepackage{fancyhdr}
{}
\usepackage[cp1252]{inputenc}
\usepackage{amsthm, amsmath, amssymb, amscd, latexsym, multicol, verbatim, enumerate, graphicx,xy, color}
\usepackage{amstext,amsfonts,amssymb,amscd,amsbsy,amsmath}

\usepackage{tikz,float}
\usetikzlibrary{patterns}
\usepackage[english]{babel}
\usepackage{cite}
\usepackage[colorlinks=true, pdfstartview=FitV, linkcolor=blue, citecolor=blue, urlcolor=blue, breaklinks=true]{hyperref}
\usepackage{helvet}         % selects Helvetica as sans-serif font
\usepackage{courier}        % selects Courier as typewriter font
\usepackage{type1cm}  
\usepackage{multicol}        % used for the two-column index
\usepackage{longtable}

\usepackage{hyperref}
\usepackage{graphicx}       

\makeindex             

\usepackage{amstext,amsfonts,amssymb,amscd,amsbsy,amsmath}
\usepackage[ruled,vlined]{algorithm2e}

\usepackage{algpseudocode}
\usepackage{tikz-cd}	
\usepackage{tikz}

\numberwithin{equation}{section}

\theoremstyle{definition}
\newtheorem{definition}{Definition}
\newtheorem{proposition}{Proposition}
\newtheorem{theorem}{Theorem}
\newtheorem{lemma}{Lemma}
\newtheorem{corollary}{Corollary}
\newtheorem{remark}{Remark}
\newtheorem{example}{Example}
\newtheorem{question}{Question}
\newtheorem*{question*}{Question}

\newtheorem*{corollary*}{Corollary}
\newtheorem*{lemma*}{Lemma}
\newtheorem*{theorem*}{Theorem}
\newtheorem*{proposition*}{Proposition}

\newcommand{\w}{\underline{w}}
\newcommand{\p}{\mathbf{p}}

\newcommand{\boundellipse}[3]% center, xdim, ydim
{(#1) ellipse (#2 and #3)
}

\newcommand{\X}{\mathcal{X}}
\newcommand{\A}{\mathcal{A}}

\newcommand{\area}{\mathsf{A}}

\newcommand{\pa}{\mathsf{pa}}

\DeclareMathOperator{\trop}{trop}

\DeclareMathOperator{\res}{res}
\DeclareMathOperator{\ord}{ord}
\DeclareMathOperator{\ind}{ind}

\newcommand{\Hom}{\operatorname*{Hom}}
\newcommand{\Lie}{\operatorname*{Lie}}
\newcommand{\Spec}{\operatorname*{Spec}}

\newcommand{\lie}{\mathfrak}

\newcommand{\mut}{{\operatorname*{mut}}}

\DeclareMathOperator{\init}{in}

\DeclareMathOperator{\conv}{conv}

\title[String cone and Superpotential]{String cone and Superpotential combinatorics for flag and Schubert varieties in type A}

\author[L.~Bossinger]{L.~Bossinger}
\address{L.~Bossinger: Instituto de Matem\'aticas UNAM Unidad Oaxaca\\
Le\'on 2, altos, Oaxaca de Ju\'arez\\
Centro Hist\'orico\\
68000 Oaxaca\\
Mexico}
\email{lara@im.unam.mx}

\author[G.~Fourier]{G.~Fourier}
\address{G.~Fourier: RWTH Aachen University, Pontdriesch 10--16, 52062 Aachen}
\email{fourier@mathb.rwth-aachen.de}

\thanks{The first author was supported by "Programa de Becas Posdoctorales en la UNAM 2018" Instituto de Matem\'aticas, Universidad Nacional Aut\'onoma de M\'exico}

\begin{document}

\begin{abstract}
We study the combinatorics of pseudoline arrangements and their relation to the geometry of flag and Schubert varieties. We associate to each pseudoline arrangement two polyhedral cones, defined in a dual manner. 
We prove that one of them is the weighted string cone by Littelmann and Berenstein-Zelevinsky. 
For the other we show how it arises in the framework of cluster varieties and mirror symmetry by Gross-Hacking-Keel-Kontsevich: for the flag variety the cone is the tropicalization of their superpotential while for Schubert varieties a restriction of the superpotential is necessary.

We prove that the two cones are unimodularly equivalent. As a corollary of our combinatorial result we realize Caldero's toric degenerations of Schubert varieties as GHKK-degeneration using cluster theory.
\end{abstract}

\maketitle

\section{Introduction}\label{chap:intro-old}
Toric degenerations of Grassmannians (and later of flag and Schubert varieties) have gained much attention from the first works of Hodge \emph{et al}, realized by Gonciulea and Lakshmibai \cite{GL96}. 
Today it is a fast growing subject with strong links to representation theory, symplectic geometry, differential geometry, algebraic geometry. 
In the past twenty something years several approaches to construct toric degenerations have been provided: representation theory using dual canonical bases \cite{Cal02}, or more general using birational sequences \cite{FFL15}, Newton-Okounkov bodies using valuations \cite{An13,KK12,Kav15}, tropical geometry using Gr\"obner theory \cite{M-S,KM16}, and last but not least, cluster varieties and superpotentials \cite{GHKK14,RW17}. 

We briefly recall two of these approaches in more detail focusing on the flag variety of full flags of subspaces in $\mathbb{C}^n$.

\subsection{String cones}
Littelmann in \cite{Lit98} and Berenstein-Zelevinsky in \cite{BZ01} introduced string cones to parametrize Lusztig's dual canonical basis of the quantum group.
The defining hyperplanes of the cones are described inductively by Littelmann and explicitly by Gleizer-Postnikov \cite{GP00} using \emph{pseudoline arrangements} (see Definition~\ref{def:pseusoline arr}). 
In these arrangements, associated to a fixed reduced expression of the longest Weyl group element, Gleizer and Postnikov describe the defining inequalities using the turning points of \emph{GP-paths} (see Definition~\ref{def:GPpath}).
This can be naturally generalized to Schubert varieties as done by Littelmann and generalized in Theorem~\ref{thm:wt GP is wt string} below.

We consider the flag variety as $SL_n/B$ and fix an embedding into projective space associated to a highest weight representation $V(\lambda)$ of $\lie{sl}_n$. 
Then the string cones also parametrize a basis of the homogeneous coordinate ring of the flag variety. In this setup, Caldero used the string parametrization and the approach via Rees algebras to define toric degenerations of flag and Schubert varieties for reduced expressions of corresponding Weyl group elements in \cite{Cal02}. 
The polytope associated to the toric fiber of the degeneration is the string polytope.

%%%%%%%%%%%%%%%%%% CLUSTER
\subsection{Cluster algebras}\label{subsec:cluster}
Cluster algebras were introduced in \cite{FZ02} by Fomin and Zelevinsky and quickly grew to become a research area on their own. They are commutative rings endowed with \emph{seeds} (maximal sets of algebraically independent generators) related by \emph{mutation} (local transformations exchanging one seed by another). At their origin they are closely related to the representation theory of finite dimensional algebras, but also many objects related to algebraic groups have a cluster structure. For example, the homogeneous coordinate ring of Grassmannians (see \cite{FZ02,Sco06}), double Bruhat cells (see \cite{BFZ05}), (partial) flag varieties (see \cite{GLS13}) or Richardson varieties (see \cite{Lec16}). 

A geometric approach to cluster algebras was introduced by Fock and Goncharov in \cite{FG06}. In this setting they work with \emph{cluster varieties}, schemes glued from algebraic tori (one for every seed) with gluing given by the birational transformations induced by mutation. They come in two flavours, $\mathcal A$- and $\mathcal X$-cluster varieties, one being the \emph{mirror dual} to the other as developed by Gross, Hacking, Keel and Kontsevich in \cite{GHKK14}. Among other things, they define \emph{$\vartheta$-bases} for cluster algebras and provide toric degenerations of (partial compactifications of) cluster varieties. 
The $\mathcal X$-cluster variety comes endowed with a Laurent polynomial, the \emph{superpotential}, whose tropicalization gives a polyhedral cone and a polytope as a slice of this cone. The superpotential polytope is the polytope associated to the special fibre of the toric degeneration.

\begin{question}\label{Q:cluster}
Can Caldero's toric degenerations of flag and Schubert varieties be recovered by tropicalizing a superpotential?
\end{question}

A first hint towards a positive answer to this question was given by Magee in \cite{Mag15}. He recovers the Gelfand-Tsetlin polytope\footnote{Littelmann showed in \cite{Lit98} that the Gelfand-Tsetlin polytope is unimodular equivalent to a certain string polytope.} as a superpotential polytope in a particular seed.
Further results in this direction are obtained by Genz-Koshevoy-Schumann in \cite{GKS} and \cite{GKS2}, who generalize Magee's result to flag varieties of simple, simply connected, simply laced algebraic groups. They recover the classical string and Lusztig parametrizations from the superpotential.

\subsection{Summary of the results}

In \cite{GP00} Gleizer and Postnikov use pseudoline arrangements associated to reduced expressions $\w_0$ of $w_0\in S_n$ and GP-paths in these to parametrize the inequalities for string cones $C_{\w_0}\subset \mathbb R^N$. 
We extend their result by adding weight inequalities encoded combinatorially in the pseudoline arrangement and obtain the \emph{weighted string cones} $\mathcal C_{\w_0}\subset \mathbb R^{N+n-1}$ as defined in \cite{Lit98}. 
Intersecting $\mathcal C_{\w_0}$ with the preimage of a weight $\lambda\in \mathbb R^{n-1}$ of an appropriate projection $\pi:\mathbb R^{N+n-1}\to \mathbb R^{n-1}$ yields the string polytope $\pi^{-1}(\lambda)\cap\mathcal C_{\w_0}$. 
Generalizing to arbitrary $w\in S_n$ and following Caldero \cite{Cal02} we obtain similarly the string cone, weighted string cone and string polytope for the Schubert variety $X(w)$.

We introduce a second polyhedral cone $S_{\w_0}\subset\mathbb R^N$ associated to a pseudoline arrangement in a dual way: the variables are indexed by the faces of the diagram as opposed to the vertices in case of the string cone.
From additional weight inequalities we get a weighted cone $\mathcal S_{\w_0}$ and from a second projection $\tau:\mathbb R^{N+n-1}\to\mathbb R^{n-1}$ we obtain polytopes $\tau^{-1}(\lambda)\cap\mathcal S_{\w_0}$ for $\lambda\in \mathbb R^{n-1}$.
As in the case of string cones we provide these also for arbitrary $w\in S_n$. 
 
The combinatorial result of our study is the following (see Theorem~\ref{thm:unimod}).

\begin{theorem}\label{thm:dual cones intro}
For every $\w\in S_n$, the two cones $\mathcal{C}_{\w}$ and $\mathcal{S}_{\w}$ are unimodularly equivalent and the lattice-preserving linear map is given by the duality of faces and vertices in the pseudoline arrangement. 
Moreover, this linear map restricts to linear bijections between the polytopes
$\pi^{-1}(\lambda) \cap \mathcal{C}_{\w} \cong \tau^{-1}(\lambda) \cap \mathcal{S}_{\w}$ and the cones $S_{\w} \cong C_{\w}$.
\end{theorem}

The cone $\mathcal S_{\w_0}$ appears in the framework of mirror symmetry for cluster varieties \cite{GHKK14}: 
denote by $B^-\subset SL_n$ the Borel subgroup of lower triangular matrices and by $U\subset B$ (resp. $U^-\subset B^-$) the unipotent radical with all diagonal entries being $1$. 
The double Bruhat cell $G^{e,w_0}= B^- \cap Bw_0B$ is an $\mathcal A$-cluster variety (see \cite{BFZ05}) and can be identified with an open subset of $Bw_0B/U$.
Let $\mathcal{X}$ the be dual of the $\mathcal{A}$-cluster variety $G^{e,w_0}$ and let $s_0=s_{\hat\w_0}$ be the seed of the cluster algebra $\mathbb{C}[G^{e,w_0}]$ corresponding to the reduced expression $\hat\w_0 = s_1 \, s_2s_1 \, \cdots s_{n-1}\cdots s_2s_1$. 
Let $W$ be the superpotential defined by the sum of the $\vartheta$-functions for  frozen variables in $s_0$ as introduced in \cite{GHKK14}. 
Then $W^{\trop}$ denotes the tropicalization of the superpotential. Magee has shown in \cite{Mag15} (see also Goncharov-Shen in \cite{GS15}) that 
\[
\mathcal{S}_{\hat\w_0} = \{ x \in \mathbb{R}^{N+n-1} \mid W^{\trop}\vert_{\mathcal{X}_{s_{0}}}(x)\geq 0 \} =: \Xi_{s_0}.
\]
We show that mutation of the pseudoline arrangement and hence of the cone $\mathcal{S}_{\w_0}$, is compatible with mutation of the superpotential \cite{GHK15} by introducing mutation of GP-paths. 
We obtain the following result\footnote{Genz-Koshevoy-Schumann \cite{GKS} obtain a similar result in the context of crystal graphs.} (see Corollary~\ref{cor:area is trop super}):
\begin{theorem}\label{thm:intro2}
Let $\w_0$ be an arbitrary reduced expression of $w_0 \in S_{n}$ and $s_{\w_0}$ be the seed corresponding to the pseudoline arrangement, $\mathcal{X}_{s_{\w_0}}$ the toric chart of the seed $s_{\w_0}$. Then 
\[
\mathcal{S}_{\w_0} =  
\Xi_{s_{\w_0}},
\]
the polyhedral cone defined by the tropicalization of $W$ expressed in the seed $s_{\w_0}$.
\end{theorem}

Consider $w \in S_{n}$ arbitrary and $\w$ a reduced expression of $w$. Let $W$ be as above and consider its restriction $\res_{\w}(W\vert_{\mathcal X_{s_{\w_0}}})$ to the mirror dual of the $\mathcal{A}$-cluster variety $G^{e,w}$.  
Let $s_{\w}$ be the corresponding seed in the cluster algebra. 
Then the tropicalization of the restriction yields again a cone $\Xi_{s_{\w}}$. 
The last result establishes an answer to the question above for Schubert varieties (see Theorem~\ref{thm: application unimod w}).

\begin{theorem}\label{thm:intro3}
Let $\w \in S_{n}$, and fix $\w_0=\w s_{i_{\ell(w)+1}}\dots s_{i_N}$ a reduced expression of $w_0\in S_n$. 
Let $s_{\w}$ resp. $s_{\w_0}$ be the corresponding seeds, then
\[
\mathcal{S}_{\w} = \Xi_{s_{\w}}.
\]
\end{theorem}

The article is structured as follows: after introducing relevant notation, we recall pseudoline arrangements and define the two collections of polyhedral objects and unimodular equivalences among them in \S\ref{sec:pa and gp}.
In \S\ref{subsec:string} we show that one of the cones is the weighted string cone and in \S\ref{subsec:super} we show how the other arises from the superpotential. 
Then in \S\ref{subsec:apply} we apply our combinatorial result and relate to toric degenerations.

\medskip

\textbf{Acknowledgements.} Both authors would like to thank X.~Fang for the stimulating and helpful discussions. In particular, we would like to thank T.~Magee for explaining his results over and over again. We are grateful to P.~Littelmann and A.~Berenstein for inspiring conversations.

\section{Preparation}\label{chap:prep}

\subsection{Notation for \texorpdfstring{$SL_n(\mathbb C)$}{}}\label{sec:pre rep theory}
We fix as Borel subgroup the upper triangular matrices $B\subset SL_n$ and diagonal matrices as maximal torus $T\subset B$. 
We denote the Borel subgroup of lower triangular matrices $B^-$ (it is also called the \emph{opposite} Borel subgroup of $B$). 
Inside of $B$ (resp. $B^-$) we have the subgroup of unipotent matrices $U$ (resp. $U^-$) with all diagonal entries being $1$. They are the \emph{unipotent radical} of $B$ (resp. $B^-$).
Consider the Lie algebra $\Lie(SL_n)=\lie{sl}_n$ and fix the Cartan decomposition $\lie{sl}_n=\lie n^-\oplus \lie h\oplus \lie n^+$.
Let us denote the \emph{root system} of $SL_n$ by $R\subset \mathbb R^{n}$. 
Denoting the standard basis of $\mathbb R^n$ by $\{\epsilon_i\}_{i=1,\dots,n}$ we fix the the simple roots of $R$ to be $\alpha_i=\epsilon_i-\epsilon_{i+1}$ for $i=1,\dots,n-1$. They generate the \emph{root lattice}. 
The positive roots are denoted $R^+=\{\beta\in R\mid \beta>0\}$. They are of form $\alpha_{i,j}:=\alpha_i+\dots+\alpha_j$ for $i\le j<n$. 
With our choice of simple roots we have $\alpha_{i,j}=\epsilon_i-\epsilon_{j+1}$. The number of positive roots is denoted by $N=\frac{n(n-1)}{2}$. 
For a positive root $\beta=\alpha_{i,j}$ let $f_{\beta}\in \lie n^-$ be the root vector of weight $-\beta$. 
For the weight lattice we choose the notation $\Lambda$ with generators the fundamental weights being $\omega_1,\dots,\omega_{n-1}$. 
Let $\Lambda^+$ denote the dominant integral weights in $\Lambda$.
By $\Lambda^{++}$ we denote the set of \emph{regular dominant weights}, i.e. those $\lambda=\sum_{i=1}^{n-1}a_i\omega_i$ with $a_i\in\mathbb Z_{>0}$.
For every $\lambda\in\Lambda^+$ there is a (finite-dimensional) irreducible representation of $\lie{sl}_n$ of highest weight $\lambda$, denote it by $V(\lambda)$. It is cyclically generated by a highest weight vector $v_\lambda\in V(\lambda)$ (unique up to scaling) over $U(\lie n^-)$, the \emph{universal enveloping algebra} of $\lie n^-$. 
The Weyl group of $SL_n$ is the symmetric group $S_n$ generated by the simple transpositions $s_i=(i,i+1)$ for $1\le i<n$. 
By $w_0$ we denote the longest element in $S_n$.
For every $w \in S_{n}$, we denote by $\ell(w)$ the minimal length of $w$ as a word in the generators $s_i$. Further, $\underline{w}$ denotes a reduced expression  
$\underline{w} = s_{i_1} \cdots s_{i_{\ell(w)}}$.
Such an expression is not unique. For any two reduced expressions of $w$ there is a sequence of local transformations leading from one to the other. These local transformations are either swapping orthogonal reflections $s_i s_j = s_j s_i$ if $|i - j | > 1$ or exchanging consecutive $s_is_{i+1}s_i = s_{i+1}s_is_{i+1}$.
The symmetric group acts on the weight lattice. % as follows. 
Fix $w\in S_n$ and $\lambda\in\Lambda^+$, then the weight space of weight $w(\lambda)$ in $V(\lambda)$, denoted $V(\lambda)_{w(\lambda)}$, is called \emph{extremal} and it is one-dimensional.

\begin{definition}\label{def: Demazure module}
For $w\in S_n$ and $\lambda\in \Lambda^+$ we fix $v_{w\lambda}\in V(\lambda)_{w\lambda}$ we consider
$U(\mathfrak b)\cdot v_{w\lambda}=:V_w(\lambda)$.
This is a $\mathfrak b$-module called the \emph{Demazure module}. 
\end{definition}
Note that though $V_w(\lambda)$ is a $\mathfrak b$-submodule of $V(\lambda)$, it is not an $\mathfrak{sl}_n$-module. 
For a Demazure module $V_w(\lambda)$ we denote by $V_w(\lambda)^{\perp}$ its orthogonal complement in $V(\lambda)^*$.

For $\lambda\in\Lambda^+$ let $L_\lambda$ be the total space of the \emph{homogeneous line bundle associated to the weight $\lambda$} over $SL_n/B$. 
These line bundles satisfy $L_{m\lambda}=L_\lambda^{\otimes m}$ for $m\ge 1$ and are ample, if $\lambda\in\Lambda^{++}$.
By the \emph{Borel-Weil-Theorem} we have
\[
H^0(SL_n/B,L_\lambda)^*\cong V(\lambda).
\]
This correspondence induces an embedding $SL_n/B\hookrightarrow \mathbb P(V(\lambda))$, by $ gB \mapsto g[v_\lambda]$.
In particular, we can realize the homogeneous coordinate ring of the flag variety as $\mathbb C[SL_n/B]=\bigoplus_{k\ge 0}V(k\lambda)^*$. Similarly, we obtain $\mathbb C[SL_n/U]=\bigoplus_{\lambda\in\Lambda^+} V(\lambda)$ which is a consequence of the \emph{Peter-Weyl-Theorem}. The quasi-affine variety $SL_n/U$ is also called \emph{base affine space}.

We consider for $w\in S_n$ the \emph{Bruhat cell} $BwB\subset SL_n$.
The quotient $BwB/B$ is called \emph{Schubert cell} and its Zariski closure $X_w:=\overline{BwB/B}\subset SL_n/B$ is the \emph{Schubert variety}.
Schubert varieties are normal, not necessarily smooth (but if singular having only rational singularities) subvarieties of the flag variety.
Their dimension equals the length of the associated Weyl group element, i.e. $\dim X_w=\ell(w)$.
The line bundles $L_\lambda$ can be restricted to Schubert varieties and the Borel-Weil Theorem generalizes as follows.
Fix $w\in S_n$ and $\lambda\in\Lambda^+$, then
\[
H^0(X_w,L_\lambda)^*\cong V_w(\lambda).
\]

\subsection{Cluster algebras}\label{sec: prep cluster}
For a reminder on cluster algebras we refer the reader to \cite{FZ02} and \cite{BFZ05}.
We review below $\mathcal A$- and $\mathcal X$-cluster varieties as introduced in \cite{FG06}, following \cite[\S2]{GHK15}.
As the cluster algebra we are interested in is skew-symmetric and of geometric type we restrict to this case.

\begin{definition}
The skew-symmetric \emph{fixed data} $\Gamma$ consists of 
\begin{itemize}
    \item a finite set $I=\{1,\dots,m\}$ of \emph{directions} with a subset of \emph{unfrozen directions} $I_{\text{uf}}=\{1,\dots,n\}$; 
    \item a lattice $N$ of rank $ |I|$;
    \item a saturated sublattice $ N_{\text{uf}} \subseteq N$ of rank $|I_{\text{uf}}|$;
    \item a skew-symmetric bilinear form $\lbrace \cdot , \cdot \rbrace: N\times N \to \mathbb Q$;
    \item a sublattice $N^{\circ}\subseteq N$ of finite index satisfying
    \[\lbrace N_{\text{uf}}, N^{\circ}\rbrace \subset \mathbb Z  \ \text{ and }\  \lbrace N,N_{\text{uf}}\cap N^{\circ}\rbrace \subset \mathbb Z;
    \]
    \item $M=\Hom(N,\mathbb Z)$ and $M^{\circ}=\Hom(N^{\circ},\mathbb Z)$.
\end{itemize}
\end{definition}

A \emph{seed} $s$ is a basis $\{e_{1,s},\dots,e_{m,s}\}$ of $N=\mathbb Z^m$.
We sometimes write $N_s$ to refer to $N$ with the associated basis.
For a fixed seed $s$ the bilinear form induces a quiver $Q_s$ given by its exchange matrix $\varepsilon_s:=(\{e_{i,s},e_{j,s}\})_{ij}$.
Let $\{f_{1,s},\dots,f_{m,s}\}$ be the dual basis for $M$.
To each seed we associate two tori $T_N\cong (\mathbb C^*)^m\cong T_M$ by
\begin{align}\label{def: seed tori}
\X_s:=T_{M,s}=\Spec(\mathbb C[N]) \text{ and } \A_s:=T_{N,s}=\Spec(\mathbb C[M]).
\end{align}

We denote the coordinates on $\X_s$ by $X_{1,s},\dots,X_{m,s}$. Corresponding to the basis of the lattice we have $X_{i,s}:=z^{e_{i,s}}$. When the seed we are working in is clear we drop it from the notation. We define \emph{mutation} at $k\in I_{\text{uf}}$ on the basis $\{e_{i,s}\}$ of the lattice $N$ for seed $s$ by
\begin{align}\label{eq: def mutation lattice basis}
e_{i,s'}:=\left\{
    \begin{matrix}
    e_{i,s}+\max\{\epsilon_{ik},0\}e_{k,s}, & \text{ for } i\not =k,\\
    -e_{k,s},& \text{ for } i=k.
    \end{matrix}
\right.
\end{align}
Then $\{e_{1,s'},\dots,e_{m,s'}\}$ forms again a basis for $N$ associated with the seed $s'=\mu_k(s)$. The dual basis for $M$ transforms as
\[
f_{i,s'}:=\left\{
    \begin{matrix}
        -f_{i,s},& \text{ for } i\not=k, \\
    f_{k,s}+\sum_{j}\max\{-\epsilon_{kj},0\}f_{j,s}, & \text{ for } i =k.
    \end{matrix}
\right.
\]
Then $\{f_{1,s'},\dots,f_{m,s'}\}$ is the dual basis for $M$ associated with $s'=\mu_k(s)$. Mutation induces birational maps between the tori
\[
\mu_k:\X_s\to \X_{\mu_k(s)} \text{ and } \mu_k:\A_s\to \A_{\mu_k(s)}.
\]
defined by the pullback of functions. We have for $\X$-tori
\begin{align}\label{eq: def pullback X-mut}
\mu_k^*(z^n):=z^n(1+z^{e_{k,s}})^{-{\{n,e_{k,s}\}_s}}, \text{ for } n\in N. 
\end{align}
For the $\A$-tori the birational map is induced by
\begin{align}\label{eq: A mutation}
\mu_k^*(A_{k,s'})= \left\{\begin{matrix}
    A_{i,s}, &\text{ for } i\not=k, \\
    \frac{\prod_{i\to k\in Q_s} A_{i,s} + \prod_{k\to j\in Q_s} A_{j,s}}{A_{k,s}}, &\text{ for } i=k.
    \end{matrix}\right.
\end{align}
Note that the variables $A_{i,s}$ for $n+1\le i\le m$ never change under mutation, we can therefore drop the $s$ from their index. To be consistent with the $\X$-notation, we set $A_{i,s}=z^{f_{i,s}}$.
We denote by $s_0$ a fixed initial seed for $\Gamma$ and for any other seed obtained from $s_0$ by a sequence of mutations we write $s\sim s_0$.
By \cite[Proposition~2.4]{GHK15} we can give the following definition (see \cite[\S3]{BFMNC} for a more detailed description).

\begin{definition}\label{def:cluster variety}
Given fixed data and a fixed initial seed $s_0$ the $\X$- (resp. $\A$-) cluster variety is defined as the scheme
\begin{align}\label{def: A and X variety}
\X := \bigcup_{s\sim s_0 }\X_{s} \ (\text{resp.  } \A:=\bigcup_{s\sim s_0} \A_{s})
\end{align}
obtained by gluing the tori $\X_{s}$ (resp. $\A_{s}$) along the birational maps induced by mutation.
\end{definition}

The relation to cluster algebras is the following. The global sections of the structure sheaf on $\A$ are related to the upper cluster algebra $\overline{\mathcal{Y}}(s_0)$ (see \cite[Definition~1.1]{BFZ05}) associated to the initial data with initial seed $s_0$:
\[
H^0(\A,\mathcal O(\A))=\overline{\mathcal{Y}}(s_0)\otimes_{\mathbb Z}\mathbb C.
\]
A natural (partial) compactification $\bar \A$ of $\A$ (an $\A$-cluster variety) is given by allowing the frozen variables $A_{n+1},\dots,A_m$ to vanish. 
We denote the resulting \emph{boundary divisor} in $\bar \A$ by
\begin{align}\label{eq: def boundary divisor}
D:=\sum_{f=n+1}^m D_f, \text{ where }D_{f}:=\{A_f=0\}\subset \bar \A.
\end{align}
Every component $D_f$ of the boundary divisor induces a (rank 1) valuation $\ord_{D_f}:\mathbb C[\A]\to\mathbb Z$ by sending a function $g\in\mathbb C[\A]$ to its order of vanishing along $D_f$. 
If $g$ has a pole along $D_f$, then $\ord_{D_f}(g)<0$ is the order of the pole.
These valuations are called \emph{divisorial discrete valuations} in \cite{GHKK14}.

A main result of \cite{GHKK14} is the definition and parametrization of the $\vartheta$-basis for $\mathbb C[\A]$. One central question is:
%The superpotential plays a crucial role in answering the question: 
\emph{When is a basis element of $\mathbb C[\A]$ also a basis element for $\mathbb C[\bar \A]$?}

The \emph{full Fock-Goncharov conjecture} (see \cite[Definition~0.6]{GHKK14}) suggests that basis elements for $\mathbb C[\A]$ are parametrized by \emph{tropical points} in $\X^{\trop}(\mathbb Z)$ (see \cite[\S2]{GHKK14}). 
We will not go into detail about this tropical space due to the following fact: 
fixing a seed $s$ we have an isomorphism 
\[
\X^{\trop}(\mathbb Z)\vert_{s}\cong N_s\cong \mathbb Z^{m}.
\]
For the purpose of this paper we always work in a fixed seed and therefore have an identification of lattice points in $N_s$ with basis elements for $\mathbb C[\A]$.
From now on we assume that the cluster variety $\A$ satisfies the full FG-conjecture, as this is the case for the cluster variety we are interested in: 
Magee showed in \cite{Mag15} that the full FG-conjecture is satisfied by the cluster variety inside $SL_n/U$.
A number of criteria for the full Fock-Goncharov conjecture to hold are discussed in \cite[\S8.4]{GHKK14} and we refer the interested reader there for more details.

Associated to each component of the boundary divisor there exists a function $\vartheta_f$ on the dual cluster variety $\X$.
Assuming the full FG-conjecture we can compute and expression for $\vartheta_f$ in $\X_{s_0}$ ($s_0$ being a fixed initial seed) as described by the Algorithm~\ref{alg:superpot via opt seeds}.

\begin{definition}\label{def:superpotential}
Let $\A$ be an $\A$-cluster variety satisfying the full Fock-Goncharov conjecture. 
Then we define the \emph{superpotential} $W:\X\to \mathbb C$ on the dual cluster variety $\X$ as
\[
W:=\sum_{f \text{ frozen vertex in }Q }\vartheta_f.
\]
\end{definition}

\begin{remark}
A seed $s_f$ for which a frozen vertex $f$ is a sink (as in the first step of Algorithm~\ref{alg:superpot via opt seeds}) is called \emph{optimized} for $f$.
Finding an optimized seed for a frozen vertex is in general a hard problem as there might be infinitely many seeds. Further, doing these computations by hand is already after a few mutation quite frustrating due to the recursive formulas. An excellent tool for such computations is provided by Keller's \emph{quiver mutation applet} \cite{QuiverApp}.
\end{remark}

\begin{algorithm}[ht]
\SetAlgorithmName{Algorithm}{} 
\KwIn{\medskip {\bf Input:\ }  A cluster variety $\A$ with initial seed $s_0$ satisfying the full FG-conjecture.}
\BlankLine
\For{every frozen vertex $f\in Q_{s_0}$}{find a sequence of mutations $\overline{\mu}$ from $s_0$ to a seed $s_f$ where $f$ is a sink.\\
    \If{$s_f=s_0$}{{\bf Output:\ } $\vartheta_f\vert_{\X_{s_0}}=z^{-e_{f,s_0}}$.}
    \Else{apply the pullback of the reverse mutation sequence to $z^{-e_{f,s_f}}$. \\
    {\bf Output:\ } $\vartheta_f\vert_{\X_{s_0}}=(\overline{\mu})^*(z^{-e_{f,s_f}})$.}}
\BlankLine
{\bf Output:\ } The superpotential $W\vert_{\X_{s_0}}=\sum_{f\text{ frozen in }Q_{s_0}} \vartheta_f\vert_{\X_{s_0}}$. 
\label{alg:superpot via opt seeds}
\caption{Computing an expression for the superpotential in a given initial seed.}
\end{algorithm}

Coming back to $\mathbb C[\A]$, note that a basis element $\vartheta\in\mathbb C[\A]$ gives an element in $\mathbb C[\bar \A]$ if $\ord_{D_f}(\vartheta)\ge 0$ for every component $D_f$ of the boundary divisor.
In particular,
\[
\vartheta \in \mathbb C[\bar \A] \ \Leftrightarrow\ \min_{f \text{ frozen}}\{\ord_{D_f}(\vartheta)\}\ge 0.
\]
We need the following notion of tropicalization for Laurent polynomials.
%Let $g=\sum_{\textbf{u}\in\mathbb Z^m} a_{\textbf u}z^{\textbf u}\in \mathbb C[z_1^{\pm1},\dots,z_m^{\pm1}]$. 
The \emph{tropicalization} of $g=\sum_{\textbf{u}\in\mathbb Z^m} a_{\textbf u}z^{\textbf u}\in \mathbb C[z_1^{\pm1},\dots,z_m^{\pm1}]$ is the function $g^{\trop}:\mathbb R^{m}\to \mathbb R$ given by 
\[
g^{\trop}(\textbf x):=\min\{\textbf x\cdot \textbf u \mid \textbf u\in \mathbb Z^{m} \text{ and } a_{\textbf u}\neq 0\}.
\]
Let $g_\vartheta\in N_s$ be the lattice point associated to $\vartheta$ for a fixed seed $s$. Then using the fact that $\vartheta_f^{\trop}(g_{\vartheta})=\ord_{D_f}(\vartheta)$, this translates to
\begin{align}\label{eq:pts in superpot cone}
\vartheta \in \mathbb C[\bar{\A}]\   \Leftrightarrow \ g_\vartheta\in \{\mathbf x\in \mathbb R^{m}\mid W\vert_{\X_s}^{\trop}(\mathbf x)\ge 0\}\cap N_s.
\end{align}
In particular, the lattice points in $\{\mathbf x\in \mathbb R^{m}\mid W\vert_{\X_s}^{\trop}(\mathbf x)\ge 0\}$ parametrize a basis for $\mathbb C[\bar \A]$.

\section{Pseudoline arrangements and Gleizer-Postnikov paths}\label{sec:pa and gp}

In the following section we develop the combinatorial background and prove the combinatorial Theorem~\ref{thm:dual cones intro}, whose main applications are Theorems~\ref{thm:intro2}\&\ref{thm:intro3} stated in the introduction. 
In particular, with Theorem~\ref{thm:dual cones intro} we Question~\ref{Q:cluster}, which is the aim of this article.

Recall our notation for the symmetric group $S_n$ from \S\ref{sec:pre rep theory}. 
We associate for each $w\in S_n$ a diagram called a pseudoline arrangement to every reduced expression $\w$. 
These diagrams turn out to be closely related to cluster algebras. In fact, to every pseudoline arrangement one can associate a quiver and then using the construction summarized in \S\ref{sec: prep cluster} define a cluster algebra. 
We start by introducing the combinatorial tools: to a pseudoline arrangement we associate two weighted cones and give a unimodular equivalence between them.

\begin{definition}\label{def:pseusoline arr}
A \emph{pseudoline arrangement} $\pa(\w)$ associated to a reduced expression $\w=s_{i_1} \cdots s_{i_{\ell(w)}}$ is a diagram consisting of $n$ horizontal \emph{pseudolines} $l_1,\dots,l_{n}$ (or short \emph{lines}) labelled at the left end from bottom to top, with crossings indicated by the reduced expression. A reflection $s_i$ indicates a crossing at \emph{level} $i$ (see e.g. Figure~\ref{fig:pa 121}).
\end{definition}
For a given reduced expression $\w=s_{i_1} \cdots s_{i_{l(w)}}$, we associate to each $s_{i_j}$ the positive root $\beta_{i_j}:=s_{i_1}\cdots s_{i_{j-1}}(\alpha_{i_j})$.
Then $\beta_{i_j}=\alpha_{k,m-1}$ for $k,m< n$ and $s_{i_j}$ induces the crossing of the lines $l_k$ and $l_m$ in $\pa(\w)$.
The crossing point is a vertex in the diagram and it is labelled ${(k,m)}$.
As two lines $l_k,l_m$ cross at most once, there is at most one position with label $(k,m)$. 
For a given $w$ the pairs appearing as labels for crossing points are exactly those for which $w(\alpha_{k,m-1})<0$.
Further, the right end of a pseudoline $l_i$ is a vertex labelled $L_i$. 
Let $\pa(\w)_0$ be the set of all vertices in $\pa(\w)$.

\begin{definition}\cite[Definition~2.2]{BFZ05}\label{def:quiver pa}
Let $w\in S_n$ with reduced expression $\w$. 
Then the quiver $Q_{\w}$ associated to $\pa(\w)$ has \emph{vertices} $w_F$ associated to faces $F$ of $\pa(\w)$ and \emph{arrows}:

(1) if two faces are at the same level separated by a crossing then there is an arrow from

left to right (see Figure~\ref{fig:arrow}a);

(2) if two faces are on consecutive levels separated by two crossings then there is an arrow 

from right to left (either upwards or downwards, see Figure~\ref{fig:arrow}b, \ref{fig:arrow}c).

\noindent
Vertices corresponding to unbounded faces are \emph{frozen} and we disregard arrows between them. All the other vertices are called \emph{mutable}.
\end{definition}

\begin{center}
\begin{figure}[ht]
\centering
\begin{tikzpicture}[scale=.7]
\draw[rounded corners] (0,.5) -- (1,.5) -- (2,1.5) -- (3,1.5);
\draw[rounded corners] (0,1.5) -- (1,1.5) -- (2,.5) -- (3,.5);
\draw[thick, ->] (.5,1) -- (2.5,1);
\node at (1.5,-.8) {$a$};

\begin{scope}[xshift=4cm]
\draw[rounded corners] (0,2) -- (1,2) -- (2,1) -- (3,1) -- (4,0) -- (5,0);
\draw[rounded corners] (0,1) -- (1,1) -- (2,2) -- (3,2);
\draw[rounded corners] (2,0) -- (3,0) -- (4,1) -- (5,1);
\draw[thick, ->] (3.5,1.5) -- (1.5,0.5);
\node at (2.5,-.8) {$b$};
\end{scope}

\begin{scope}[xshift=10cm]
\draw[rounded corners] (0,0) -- (1,0) -- (2,1) -- (3,1) -- (4,2) -- (5,2);
\draw[rounded corners] (0,1) -- (1,1) -- (2,0) -- (3,0);
\draw[rounded corners] (2,2) -- (3,2) -- (4,1) -- (5,1);
\draw[thick,->] (3.5,0.5) -- (1.5,1.5);
\node at (2.5,-.8) {$c$};
\end{scope}
\end{tikzpicture}
\caption{Arrows of the quiver arising from the pseudoline arrangement.}\label{fig:arrow}
\end{figure}
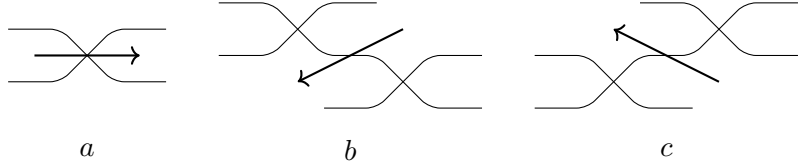
\end{center}

\begin{figure}[ht]
    \centering
\begin{tikzpicture}[scale=.8]

\node at (5.5,1.5) {\small $F_{\init_2}$};
\node at (6,.5) {\small $F_{\text{out}_2}$};
\node at (3.5,.5) {\small $F$};
\node at (1.5,1.5) {\small $F_{\text{out}_1}$};
\node at (1,.5) {\small $F_{\init_1}$};

\draw[thick, red, ->] (1.5,.5) -- (3,.5);
\draw[thick, blue, ->] (4.25,.5) -- (5.25,.5);
\draw[thick, red, ->] (4.75,1.5) -- (3.75,.75);
\draw[thick, blue, <-] (2,1.5) -- (3.25,.75);

\draw[rounded corners] (0.5,0) --(1.25,0) -- (2,0)-- (3,1) --(4,2) -- (6.5,2);
\draw[rounded corners] (.55,1) -- (2,1) -- (3,0) -- (3.5,0)-- (4,0) -- (5,1) -- (5.57,1)-- (6.5,1);
\draw[rounded corners] (.55,2) -- (3,2) -- (4,1) -- (5,0) -- (5.75,0) -- (6.5,0);

\node at (3.5,-.75) {$s_{r+1}$};
\node at (2.5,-.75) {$s_r$};
\node at (4.5,-.75) {$s_r$};

\draw[->] (7,1) -- (8,1);
\node[above] at (7.5,1) {\small $\mu_F$};
    
\begin{scope}[xshift=8cm]

\node at (6,1.5) {\small $F'_{\init_2}$};
\node at (5.5,.5) {\small $F'_{\text{out}_2}$};
\node at (3.5,1.5) {\small$F'$};
\node at (1,1.5) {\small $F'_{\text{out}_1}$};
\node at (1.5,.5) {\small $F'_{\init_1}$};

\draw[thick, red, ->] (1.5,1.5) -- (3,1.5);
\draw[thick, blue, ->] (4.25,1.5) -- (5.25,1.5);
\draw[thick, blue, ->] (3,1.25) -- (2,.5);
\draw[thick, red, ->] (5,.5) -- (4,1.25);

\draw[rounded corners] (0.55,0) --(1.25,0) -- (3,0)-- (4,1) -- (5,2) -- (6.5,2);
\draw[rounded corners] (0.55,1) -- (2,1) -- (3,2) -- (4,2) -- (5,1) -- (6.5,1);
\draw[rounded corners] (.55,2) -- (2,2) -- (3,1)-- (4,0) -- (5.75,0) -- (6.5,0);

\node at (3.5,-.75) {$s_{r}$};
\node at (2.5,-.75) {$s_{r+1}$};
\node at (4.5,-.75) {$s_{r+1}$};
\end{scope}    
\end{tikzpicture}
\caption{Mutation of pseudoline arrangements.}
\label{fig:pseudo.mut}
\end{figure}
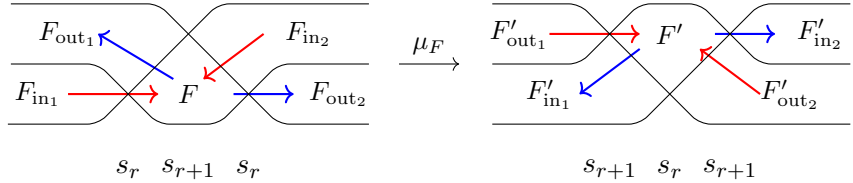

\begin{definition}\label{defn:mutation pa}
Let $w\in S_n$ with reduced expression $\w$. 
A \emph{mutation} of $\pa(\w)$ (resp. of $\w$) is a change of consecutive $s_r s_{r + 1} s_r$ in $\w$ to $s_{r+1} s_r s_{r+1}$ (or vice versa) (see Figure \ref{fig:pseudo.mut}). 
We call a face $F$ of $\pa(\w)$ \emph{mutable} if it corresponds to  $s_rs_{r+1}s_r$ (or $s_{r+1}s_r s_{r+1}$) and denote the corresponding mutation by $\mu_F$. 
The resulting pseudoline arrangement is associated to the reduced expression $\mu_F(\w)$ of $w$ and denoted by $\pa(\mu_F(\w))$.
\end{definition}

Mutation of pseudoline arrangements is a special case quiver mutation (see e.g. \cite[Definition~2.1.1]{FWZ})

Note, that the quivers $Q_{\w}$ and $Q_{\mu_{F(\w)}}$ are related by quiver mutation at the vertex $w_F$. However, $Q_{\w}$ has more mutable vertices than $\pa(\w)$ has mutable faces. 
When mutating $Q_{\w}$ at a vertex $w_{F'}$ with $F'$ not mutable in $\pa(\w)$, then for $\mu_{F'}(Q_{\w})$ there is no reduced expression of $w$ that would give rise to this quiver via a pseudoline arrangement.

Consider $\w_0\in S_n$ with reduced expression $\hat\w_0:=s_1s_2s_1s_3s_2s_1\dots s_{n-1}s_{n-2}\dots s_3s_2s_1$ and the quiver $Q_{\hat\w_0}$. 
We label the vertices for faces $F_{(i,j)}$ bounded to the left by the crossing of lines $l_i$ and $l_j$ by $w_{(i,j)}$. 
In particular, the frozen vertices at the right boundary are labelled $w_{(n-1,n)},\dots,w_{(1,n)}$ from bottom to top.
Referring to their level, the frozen vertices on the left boundary are labelled by $w_{1},\dots,w_{n-1}$ from bottom to top. 
In the following example we describe the quiver corresponding to this \emph{initial} reduced expression $\hat\w_0$ for $n=5$.

\begin{example}\label{exp:initial seed S_5}
Consider $\hat \w_0=s_1s_2s_1s_3s_2s_1s_4s_3s_2s_1\in S_5$. The pseudoline arrangement and the corresponding quiver are depicted in Figure~\ref{fig:initial}. 

\begin{center}
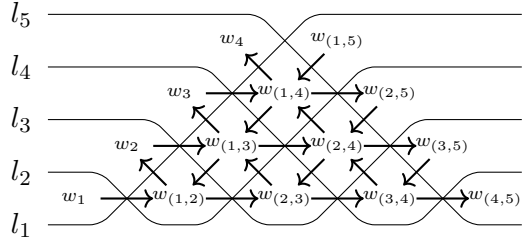
\begin{figure}[ht]
\centering
\begin{tikzpicture}[scale=.7]
%frozen
\node at (-.5,0) {$l_1$};
\node at (-.5,1) {$l_2$};
\node at (-.5,2) {$l_3$};
\node at (-.5,3) {$l_4$};
\node at (-.5,4) {$l_5$};

\node at (.5,.5) {\tiny $w_{1}$};
\node at (1.5,1.5) {\tiny $w_{2}$};
\node at (2.5,2.5) {\tiny $w_{3}$};
\node at (3.5,3.5) {\tiny $w_{4}$};
\node at (5.5,3.5) {\tiny $w_{(1,5)}$};
\node at (6.5,2.5) {\tiny $w_{(2,5)}$};
\node at (7.5,1.5) {\tiny $w_{(3,5)}$};
\node at (8.5,.5) {\tiny $w_{(4,5)}$};
%mutable
\node at (4.5,2.5) {\tiny $w_{(1,4)}$};
\node at (3.5,1.5) {\tiny $w_{(1,3)}$};
\node at (5.5,1.5) {\tiny $w_{(2,4)}$};
\node at (2.5,.5) {\tiny $w_{(1,2)}$};
\node at (4.5,.5) {\tiny $w_{(2,3)}$};
\node at (6.5,.5) {\tiny $w_{(3,4)}$};

\draw[rounded corners] (0,0) --(1,0) -- (5,4) -- (9,4);
\draw[rounded corners] (0,1) -- (1,1) -- (2,0) -- (3,0) -- (6,3) -- (9,3);
\draw[rounded corners] (0,2) -- (2,2) -- (4,0) -- (5,0) -- (7,2) -- (9,2);
\draw[rounded corners](0,3) -- (3,3) -- (6,0) -- (7,0) -- (8,1) -- (9,1);
\draw[rounded corners] (0,4) -- (4,4) -- (8,0) -- (9,0);

\draw[thick,->] (5.25,3.25) -- (4.75,2.75);
\draw[thick,->] (4.25,2.25) -- (3.75,1.75);
\draw[thick,->] (3.25,1.25) -- (2.75,.75);

\draw[thick,<-] (3.75,3.25) -- (4.25,2.75);
\draw[thick,<-] (4.75,2.25) -- (5.25,1.75);
\draw[thick,<-] (5.75,1.25) -- (6.25,.75);

\draw[thick,<-] (2.75,2.25) -- (3.25,1.75);
\draw[thick,<-] (3.75,1.25) -- (4.25,.75);

\draw[thick,<-] (1.75,1.25) -- (2.25,.75);

\draw[thick,->] (6.25,2.25) -- (5.75,1.75);
\draw[thick,->] (5.25,1.25) -- (4.75,0.75);

\draw[thick,->] (7.25,1.25) -- (6.75,0.75);

\draw[thick,->] (3,2.5) -- (4,2.5);
\draw[thick,->] (5,2.5) -- (6,2.5);

\draw[thick,->] (2,1.5) -- (3,1.5);
\draw[thick,->] (4,1.5) -- (5,1.5);
\draw[thick,->] (6,1.5) -- (7,1.5);

\draw[thick,->] (1,0.5) -- (2,0.5);
\draw[thick,->] (3,0.5) -- (4,0.5);
\draw[thick,->] (5,0.5) -- (6,0.5);
\draw[thick,->] (7,0.5) -- (8,0.5);
\end{tikzpicture}
\caption{$\pa(\hat\w_0)$ and $Q_{\hat\w_0}$ with $\hat\w_0=s_1s_2s_1s_3s_2s_1s_4s_3s_2s_1\in S_5$.}
\label{fig:initial}
\end{figure}
\end{center}
\end{example}

\subsection{Orientation and paths.}
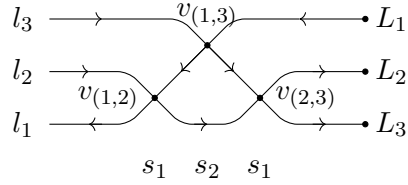
\begin{figure}[ht]
\centering
\begin{center}
\begin{tikzpicture}[scale=.7]

\node at (0,2) {$l_3$};
\node at (0,1) {$l_2$};
\node at (0,0) {$l_1$};

\node at (7,0) {$L_3$};
    \draw [fill] (6.5,0) circle [radius=0.05];
\node at (7,1) {$L_2$};
    \draw [fill] (6.5,1) circle [radius=0.05];
\node at (7,2) {$L_1$};
    \draw [fill] (6.5,2) circle [radius=0.05];

\node[below] at (2.5,-.5) {$s_{1}$};
\node[below] at (3.5,-.5) {$s_{2}$};
\node[below] at (4.5,-.5) {$s_{1}$};

\node[above] at (3.5,1.65) {$v_{(1,3)}$};
    \draw [fill] (3.5,1.5) circle [radius=0.05];
\node[right] at (4.6,0.5) {$v_{(2,3)}$};
    \draw [fill] (4.5,0.5) circle [radius=0.05];
\node[left] at (2.4,0.5) {$v_{(1,2)}$};
    \draw [fill] (2.5,0.5) circle [radius=0.05];

\draw[rounded corners] (0.5,0) --(1.25,0) -- (2,0)-- (3,1) --(4,2) -- (5.25,2);
    \draw[->, rounded corners] (1.5,0) -- (1.25,0);
    \draw[->, rounded corners] (3.5,1.5) -- (3,1);
    \draw[->, rounded corners] (6.5,2) -- (5.25,2);
\draw[rounded corners] (1.25,1) -- (2,1) -- (3,0) -- (3.5,0)-- (4,0) -- (5,1) -- (5.57,1)-- (6.5,1);
    \draw[->] (0.5,1) -- (1.25,1);
    \draw[->] (3.25,0) -- (3.5,0);
    \draw[->] (5.5,1) -- (5.75,1);
\draw[rounded corners] (1.5,2) -- (3,2) -- (4,1)-- (5,0) -- (5.75,0) -- (6.5,0);
    \draw[->] (0.5,2) -- (1.5,2);
    \draw[->] (3.5,1.5) -- (4,1);
    \draw[->] (5.5,0) -- (5.75,0);

\end{tikzpicture}
\end{center}
\caption{$\pa(\w_0)$ for $\underline w_0=s_1s_2s_1\in S_3$ with orientation for $(l_1,l_2)$.}\label{fig:pa 121} 
\end{figure}

\begin{figure}[ht]
\centering
\begin{center}
\begin{tikzpicture}[scale=.8]

\draw[rounded corners] (3,0) -- (2,0) -- (1,1) -- (0,1);
        \draw[->] (3,0) -- (2.5,0);
        \draw[->] (0.8,1) -- (.5,1);
\draw[rounded corners] (3,1) -- (2,1) -- (1,0) -- (0,0);
    \draw[->] (3,1) -- (2.5,1);
    \draw[->] (0.8,0) -- (.5,0);
\draw[->, ultra thick, red] (1.9,0.9) -- (1.1,0.1);

\begin{scope}[xshift=5cm]
  \draw[rounded corners] (3,0) -- (2,0) -- (1,1) -- (0,1);
        \draw[->] (.5,1) -- (.8,1);
        \draw[->] (2.2,0) -- (2.5,0);
\draw[rounded corners] (3,1) -- (2,1) -- (1,0) -- (0,0);
    \draw[->] (.5,0) -- (.8,0);
    \draw[->] (2.2,1)-- (2.5,1);
\draw[<-, ultra thick, red] (1.9,0.1) -- (1.1,0.9);
\end{scope}
\end{tikzpicture}
\end{center}
\caption{The two red arrows are forbidden in rigorous paths.}\label{fig:rigorous}
\end{figure}
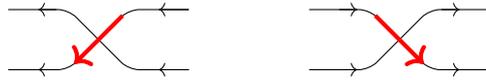

For every pair $(l_i,l_{i+1})$ with $1\le i\le n-1$ we give an orientation to a pseudoline arrangement by orienting lines $l_1,\dots,l_i$ from right to left and lines $l_{i+1},\dots,l_{n}$ from left to right, see Figure~\ref{fig:pa 121}. 
Consider an oriented path with three consecutive crossings $v_{k-1}\to v_{k}\to v_{k+1}$ belonging to the same pseudoline $l_i$. Then $v_k$ is the intersection of $l_i$ with some line $l_j$, i.e. $v_k=v_{(i,j)}$.
If either $i<j$ and both lines are oriented to the left, or $i>j$ and both lines are oriented to the right, the path is called \emph{non-rigorous}. Figure~\ref{fig:rigorous} shows these two situations. A path is called \emph{rigorous} if it is not non-rigorous. 

\begin{definition}\label{def:GPpath} 
Let $\underline{w}$ be a fixed reduced expression of $w \in S_{n}$. 
A \emph{Gleizer-Postnikov path} (or short \emph{GP-path}) is a rigorous path $\p$ in $\pa(\w)$ endowed with some orientation $(l_i,l_{i+1})$ for $i\in[n-1]:=\{1,\dots,n-1\}$. 
It has source $L_p$ and sink $L_q$ for $p\le i$ and $ q\ge i+1$. 
Further, $w(i+1)\le w(p)\le w(i)$ and $w(i+1)\le w(q)\le w(i)$. %We then say $\p$ is of \emph{shape} $(l_i,l_{i+1})$. 
The set of all GP-paths for all orientations in the pseudoline arrangement associated to $\underline{w}$ is denoted by $\mathcal P_{\w}$.
\end{definition}

Note that if $w(i)<w(i+1)$ there are no GP-paths of shape $(l_i,l_{i+1})$ and in case $w(p)\le w(q)$ there are no GP-paths with source $L_p$ and sink $L_q$.

\begin{proposition}\label{prop:in stream} 
Let $w\in S_n$ with reduced expression $\w$.
Consider $\mathbf{p}\in \mathcal P_{\w}$ of shape $(l_i,l_{i+1})$.
Then $\mathbf{p}$ is either the empty path or does not cross the lines $l_{i+1}$ and $l_i$.
In particular, $\p$ does not leave the area in $\pa(\w)$ bounded by $l_i$ and $l_{i+1}$ to the left. 
\end{proposition}

\begin{proof}
Without loss of generality we assume $w(i)<w(i+1)$, otherwise $\mathbf{p}$ is empty and we are done. 
Further, let $L_p$ be the source of $\mathbf p$ and $L_q$ the sink.
We assume $w(p)\le w(q)$, otherwise, again, $\mathbf p$ is empty. 
We focus on the part of $\pa(\w)$ to the \emph{right} of the crossing of $l_i$ and $l_{i+1}$ (which exists as $w(i)<w(i+1)$).
Observe the following:

$\bullet$ all lines crossing $l_i$ do so oriented from top to bottom.

$\bullet$ all lines crossing $l_{i+1}$ do so oriented from bottom to top.

\noindent
As $L_p$ and $L_q$ lie in between the lines $l_i$ and $l_{i+1}$ this observation implies that $\mathbf p$ can not cross $l_i$ and if it was to cross $l_{i+1}$ it could not return to $L_q$, a contradiction. 
The only possibility that is left, is if $\mathbf p$ was to follow $l_i$ through the crossing with $l_{i+1}$, but then again, it could not return to $L_q$.
\end{proof}

\subsection{Cones and polytopes arising from pseudoline arrangements}
We define two weighted cones, two cones, and two families of polytopes that arise from $\mathcal P_{\w}$ for $\w$ reduced expression of $w\in S_n$. 
We relate the two cones in the forthcoming sections, one to the weighted string cone (introduced by Littelmann \cite{Lit98} and Berenstein-Zelevinsky \cite{BZ01}), the other to the tropicalization of the (restriction of the) superpotential for a double Bruhat cell (see Magee \cite{Mag15}).

\subsection{The (weighted) GP-cone} For $\w=s_{i_1}\dots s_{i_{\ell(w)}}$ we label the standard basis of $\mathbb R^{\ell(w)}$ by crossing points in $\pa(\w)$, i.e. $\{c_{(k,m)}\mid w(\alpha_{k,m-1})<0\}$. 
Sometimes it is convenient to use the notation $c_{i_j}:=c_{(k,m)}$, when $s_{i_j}$ induces the crossing of $l_k$ and $l_m$ in $\pa(\w)$.
Consider $\mathbf{p}\in\mathcal P_{\underline{w}}$. It is uniquely determined by those vertices in $\pa(\w)_0$ where $\p$ changes from one line to another. 
For some $1\le  p \le i < q \le n$ we can therefore write $\p$ as
\[ 
\mathbf{p}=L_{p}\to v_{(p,j_1)}\to v_{(j_1,j_2)}\to\dots\to v_{(j_{k},q)}\to L_{q}.
\]
Set $j_0:=p$ and $j_{k+1}:=q$, then we associate to $\p$ the vector
\begin{align}\label{eq:def c_p}
c_{\p} := \sum_{s=0}^{k} c_{(j_s, j_{s+1})} \in \mathbb R^{\ell(w)},
\end{align}
where we set $c_{(i,j)} := - c_{(j,i)}$ if $i > j$ and $c_{(i,i)}:=0$.

\begin{definition}\label{def:gp-cone}
The following polyhedral cone is called \textit{GP-cone} (due to Gleizer-Postnikov \cite{GP00} who call it \emph{principal cone}):
\begin{align}\label{eq:def GP cone}
C_{\w} = \{ {\mathbf x} \in \mathbb{R}^{\ell(w)} \mid (c_{\p})^t({\mathbf x}) \geq 0, \forall \p \in \mathcal P_{\w} \}.
\end{align}
\end{definition}

\begin{example}
Consider the reduced expression $\w_0=s_1s_2s_1\in S_3$. 
We endow $\pa(s_1s_2s_1)$ with the orientation for $(l_1,l_2)$, i.e. $l_1$ is oriented to the left and $l_2,l_3$ are oriented to the right (see Figure~\ref{fig:pa 121}). There are two paths in $\mathcal P_{s_1s_2s_1}$ from $L_1$ to $L_2$,
\[
\mathbf p_1=L_1\to v_{(1,3)}\to v_{(1,2)}\to v_{(2,3)}\to L_2 \text{ and } \mathbf p_2=L_1\to v_{(1,3)}\to v_{(2,3)} \to L_2.
\]
They yield $c_{\mathbf p_1}= c_{(1,2)}$ and $c_{\mathbf p_2}=c_{(1,3)}-c_{(2,3)}$. Similarly for the orientation $(l_2,l_3)$ we find a path $\mathbf p_3=L_2\to v_{(2,3)}\to L_3$ with $c_{\mathbf p_3}=c_{(2,3)}$. 
Then
\[
C_{s_1s_2s_1}=\{(x_{(1,2)},x_{(1,3)},x_{(2,3)})\in \mathbb R^3\mid x_{(1,2)}\ge 0, x_{(1,3)}\ge x_{(2,3)}\ge 0\}.
\]
\end{example}

We are interested in a weighted version of this cone to relate it to string polytopes in the next section.
The weighted cone lives in $\mathbb{R}^{\ell(w)+n-1}$, where the additional basis elements are indexed $c_1, \ldots, c_{n-1}$. 
By some abuse of notation we denote by $c_{\p}$ also the vector $(c_{\p},0\dots,0)\in \mathbb R^{\ell(w)}\times \{0\}^{n-1}\subset \mathbb R^{\ell(w)+n-1}$.

For every $i\in[n-1]$ we define the following subset of $[\ell(w)]$
\begin{align}\label{eq:def J(i) and n_i}
J(i):=\{k\in[\ell(w)]\mid s_{i_k}=s_i\} \text{ with } n_i:=\#J(i).
\end{align}
Let $J(i)=\{j_1,\dots,j_{n_i}\}$, then we set $c_{[i:0]}:=c_i$ and for $1\le k\le n_i$ we define 
\begin{align}\label{eq: def wt ineq GP}
c_{[i:k]} := c_i - c_{{i_{j_k}}} - 2 \sum_{j\in J(i),j>j_k} c_{{i_j}} + \sum_{l\in J(i-1)\cup J(i+1), l>j_k} c_{{i_l}}.
\end{align}
These vectors are normal vectors to the faces of the following weighted cone.

\begin{definition}\label{def:wgp-cone}
The \emph{weighted Gleizer-Postnikov cone} $\mathcal{C}_{\w} \subset \mathbb{R}^{\ell(w)+n-1}$ is defined as
\begin{align}\label{eq:def weighted GP cone}
\mathcal{C}_{\w}  := \left\{ 
{\mathbf x} \in \mathbb{R}^{\ell(w)+n-1} \left|
\begin{matrix}
(c_\mathbf{p})^t({\mathbf x}) \ge 0 \; , \; &\forall \; \p \in \mathcal P_{\w},&\\
(c_{[i:k]})^t ({\mathbf x}) \ge 0, \; & \forall  \; i\in[n-1], 0\le k \le n_i&
\end{matrix}
\right.\right\}.
\end{align} 
\end{definition}

\begin{example}\label{exp:pathGT}
Consider $w_0\in S_n$ and the reduced expression $\hat\w_0$ defined above.
For $ i\in [n-1]$ all GP-paths in $\pa(\w_0)$ with orientation $(l_i,l_{i+1})$ are of form
\[
\p_{i,j}:=L_i\to v_{(i,n)} \to v_{(i,n-1)}\to \dots \to v_{(i,j)}\to v_{(i+1,j)}\to \dots \to v_{(i+1,n)}\to L_{i+1}.
\]
In particular, the GP-cone $C_{\w_0}$ is described by inequalities defined by the normal vectors $c_{(i,j+1)}- c_{(i+1,j+1)}$ and $c_{(i,i+1)}$ for $i\in[n-1]$ and $j \in[i+1,n-1]$.
The vectors defining weight inequalities are (for all $i < j$):
\[
c_{j-i} - c_{(i,j)} - 2 \sum_{k = 1}^{n-j} c_{(i+k, j+k)} + \sum_{k = 0}^{n-j-1} c_{(i+k, j+1+k)} + \sum_{k = 0}^{n-j} c_{(i+1+k, j+k)}.
\]
\end{example}

\subsection{The (weighted) area cone}
We associate to the set of all GP-paths $\mathcal P_{\w}$ a second cone. In this setup, the standard basis of $\mathbb{R}^{\ell(w)+n-1}$ is indexed by the faces of the pseudoline arrangement $\{e_F\mid F \text{ face of }\pa(\w)\}$. 
Namely, there are basis vectors associated to faces $F_{(i,j)}$ bounded to the left by a crossing $(i,j)$, and to faces $F_l$ unbounded to the left for every $l\in[n-1]$. 
Let $\p \in \mathcal P_{\w}$.
We denote by $\area_{\p}$ the area to the left of $\p$ (with respect to the orientation), i.e. the area enclosed by $\p$. 
Note that for non-trivial $\p$, $\area_{\p}$ is a non-empty union of faces $F$ in the pseudoline arrangement. We associate to $\p$ the vector
\begin{align}\label{eq:def area ineq}
e_{\p} := - \sum_{F \subset \area_{\p}} e_F \in \mathbb{R}^{\ell(w)+n-1}.
\end{align}
With a little abuse of notation we denote by $e_{\p}$ also the vector in $\mathbb R^{\ell(w)}$ obtained by projecting onto the first $\ell(w)$ coordinates (forgetting the coordinates belonging to the faces that are unbounded to the left, which equal $0$ in $e_{\p}$).

\begin{definition}\label{def:s-cone} For a reduced expression $\underline{w} \in S_{n}$, we define the \emph{area cone}
\begin{align}\label{eq:def area cone}
S_{\w} := \{ {\mathbf x} \in \mathbb{R}^{\ell(w)} \mid (e_{\p})^t({\mathbf x}) \geq 0, \forall\; \p \in \mathcal P_{\w} \}.
\end{align}
\end{definition}

Again, we are interested in a weighted extension of this cone.
For this, we associate to every level $ i\in[n-1]$ a union of faces. 
Consider $F_{i}$, the face of $\pa(\w)$ that is unbounded to the left at level $i$.
As before for crossings we set $F_{i_j}:=F_{(k,m)}$ if $s_{i_j}$ in $\w$ induces the crossing of $l_k$ and $l_m$ in $\pa(\w)$.
We define $\area_{i}:=F_{i}\cup \bigcup_{k=1}^{n_i}F_{i_k}$, then $\area_i \cap \area_{i'} = \varnothing$ if $i \neq i'$. 
It is called the \emph{weight area} associated to the level $i$. 
For each $k$ with $0\le k\le n_i$, we define a vector
\begin{align}\label{eq:def area wt ineq}
e_{[i:k]} := -e_{F_i} - \sum_{j\in J(i),j\le j_k}  e_{F_{i_j}} \in \mathbb{R}^{\ell(w)+n}.
\end{align}
Note that $e_{[i:0]}=-e_{F_i}$ and $e_{[i:n_i]}=-\sum_{F\subset \area_i} e_F$.

\begin{definition}\label{def:ws-cone}
The \emph{weighted area cone} $\mathcal S_{\w}\subset \mathbb R^{\ell(w)+n-1}$ associated to the reduced expression $\w$ of $w\in S_n$ is defined as 
\begin{align}\label{eq: def wt area cone}
\mathcal S_{\w} := 
\left\{  {\mathbf x} \in \mathbb{R}^{\ell(w)+n-1} \left|
\begin{matrix} (e_{\p})^t({\mathbf x}) \ge 0 \; , \; &\forall \;  \p \in \mathcal P_{\w},& \\ 
(e_{[i:k]})^t({\mathbf x}) \ge 0 \; , \; &\forall \; i\in[n-1], 0\le k \le n_i&   
\end{matrix} 
\right. \right\}. 
\end{align}
The additional inequalities induced by the $e_{[i:k]}$ are called \textit{weight inequalities}.
\end{definition}

\begin{remark}
In all four cases, $C_{\w}, \mathcal{C}_{\w}, S_{\w}$ and $\mathcal{S}_{\w}$, some of the inequalities might be redundant and these cones are far from being simplicial in general. The vectors $e_{\mathbf p},c_{\mathbf p},e_{[i:k]}$ and $c_{[i:k]}$ are normal vectors to the defining hyperplanes of the cones $S_{\w},C_{\w},\mathcal S_{\w}$ and $\mathcal C_{\w}$ respectively. Not all of them are normal vectors to facets of these cones in general.
\end{remark}

\begin{example}\label{exp:area}
Consider the reduced expression $\hat\w_0\in S_5$. We have seen all GP-paths in $\pa(\w)$ in Example~\ref{exp:pathGT}. 
Take the path ${\p}= L_1\to v_{(1,5)}\to v_{(1,4)}\to v_{(1,3)} \to v_{(1,2)}\to v_{(2,3)}\to v_{(2,4)}\to v_{(2,5)}\to L_2$. 
The area $\area_{\p}$ associated to this path is shaded blue in Figure~\ref{fig:GTstreams}. 
The weight area $\area_2$ corresponding to level $2$ is also shown in Figure~\ref{fig:GTstreams} dotted in red.
\begin{center}
\begin{figure}[ht]
\centering
\begin{tikzpicture}[scale=.75]

\draw[red] (0,1) -- (9,1);
\draw[red] (0,2) -- (9,2);

\draw[rounded corners] (0,0) --(1,0) -- (5,4) -- (9,4);
\draw[rounded corners] (0,1) -- (1,1) -- (2,0) -- (3,0) -- (6,3) -- (9,3);
\draw[rounded corners] (0,2) -- (2,2) -- (4,0) -- (5,0) -- (7,2) -- (9,2);
\draw[rounded corners] (0,3) -- (3,3) -- (6,0) -- (7,0) -- (8,1) -- (9,1);
\draw[rounded corners] (0,4) -- (4,4) -- (8,0) -- (9,0);

\draw [fill=blue, rounded corners, opacity=.2] (1.5,0.5) --  (2,0) -- (3,0)-- (6,3) -- (9,3) -- (9,4)-- (5,4) -- (1.5,0.5);
\draw [blue, opacity=.5, thick, rounded corners] (1.5,0.5) --  (2,0) -- (3,0)-- (6,3) -- (9,3) -- (9,4)-- (5,4) -- (1.5,0.5);
\draw[pattern=dots, pattern color=red, opacity=.5, thick] (0,1) rectangle (9,2);
\draw[red, opacity=.5, thick] (0,1) rectangle (9,2);

\node[blue] at (9.5,3.5) {$\area_{\p}$};
\node at (.5,.5) {$F_{1}$};
\node at (1.5,1.5) {$F_{2}$};
\node at (2.5,2.5) {$F_{3}$};
\node at (3.5,3.5) {$F_{4}$};
\node at (5.5,3.5) {$F_{(1,5)}$};
\node at (6.5,2.5) {$F_{(2,5)}$};
\node at (7.5,1.5) {$F_{(3,5)}$};
\node at (8.5,.5) {$F_{(4,5)}$};
%mutable
\node at (4.5,2.5) {$F_{(1,4)}$};
\node at (3.5,1.5) {$F_{(1,3)}$};
\node at (5.5,1.5) {$F_{(2,4)}$};
\node at (2.5,.5) {$F_{(1,2)}$};
\node at (4.5,.5) {$F_{(2,3)}$};
\node at (6.5,.5) {$F_{(3,4)}$};
\node[red] at (-.5,1.5) {$\area_2$};

\end{tikzpicture}
\caption{The area $\area_{\p}$ for ${\p}$ as in Example~\protect{\ref{exp:area}} shaded in blue and the weight area $\area_2$ dotted in red.}
\label{fig:GTstreams}
\end{figure}
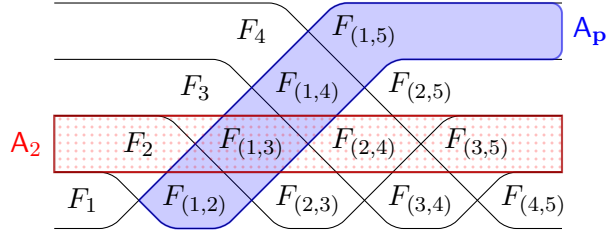
\end{center}
\end{example}

\begin{example}\label{exp:areaGT}
Consider $\hat\w_0\in S_n$ as in Example~\ref{exp:pathGT} and recall $\p_{i,j}\in\mathcal P_{\w_0}$ with $i\in[n-1]$ and $j\in[i+1,n-1]$.
The assigned area is $\area_{\p_{i,j}}=F_{(i,j)}\cup F_{(i,j+1)}\cup \dots\cup F_{(i,n)}$ for $F_{(i,k)}$ the area bounded by $v_{(i,k)}$ to the left. Hence, the cone $S_{\w_0}$ is given by inequalities defined by
\begin{align}\label{eq: e_p in GT w_0}
e_{\p_{i,j}}=-e_{F_{(i,j)}}-e_{F_{(i,j+1)}}-\dots -e_{F_{(i,n)}}-e_{F_{(i,n+1)}}.
\end{align}
The additional weight inequalities defining the cone $\mathcal S_{\w_0}$ are given by the normal vectors
\begin{align}\label{eq: f_i,k in GT w_0}
e_{[i:k]}=-e_{F_i}-e_{F_{(1,i+1)}}-e_{F_{(2,i+2)}}-\dots -e_{F_{(k,i+k)}},
\end{align}
for $i\in[n-1]$ and $0\le k\le n-i$.
\end{example}

\subsection{The polytopes} Let $\pi: \mathbb{R}^{\ell(w) + n-1} \to \mathbb{R}^{n-1}$ be the projection onto the last $n-1$ coordinates, also called \emph{weight coordinates}. 
We are interested in the preimage $\pi^{-1}(\lambda)$ for $\lambda\in\mathbb R^{n-1}$.
It is the intersection of the following hyperplanes for each $i\in[n-1]$ defined by
\begin{align}\label{eq: wt hyperplanes gp}
(c_{[i:0]})^t(\mathbf x) =\lambda_i,  \ \forall \ \mathbf x\in\mathbb R^{\ell(w)+n-1}.
\end{align}
Fix $w\in S_n$ with reduced expression $\w$. We define a second map $\tau_{\w}:\mathbb R^{\ell(w)+n-1}\to \mathbb R^{n-1}$ by
$\tau_{\w}(\mathbf x)=((e_{[i:n_i]})^t(\mathbf x))_{i=1,\dots,n-1}$.
The preimage of $\lambda\in \mathbb R^{n-1}$ with respect to $\tau_{\w}$ is also an intersection of hyperplanes in $\mathbb R^{\ell(w)+n-1}$. For each $i\in[n-1]$ they are defined by
\begin{align}\label{eq:wt hyperplanes area}
(e_{[i:n_i]})^t(\mathbf x)=\lambda_i, \  \forall \  \mathbf x\in\mathbb R^{\ell(w)+n-1}.
\end{align}

\begin{definition}\label{def: polytopes from cones}
For $w\in S_n$ with reduced expression $\w$ and for $\lambda\in\mathbb R^{n-1}$ we define the following polytopes in $\mathbb R^{\ell(w)+n-1}$
\begin{align}\label{eq:def polytopes}
\mathcal S_{\w}(\lambda):=\mathcal{S}_{\w}\cap \tau_{\w}^{-1}(\lambda)
\text{   and   } \mathcal C_{\w}(\lambda):=\mathcal C_{\w}\cap \pi^{-1}(\lambda).
\end{align}
\end{definition}

Note that by \eqref{eq:wt hyperplanes area} (resp. \eqref{eq: wt hyperplanes gp}) we obtain a description of $\mathcal S_{\w}(\lambda)$ (resp. $\mathcal C_{\w}(\lambda)$) in terms of defining equalities and inequalities
by replacing the weight inequalities $e_{[i:n_i]}^t({\mathbf x})\ge 0$ in \eqref{eq: def wt area cone} (resp. $(c_i)^t(\mathbf x)\ge 0$ in \eqref{eq:def weighted GP cone}) by $(e_{[i:n_i]})^t(\mathbf x)=\lambda_i$ (resp. $(c_i)^t(\mathbf x) =\lambda_i$).
In particular, the defining normal vectors for $\mathcal S_{\w}$ (resp. $\mathcal C_{\w}$) coincide with those for $\mathcal S_{\w}(\lambda)$ (resp. $\mathcal C_{\w}(\lambda)$).
This observation is important in the proof of Theorem~\ref{thm:unimod}.

\subsection{A unimodular equivalence} The above pairs of cones (resp. polytopes) $(S_{\w},C_{\w})$ and $(\mathcal S_{\w},\mathcal C_{\w})$ ( resp. $(\mathcal S_{\w}(\lambda),\mathcal C_{\w}(\lambda))$) have in fact more in common than the combinatorics defining them. 
To make this statement precise we need to introduce the notion of unimodular equivalence (see e.g. \cite[\S2]{HL16}).

\begin{definition}\label{def:unimod equiv}
Two polytopes $P,Q\subset \mathbb R^d$ (resp. polyhedral cones $C,D\subset \mathbb R^{d}$) are called \emph{unimodularly equivalent} if there exists matrix $M\in GL_d(\mathbb Z)$ and $w\in \mathbb Z^d$
\[
Q=f_M(P)+w \ (\text{resp. } D=f_M(C)+w),
\]
where $f_M(x)=xM$ for $x\in \mathbb R^d$. We denote this by $Q\cong P$ (resp. $C\cong D$).
\end{definition}

This notion of equivalence is of particular interest to us because of its implication on the associated toric varieties. 
Recall the construction of a projective toric variety $X_P\subset \mathbb P^{d-1}$ associated with a polytope $P\subset \mathbb R^d$ in \cite[\S2.1 and \S2.3]{CLS11}. Then $Q$ being unimodular equivalent to $P$ implies
\begin{align}\label{eq: unimod for toric}
X_Q\cong X_P.
\end{align}

We want to construct a unimodular equivalence between $\mathcal C_{\w}$ and $\mathcal S_{\w}$ for all reduced expression $\w$ of $w\in S_n$.
The following definition is the affine lattice transformation ($f_M$ in Definition~\ref{def:unimod equiv}) that defines the unimodular equivalence.
We give it in terms of the bases $\{e_F\mid F \text{ face of }\pa(\w)\}$ and $\{c_{(k,m)},c_i\mid v_{(k,m)}\in\pa(\w)_0 ,i\in[n-1]\}$.
Morally, we send a face $F$ bounded to the left by a crossing to a linear combination of its adjacent crossings (see \eqref{eq: def psi wt}).
A face unbounded to the left is sent to the sum of all crossings at its level.

\begin{definition}\label{def: psi_w}
For $w\in S_n$ and $\w$ a reduced expression we define the linear map $\Psi_{\w}:\mathbb R^{\ell(w)+n-1}\to \mathbb R^{\ell(w)+n-1}$ on the basis $\{-e_{F}\}$ associated to faces $F$ of $\pa(\w)$.
Let $F=F_{i_{j_k}}$ be the face bounded to the left by the crossing induced from $s_{i_{j_k}}=s_i$ and $J(i)=\{j_1,\dots,j_{n_i}\}$ (see \eqref{eq:def J(i) and n_i}). Then \begin{align}\label{eq: def psi}
\Psi_{\w}(-e_{F_{i_{j_k}}}):=  c_{i_{j_k}}+c_{i_{j_{k+1}}} - \sum_{\begin{smallmatrix}j\in J(i-1)\cup J(i+1),\\ j_k<j<j_{k+1}\end{smallmatrix}} c_{i_j}.
\end{align}
For every level $i\in[n-1]$, we define
\begin{align}\label{eq: def psi wt}
\Psi_{\w}(-e_{F_i}) := c_{[i:1]}.
\end{align}
\end{definition}

\begin{example}\label{exp: psi lattice}
Consider $\pa(\w)$ for $\w=s_1s_2s_1\in S_3$ as in Figure~\ref{fig:pa 121}. The two bases for $\mathbb R^{5}$ are
\[
\mathcal B_{e}=\{-e_{F_1},-e_{F_2},-e_{F_{(1,2)}},-e_{F_{(1,3)}},-e_{F_{(2,3)}}\} \text{ and } \mathcal B_c=\{c_1,c_2,c_{(1,2)},c_{(1,3)},c_{(2,3)}\}.
\]
We compute the images of elements in $\mathcal B_e$ and express them in $\mathcal B_c$. The coefficients form the columns of the following matrix with the order of the bases as given above.
\begin{align*}
\Bigg(\begin{smallmatrix}
1 & 0 & 0 & 0 & 0 \\
0 & 1 & 0 & 0 & 0 \\
-1& 0 & 1 & 0 & 0 \\
1 & -1& -1& 1 & 0 \\
-2& 1 & 1 & -1 & 1
\end{smallmatrix}\Bigg) \in GL_5(\mathbb Z).
\end{align*}
\end{example}
The observation in the example above is true in general. We obtain the following Lemma as a straightforward consequence of the definition of $\Psi_{\w}$.

\begin{lemma}\label{lem: psi lattice}
Let $w\in S_n$ with reduced expression $\w$. Order the bases induced by the faces of $\pa(\w)$ resp. by the crossing points in $\pa(\w)$ as
\[
\mathcal B_e=\{-e_{F_1},\dots,-e_{F_{n-1}},-e_{F_{i_1}},\dots-e_{F_{i_{\ell(w)}}}\},
\text{ resp. } \mathcal B_c=\{c_1,\dots,c_{n-1},c_{i_1},\dots,c_{i_{\ell(w)}}\}.
\]
Then $\Psi_{\w}$ can be represented by a lower triangular matrix ${M}_{\w}^{e,c}$ with all diagonal entries being 1. In particular, ${M}_{\w}^{e,c}\in GL_{\ell(w)+n-1}(\mathbb Z)$.
\end{lemma}

\begin{corollary}\label{cor: res psi lattice}
With assumptions as in Lemma~\ref{lem: psi lattice} consider $\Psi_{\w}\vert_{\mathbb R^{\ell(w)}}:\mathbb R^{\ell(w)}\to\mathbb R^{\ell(w)}$. We order as before the bases for $\mathbb R^{\ell(w)}$ induced by the faces resp. crossing points in $\pa(\w)$ by
\[
\overline{\mathcal B}_e=\{-e_{F_{i_1}},\dots-e_{F_{i_{\ell(w)}}}\},
\text{ resp. }  \overline{\mathcal B}_c=\{c_{i_1},\dots,c_{i_{\ell(w)}}\}.
\]
Then $\Psi_{\w}\vert_{\mathbb R^{\ell(w)}}$ can be represented by a lower triangular matrix $\overline{M}_{\w}^{e,c}$ with all diagonal entries 1. In particular, $\overline{M}_{\w}^{e,c}\in GL_{\ell(w)}(\mathbb Z)$.
\end{corollary}

\begin{remark}
The map $\Psi_{\w}$ restricted to $\mathbb R^{\ell(w)}$ is related to the Chamber Ansatz due to Berenstein-Fomin-Zelevinsky in \cite{BFZ96} (see also \cite{GKS}).
\end{remark}

\begin{figure}
\centering
\begin{tikzpicture}
 %%%% RL one F
  \node at (-.5,.5) {1a};
  \node at (3.25,1) {$l_i$};
  \node at (3.25,0) {$l_{j}$};
  \draw [fill] (1.5,.5) circle [radius=0.05];
  \node at (.9,0.5) {$c_{(i,j)}$};

\draw [fill=red, semitransparent, red, rounded corners] (3,1) --(2,1) --(1.5,.5)-- (2,0) -- (3,0) -- (3,1);
\draw[rounded corners] (0,1) -- (1,1) -- (2,0) -- (3,0);
\draw[rounded corners] (0,0) -- (1,0) -- (2,1) -- (3,1);
    \draw[->] (0.25,1) -- (0.5,1);
    \draw[->] (2.25,0) -- (2.5,0);
\draw[->] (2.5,1) -- (2.25,1);
\draw[->] (0.75,0) -- (0.5,0);

 \begin{scope}[xshift=4.5cm]
 %%%% RR one F
   \node at (-.5,.5) {2a};
     \node at (3.25,1) {$l_i$};
  \node at (3.25,0) {$l_{j}$};
    \draw [fill] (1.5,.5) circle [radius=0.05];
   \node at (.8,0.5) {$-c_{(i,j)}$};

\draw [fill=red, semitransparent, red, rounded corners] (0,1) --(1,1) --(1.5,.5)-- (2,1) -- (3,1) -- (3,1.5) -- (0,1.5) -- (0,1);
\draw[rounded corners] (0,1) -- (1,1) -- (2,0) -- (3,0);
\draw[rounded corners] (0,0) -- (1,0) -- (2,1) -- (3,1);
    \draw[->] (0.25,1) -- (0.5,1);
    \draw[->] (2.25,0) -- (2.5,0);
\draw[->] (0.25,0) -- (0.5,0);
\draw[->] (2.25,1) -- (2.5,1);

\end{scope}

%%%%% LL one F
\begin{scope}[xshift=9cm]
  \node at (-.5,.5) {3a};
    \node at (3.25,1) {$l_i$};
  \node at (3.25,0) {$l_{j}$};
    \draw [fill] (1.5,.5) circle [radius=0.05];
   \node at (.8,0.5) {$-c_{(i,j)}$};

\draw [fill=red, semitransparent, red, rounded corners] (0,0) --(1,0) --(1.5,.5)-- (2,0) -- (3,0) -- (3,-.5) -- (0,-.5) -- (0,0);

\draw[rounded corners] (0,1) -- (1,1) -- (2,0) -- (3,0);
\draw[rounded corners] (0,0) -- (1,0) -- (2,1) -- (3,1);
    \draw[->] (0.75,1) -- (0.5,1);
    \draw[->] (2.5,0) -- (2.25,0);
\draw[->] (2.5,1) -- (2.25,1);
\draw[->] (0.75,0) -- (0.5,0);

\end{scope}

%%% RL three F
\begin{scope}[yshift=-2.5cm]
  \node at (-.5,.5) {1b};
    \node at (3.25,1) {$l_i$};
  \node at (3.25,0) {$l_{j}$};
    \draw [fill] (1.5,.5) circle [radius=0.05];
  \node at (.9,0.5) {$-c_{(i,j)}$};

\draw [fill=red, semitransparent, red, rounded corners] (0,1) --(1,1) --(1.5,.5)-- (1,0) -- (0,0) -- (0,-.5) -- (3,-.5) -- (3,1.5) -- (0,1.5) -- (0,1);
\draw[rounded corners] (0,1) -- (1,1) -- (2,0) -- (3,0);
\draw[rounded corners] (0,0) -- (1,0) -- (2,1) -- (3,1);
    \draw[->] (0.25,1) -- (0.5,1);
    \draw[->] (2.25,0) -- (2.5,0);
\draw[->] (2.5,1) -- (2.25,1);
\draw[->] (0.75,0) -- (0.5,0);
\end{scope}

 %%%% LL three F
\begin{scope}[xshift=4.5cm, yshift=-2.5cm]
  \node at (-.5,.5) {2b};
    \node at (3.25,1) {$l_i$};
  \node at (3.25,0) {$l_{j}$};
    \draw [fill, rounded corners] (1.5,.5) circle [radius=0.05];
  \node at (1.5,1) {$c_{(i,j)}$};

\draw [fill=red, semitransparent, red, rounded corners] (0,1) --(1,1) --(1.5,.5)-- (2,1) -- (3,1) -- (3,-.5) -- (0,-.5) -- (0,1);
\draw[rounded corners] (0,1) -- (1,1) -- (2,0) -- (3,0);
\draw[rounded corners] (0,0) -- (1,0) -- (2,1) -- (3,1);
    \draw[->] (0.75,1) -- (0.5,1);
    \draw[->] (2.5,0) -- (2.25,0);
\draw[->] (2.5,1) -- (2.25,1);
\draw[->] (0.75,0) -- (0.5,0);
\end{scope}

%%%% RR three F
\begin{scope}[xshift=9cm, yshift=-2.5cm]
  \node at (-.5,.5) {3b};
    \node at (3.25,1) {$l_i$};
  \node at (3.25,0) {$l_{j}$};
    \draw [fill] (1.5,.5) circle [radius=0.05];
  \node at (1.5,0) {$c_{(i,j)}$};

\draw [fill=red, semitransparent, red, rounded corners] (0,0) --(1,0) --(1.5,.5)-- (2,0) -- (3,0) -- (3,1.5) -- (0,1.5) -- (0,0);
\draw[rounded corners] (0,1) -- (1,1) -- (2,0) -- (3,0);
\draw[rounded corners] (0,0) -- (1,0) -- (2,1) -- (3,1);
    \draw[->] (0.25,1) -- (0.5,1);
    \draw[->] (2.25,0) -- (2.5,0);
\draw[->] (0.25,0) -- (0.5,0);
\draw[->] (2.25,1) -- (2.5,1);
\end{scope}

\end{tikzpicture}
\caption{A path $\p$ changing the line at a crossing $(i,j)$ and the corresponding area $\area_{\p}$.}\label{fig:path.cross}
\end{figure}
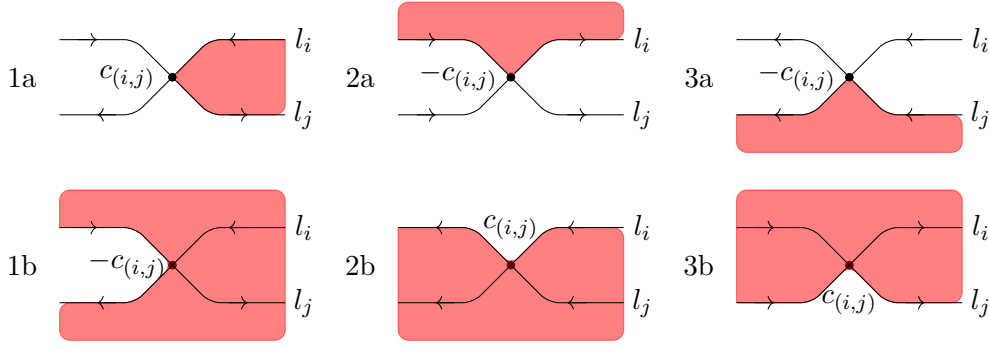

\begin{proposition}\label{prop:unimod}
Let $w\in S_n$ with reduced expression $\w$. 
For every $\p\in\mathcal P_{\w}$ we have
\[
\Psi_{\w}(e_\p)=c_\p.
\]
In particular, $\Psi_{\w}$ sends the normal vector of a defining hyperplane of $\mathcal S_{\w}$ to the normal vector of a defining hyperplane of $\mathcal C_{\w}$.
\end{proposition}

\begin{proof}
We show that for every crossing point $(i,j)$ the coefficient of $c_{(i,j)}$ coincides in $\Psi_{\w}(e_\p)$ and $c_\p$.
Recall that $\area_\p$ is the union of all faces to the left of $\p$ with respect to the given orientation. 
We distinguish three cases:

\begin{enumerate}[(i)]
\item If $(i,j)$ lies in the interior of $\area_{\p}$, then four faces $F^r\subset \area_\p$ with $r\in[4]$ are adjacent to $(i,j)$. 
For two of them in $\Psi_{\w}(-e_{F^r})$ the coefficient of $c_{(i,j)}$ is $+1$, for the other two it is $-1$.
Hence, they cancel each other and in $\Psi_{\w}(-e_\p)$ it is zero as it is in $c_{\p}$.

\item If $\p$ contains $(i,j)$ but does not change the line at $(i,j)$, then $c_{(i,j)}$ has coefficient zero in $c_{\p}$. 
For $\area_{\p}$, this means that two faces, $F^1$ and $F^2$, are adjacent to ${(i,j)}$. 
One of the two, say $F^1$, is bounded by ${(i,j)}$ to the left where for $F^2$, ${(i,j)}$ is part of the upper or lower boundary. 
In particular, $\Psi_{\w}(e_{\p})$ contains $c_{(i,j)}$ once with positive and once with negative sign, hence with the coefficient is zero.

\item Assume $\p$ changes the line at the crossing ${(i,j)}$. Figure \ref{fig:path.cross} shows the three possible orientations of $l_i$ and $l_j$.
Each yields two possibilities for the path. 
If in situation 1a, there is one face $F$ in $\area_{\p}$ bounded by ${(i,j)}$ to the left.
So $c_{(i,j)}$ has coefficient $1$ in $\Psi_{\w}(e_{\p})$. 
As $\p$ changes from $l_i$ to $l_j$ and $i<j$, also $c_{\p}$ contains $c_{(i,j)}$ with coefficient $1$.

In cases 2a and 3a, $\area_{\p}$ contains only one face bounded by ${(i,j)}$ below resp. above. 
Hence $c_{(i,j)}$ appears with coefficient $-1$ in $\Psi_{\w}(e_{\p})$. 
The same is true for $c_{\p}$: in both cases $\p$ changes from line $l_j$ to $l_i$ but $i<j$. 

Three cases remain to be checked, 1b, 2b and 3b in Figure~\ref{fig:path.cross}. 
In all of them $\area_{\p}$ contains three faces $F^1,F^2$ and $F^3$ adjacent to ${(i,j)}$. 
In case 1b, ${(i,j)}$ bounds one face to the left and the other two from above, resp. below. 
This implies that $c_{(i,j)}$ appears with coefficient $-1$ in $\Psi_{\w}(e_{\p})$. 
As $\p$ changes from $l_j$ to $l_i$ the same is true for $c_{\p}$. 
For 2b and 3b we are in the opposite case: two faces in $\area_{\p}$ are bounded to the left, resp. right, by ${(i,j)}$ and only one from above, resp. below. 
Hence, $\Psi_{\w}(e_{\p})$ contains $c_{(i,j)}$ with coefficient $1$ and the same is true for $c_{\p}$, as $\p$ changes from line $l_i$ to line $l_j$.
\end{enumerate}
\end{proof}

For our application later, it remains to show that the normal vectors defining the weight inequalities are mapped onto each other by $\Psi_{\w}$. 
Recall the weight area $\area_{i}=F_{i}\cup \bigcup_{r=1}^{n_i}F_{(i_r,j_r)}$ of level $i\in[n-1]$, with $n_i$ the number of faces bounded to the left of level $i$.

\begin{proposition}\label{prop:wunimod}
Let $w\in S_n$ with reduced expression $\w$. 
Consider $i\in [n-1]$ with $J(i)=\{j_1,\dots,j_{n_i}\}$.
Then for $ k\in [0, n_i]$ we have
\[
\Psi_{\w}(e_{[i:k]})=c_{[i:k+1]}\quad  \text{ and }\quad  \Psi_{\w}(e_{[i:n_i]})=c_{[i:0]}.
\]
In particular, $\Psi_{\w}$ sends normal vectors of defining (weight) hyperplanes of $\mathcal S_{\w}$ to normal vectors of defining (weight) hyperplanes of $\mathcal C_{\w}$.
\end{proposition}

\begin{proof}
We prove the claim by induction on $k$.
By definition we have $\Psi_{\w}(e_{[i:0]})=\Psi_{\w}(-e_{F_i}) = c_{[i:1]}$.
Let $1\le k< n_i-1$, then using induction for the third equation, we obtain
\begin{align*}
\Psi_{\w}(e_{[i:k+1]}) & \stackrel{\eqref{eq:def area wt ineq}}{=}  \Psi_{\w}(e_{[i:k]}-e_{F_{i_{j_{k+1}}}}) \\ 
& \stackrel{\eqref{eq: def psi}}{=}  c_{[i:k+1]} + c_{i_{j_{k+1}}} + c_{i_{j_{k+2}}} -\sum_{\begin{smallmatrix} j\in J(i-1)\cup J(i+1),\\ j_{k+1}<j<j_{k+2}\end{smallmatrix}} c_{i_j}\\
& \stackrel{\eqref{eq: def wt ineq GP}}{=}  c_i -c_{i_{j_{k+2}}}-2\sum_{j\in J(i),j>j_{k+2}}c_{i_j} +\sum_{\begin{smallmatrix} j\in J(i-1)\cup J(i+1),\\ j>j_{k+2}\end{smallmatrix}} c_{i_j} \ = \ c_{[i:k+2]}.
\end{align*}
Now consider $e_{[i:n_i]}=e_{[i:n_i-1]}-e_{F_{i_{j_{n_i}}}}$.  
We apply $\Psi_{\w}$ and obtain the following by induction. 
\begin{align*}
\Psi(e_{[i:n_i]})& \stackrel{\eqref{eq:def area wt ineq}}{=}  \Psi_{\w}(e_{[i:n_i-1]}-e_{F_{j_{n_i}}}) \\
& \stackrel{\eqref{eq: def psi}}{=}  c_{[i:n_i]}+ c_{i_{j_{n_i}}} - \sum_{\begin{smallmatrix} j\in J(i-1)\cup J(i+1),\\ j_{n_i}<j \end{smallmatrix}} c_{i_j} \ \stackrel{\eqref{eq: def wt ineq GP}}{=} \   c_i.
\end{align*}

\end{proof}

We can now prove the first Theorem of this section. It is a more precise formulation of Theorem~\ref{thm:dual cones intro} as stated in the introduction.

\vbox{
\begin{theorem}\label{thm:unimod}
Let $w\in S_n$ and $\w$ a reduced expression. The following polyhedral objects are unimodularly equivalent
\begin{enumerate}[(i)]
    \item $\mathcal S_{\w}\cong \mathcal C_{\w}$ via $\Psi_{\w}$,
    \item $S_{\w}\cong C_{\w}$ via $\Psi_{\w}\vert_{\mathbb R^{\ell(w)}}$,
    \item $\mathcal S_{\w}(\lambda)\cong \mathcal C_{\w}(\lambda)$ for all $\lambda\in\mathbb R^{n-1}$ via $\Psi_{\w}$.
\end{enumerate}
\end{theorem}}

\begin{proof}
We begin by proving (i). From Proposition~\ref{prop:unimod} we know that normal vectors $e_{\p}$ of $\mathcal S_{\w}$ for $\p\in\mathcal P_{\w}$ are mapped to the normal vectors $c_{\p}$ of $\mathcal C_{\w}$, i.e. $\Psi_{\w}(e_{\p})=c_{\p}$.
Further, from Proposition~\ref{prop:wunimod} we know the same is true for the normal vectors $e_{[i:k]}$ of $\mathcal S_{\w}$ for $i\in[n-1]$: 
we have $\Psi_{\w}(e_{[i:k]})=c_{[i:k+1]}$ for $k<n_i$  and  $\Psi_{\w}(e_{[i:n_i]})=c_{[i:0]}$.
As the right hand side of all defining inequalities is zero, we deduce that $\Psi_{\w}(\mathcal S_{\w})=\mathcal C_{\w}$.
By Lemma~\ref{lem: psi lattice}, $\Psi_{\w}$ is given by a matrix in $GL_{\ell(w)+n-1}(\mathbb Z)$ and hence, $\mathcal S_{\w}\cong \mathcal C_{\w}$.

To show (ii), note that by the same argument as for (i) we have $\Psi_{\w}\vert_{\mathbb R^{\ell(w)}}(S_{\w})=C_{\w}$. By Corollary~\ref{cor: res psi lattice} we deduce $S_{\w}\cong C_{\w}$.

For (iii) recall that for $\lambda\in\mathbb R^{n-1}$ the polytopes $\mathcal S_{\w}(\lambda)$ resp. $\mathcal C_{\w}(\lambda)$ are defined by the same normal vectors as $\mathcal S_{\w}$ resp. $\mathcal C_{\w}$. 
As above, it is also true that the right hand sides of the defining (in-)equalities coincides for all normal vectors being mapped onto each other. 
We therefore deduce $\mathcal{S}_{\w}(\lambda) =\mathcal{C}_{\w}(\lambda)$.
\end{proof}

%%%%%%%%%%%%%%%%%%%%%%%%%%%%%%%%%%%%%%

\section{String cones, polytopes and toric degenerations}\label{subsec:string}

Recall from \S\ref{sec:pre rep theory} our notation for $SL_n$.
In this section we recall \emph{string polytopes} and \emph{string cones} introduced by Littelmann in \cite{Lit98} and Berenstein-Zelevinsky in\cite{BZ01} as well as the \emph{weighted string cones} defined in \cite{Lit98}. 
We prove using a result from Gleizer-Postnikov in \cite{GP00} that these are exactly $\mathcal C_{\w}(\lambda)$ resp. $C_{\w}$ and $\mathcal C_{\w}$.

Littelmann \cite{Lit98} introduced in the context of quantum groups and crystal bases the so called (weighted) string cones and string polytopes $Q_{\w}(\lambda)$.
The motivation is to find monomial bases for the Demazure modules $V_{\w}(\lambda)$ for $w\in S_n$ and $\lambda\in\Lambda^+$.
Recall that $\{ f_{\alpha_{i_1}}^{m_{i_1}} \cdots f_{\alpha_{i_{\ell(w)}}}^{m_{i_{\ell(w)}}} \cdot v_\lambda \in V(\lambda) \mid m_{i_j} \geq 0 \}$
is a spanning set for $V_{\w}(\lambda)$ depending on a reduced expression $\w=s_{i_1}\cdots s_{i_{\ell(w)}}$.
Littelmann identifies a linearly independent subset of this spanning set by introducing the notion of \emph{adapted string} (see \cite[p.~4]{Lit98}) referring to a tuple $(a_1,\dots,a_{\ell(w)})\in\mathbb Z_{\ge 0}^{\ell(w)}$.
His basis for $V_{\w}(\lambda)$ consists of those elements $f_{i_1}^{a_1}\cdots f_{i_{\ell(w)}}^{a_{\ell(w)}}\cdot v_\lambda$ for which $(a_1,\dots,a_{\ell(w)})$ is adapted.

For a fixed reduced expression $\w$ of $w\in S_n$ and $\lambda\in\Lambda^+$ he gives a recursive definition of the the \emph{string polytope} $Q_{\w}(\lambda)\subset\mathbb R^{\ell(w)}$ (\cite[p.~5]{Lit98}, see also \cite{BZ01}). 
The lattice points $Q_{\w}(\lambda)\cap \mathbb Z^{\ell(w)}$ are the adapted strings for $\w$ and $\lambda$.
The \emph{string cone} $Q_{\w}\subset \mathbb R^{\ell(w)}$ is the convex hull of all $Q_{\w}(\lambda)$ for $\lambda\in\Lambda^+$. The \emph{weighted string cone} $\mathcal Q_{\w}\subset\mathbb R^{\ell(w)+n-1}$ is defined as
\[
\mathcal Q_{\w} := \conv \bigg( \bigcup_{\lambda\in\Lambda^+}Q_{\w}(\lambda)\times\{\lambda\}\bigg)\subset \mathbb R^{\ell(w)+n-1}.
\]
By definition, one obtains the string polytope from the weighted string cone by intersecting it with the hyperplanes given by $\pi^{-1}(\lambda)$ as in \eqref{eq: wt hyperplanes gp}.
The lattice points in the weighted string cone for $w=w_0$ parametrize a basis of $\mathbb{C}[SL_{n}/U]\cong \bigoplus_{\lambda\in\Lambda^+}V(\lambda)$.

\medskip

String polytopes are of great interest to us because of Caldero's work \cite{Cal02} in 2002.
He defines for a Schubert variety $X_w$ a flat family over $\mathbb A^1$ with generic fibre $X_w$ and special fibre a toric variety. 
The family is given by a construction using Rees algebras (see \cite[\S0.2]{Cal02}). Although not defined using valuations initially, it was realized this way in \cite{Kav15} and \cite{FFL15}.
His main tools are Lusztig's dual canonical basis and the string parametrization due to \cite{BZ01} and \cite{Lit98}.
We summarize his results (restricted to the case of $SL_n$) below.

Let $\w = s_{i_1} \cdots s_{i_{\ell(w)}}$ be a reduced expression of $w \in S_{n}$. 
We extend $\w$ to the \emph{right} to a reduced expression $\w_0 =\w s_{i_{\ell(w)+1}} \cdots s_{i_N}$ of $w_0$. 
This extension is not unique but the results are independent of the extension.
Caldero realizes the string cone $Q_{\w}$ for the Demazure module $V_w(\lambda)$ as a face of the string cone $Q_{\w_0}$. He deduces the following Lemma as a consequence of \cite[\S1]{Lit98}.

\begin{lemma*}(see \cite{Lit98},\cite[Lemma~3.3]{Cal02})\label{lem: restr. string cone}
Let $w\in S_n$ with $\w=s_{i_1}\cdots s_{i_{\ell(w)}}$ a reduced expression and choose a reduced expression $\w_0=\w s_{i_{\ell(w)+1}}\cdots s_{i_N}$. Then the weighted string cone $\mathcal Q_{\w}$ is obtained from the weighted string cone $\mathcal Q_{\w_0}$ by setting the variables corresponding to $s_{i_{\ell(w)+1}},\dots,s_{i_N}$ equal to zero.
\end{lemma*}

Caldero defines a filtration $(\mathcal F_{\le m})_{m\ge 1}$ on $\mathbb C[SL_n/U]$ with associated graded algebra the semi-group algebra $\mathbb C[\mathcal Q_{\w_0}\cap\mathbb Z^{N+n-1}]$. Using the Lemma, he defines a quotient filtration $(\overline{\mathcal F}_{\le m})_{m\ge 1}$ on $\mathbb C[SL_n/U]/I_{w}$, where $I_w=\bigoplus_{\lambda\in \Lambda^+} V_w(\lambda)^\perp$ (recall \S\ref{sec:pre rep theory}), i.e. 
\[
\mathbb C[SL_n/U]/I_{w}=\bigoplus_{\lambda\in\Lambda^+}V(\lambda)^*/\bigoplus_{\lambda\in\Lambda^+}V_w(\lambda)^{\perp}.
\]
The semi-group algebra $\mathbb C[\mathcal Q_{\w}\cap \mathbb Z^{\ell(\w)+n-1}]$ is the associated graded algebra of the quotient filtration.
In particular, he degenerates $X_w$ into a toric variety $Y$, whose normalization is the toric variety $X_{\mathcal Q_{\w}(\lambda)}$ associated to the string polytope $\mathcal Q_{\w}(\lambda)$ for $\lambda\in\Lambda^{++}$.

\subsection{Relation to the GP cones} Gleizer and Postnikov developed in \cite{GP00} a combinatorial model to describe string cones $Q_{\w_0}$ \emph{non-recursively} for every reduced expression $\w_0$ of $w_0\in S_n$. They use pseudoline arrangements and GP-paths to obtain the following.

\begin{corollary*}\cite[Corollary~5.8]{GP00}
Let $\w_0$ be reduced expression of $w_0\in S_n$. Then $C_{\w_0}=Q_{\w_0}$.
\end{corollary*}

On our way to showing that a toric variety isomorphic to $X_{\mathcal Q_{\w}(\lambda)}$ arises in the context of cluster varieties and mirror symmetry, we first generalize Gleizer-Postnikov's result as follows.

\vbox{
\begin{theorem}\label{thm:wt GP is wt string}
For every $w \in S_{n}$ with reduced expression $\w$ and every extension $\w_0=\w s_{i_{\ell(w)}+1}\cdots s_{i_N}$ the following polyhedral objects coincide
\begin{enumerate}[(i)]
    \item $\mathcal C_{\w}=\mathcal Q_{\w}$,
    \item $C_{\w}=Q_{\w}$,
    \item $\mathcal{C}_{\w}(\lambda)=\mathcal Q_{\w}(\lambda)$ for $\lambda\in\mathbb R^{n-1}$.
\end{enumerate}
\end{theorem}}

In order to prove Theorem~\ref{thm:wt GP is wt string} we show how to obtain $\mathcal C_{\w}$ from restricting $\mathcal C_{\w_0}$ for appropriate $\w_0$. The next subsection is dedicated to introducing restricted paths and concludes with the proof of Theorem~\ref{thm:wt GP is wt string}.

\subsection{Restriction of paths} 
We show that for $\w_0=\w s_{i_{\ell(w)+1}}\cdots s_{i_N}$ we obtain $\mathcal C_{\w}$ from $\mathcal C_{\w_0}$ by setting to zero the coordinates corresponding to crossing points $c_{i_{\ell(w)+1}},\dots,c_{i_N}$ in $\pa(\w_0)$.

\begin{definition}\label{def: res path}
Let $\w$ be a reduced expression of $w\in S_n$ and fix $\w_0=\w s_{i_{\ell(w)+1}}\cdots s_{i_N}$. 
Consider $\mathbf{p}_{\w_0}\in\mathcal P_{\w_0}$ and draw it in $\pa(\w_0)$.
Then cut $\pa(\w_0)$ in two pieces along a vertical line, such that all crossing points $v_{i_p}$ corresponding to $s_{i_p}$ with $1\le p\le \ell(w)$ are on the left of the cut and all $v_{i_q}$ corresponding to $s_{i_q}, \ell(w)<q\le N$ are on the right (see Figure~\ref{fig: res path}).
We define the \emph{restriction} $\res_{\w}(\mathbf{p}_{\w_0})$ of ${\mathbf p}_{\w_0}$ to $\pa(\w)$ as the part of $\p_{\w_0}$ to the left of the cut.
\end{definition}

We label the intersection points of the lines $l_i$ with the cutting line by $\hat L_i$.
An alternative way of describing $\res_{\w}(\p_{\w_0})$ is by removing all vertices $v_{(i,j)}$ from it for which $w(\alpha_{i,j-1})>0$. 
Denote by $\res_{\w}(\mathcal P_{\w_0})$ the set of all paths in $\pa(\w)$ that appear in a restriction of a path in $\mathcal P_{\w_0}$ (counting each path only once).

\begin{center}
\begin{figure}[ht]
\centering
\begin{tikzpicture}[scale=.6]
\draw[rounded corners] (0,0) -- (2,0) -- (3,1) -- (5,1) -- (7,3) -- (9,3);
\draw[rounded corners] (0,1) -- (1,1) -- (2,2) -- (3,2) -- (4,3) -- (6,3) -- (7,2) -- (9,2);
\draw[rounded corners] (0,2) -- (1,2) -- (3,0) -- (7,0) -- (8,1) -- (9,1);
\draw[rounded corners] (0,3) -- (3,3) -- (4,2) -- (5,2) -- (6,1) -- (7,1) -- (8,0) -- (9,0);

\draw[->] (9,3) --  (8,3);
\draw[<-] (8.5,0) -- (8.1,0);
\draw[<-] (2,3) -- (1.5,3);
\draw[<-] (1,0) -- (1.9,0);
\draw[->] (0,2) -- (.5,2);
\draw[->] (.8,1) -- (.5,1);

\node at (1.5,-.5) {$s_2$};
\node at (2.5,-.5) {$s_1$};
\node at (3.5,-.5) {$s_3$};
\node at (5.5,-.5) {$s_2$};
\node at (6.5,-.5) {$s_3$};
\node at (7.5,-.5) {$s_1$};

\node at (-.5,0) {$l_1$};
\node at (-.5,1) {$l_2$};
\node at (-.5,2) {$l_3$};
\node at (-.5,3) {$l_4$};

%% p_\w c_2
\draw[rounded corners, blue, ultra thick] (4.5,3) -- (4,3) -- (3.5,2.5);
    \draw[rounded corners, blue, ultra thick] (3.5,2.5) -- (4,2) -- (4.5,2);
   % \node[blue] at (4,2.5) {$\underline{c}_2$};
    \draw[blue, ultra thick, ->] (4.5,3) -- (4.1,3);
    \draw[blue, ultra thick, ->] (4.1,2) -- (4.2,2);
%% p'_\w   c_4 
\draw[rounded corners, blue, ultra thick] (4.5,1) -- (3,1) -- (2.5,.5);
    \draw[rounded corners, blue, ultra thick] (2.5,0.5) -- (3,0) -- (4.5,0);
    %\node[blue] at (3.5,0.5) {$\underline{c}_4$};
    \draw[blue, ultra thick, ->] (3.8,1) -- (3.5,1);
    \draw[blue, ultra thick, ->] (3.1,0) -- (3.5,0);
%% p_\w_0 c_1 and c_5    
\draw[rounded corners, ultra thick, red] (9,2) -- (7,2) -- (6,3) -- (4.5,3);
    %\node[teal] at (7,2.5) {$\underline{c}_1$};
    \draw[rounded corners, ultra thick, red] (9,1) -- (8,1) -- (7,0) -- (4.5,0);
    % \node[teal] at (7,.5) {$\underline{c}_5$};
    \draw[red, ultra thick, ->] (8,2) -- (7.8,2);
    \draw[red, ultra thick, ->] (6,0) -- (6.5,0);
%% p_\w_0 c_3    
\draw[rounded corners, ultra thick, red]    (4.5,2) -- (5,2) -- (5.5,1.5);
    \draw[rounded corners, ultra thick, red] (4.5,1) -- (5,1) -- (5.5,1.5);
    % \node[red] at (5,1.5) {$\underline{c}_{3}$};
    \draw[red, ultra thick,->] (5.1,1.9) -- (5.3,1.7);
    \draw[red, ultra thick, ->] (4.8,1) -- (4.6,1);
    
\draw[dashed, thick] (4.5,-.8) -- (4.5,3.8);    
\end{tikzpicture}
\caption{A path ${\p}_{\w_0}\in\mathcal P_{\w_0}$ for $\w_0=s_2s_1s_3s_2s_3s_1$ that restricts to two paths ${\p}_{\w},{\p}'_{\w}\in\mathcal{P}_{\w}$ (in blue to the left of the dashed cut) for $\w=s_2s_1s_3$. %The labels $\underline{c}_i$ denote the vectors corresponding to each fragment of ${\p}_{\w_0}$.
}\label{fig: res path}
\end{figure}
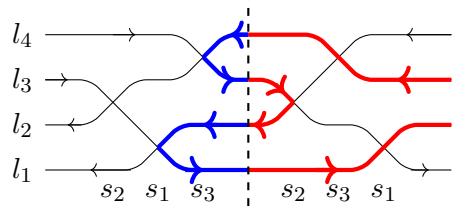
\end{center}

\begin{example}\label{exp: res path}
Consider $\w=s_2s_1s_3$ and extend it to $\w_0=s_2s_1s_3s_2s_3s_1\in S_4$. We draw $\pa(\w)$ and endow it with the orientation for $(l_2,l_3)$. 
Figure~\ref{fig: res path} shows a GP-path ${\p}_{\w_0}$.
Its restriction $\res_{\w}({\p}_{\w_0})$ consists of two GP-paths for $\w$ shown in blue to the left of the cut.  
\end{example}

\begin{proposition}\label{prop:respath}
Let $\w$ be a reduced expression of $w\in S_n$ and fix $\w_0=\w s_{i_{\ell(w)+1}}\cdots s_{i_N}$. 
Consider $\p_{\w_0}\in\mathcal P_{\w_0}$, then $\res_{\w}(\p_{\w_0})$ is either empty or a union of paths in $\mathcal P_{\w}$.
In particular, $\res_{\w}(\mathcal P_{\w_0})\subset \mathcal P_{\w}$.
\end{proposition}

\begin{proof}
Let $\p_{\w_0}$ be a path for orientation $(l_i,l_{i+1})$, i.e. of form $\p_{\w_0}=L_{i}\to v_{(i,j_1)}\to v_{(j_1,j_2)}\to\dots\to v_{(j_k,i+1)}\to L_{i+1}$.
To simplify notation we set $i=j_0$ and $i+1=j_{k+1}$.
First note that if $w(\alpha_{j_r,j_{r+1}-1})>0$ for all $0\le r\le k$ then $\res_{\w}(\p_{\w_0})=\varnothing$.
Otherwise $\res_{\w}(\p_{\w_0})$ is a union of paths
\[
\p_r=\hat L_{j_r}\to v_{(j_r,j_{r+1})}\to\dots\to v_{(j_{r+s},j_{r+s+1})}\to \hat L_{j_{r+s+1}},
\]
such that $w(\alpha_{j_{r+p},j_{r+p+1}-1})<0$ for all $0\le p\le s$, $0\le r\le k$ and $ 0\le s\le k-r$.
By definition, each $\p_r$ is rigorous and hence, in $\mathcal P_{\w}$.
\end{proof}

\begin{algorithm}
\SetAlgorithmName{Algorithm}{} 
\KwIn{\medskip {\bf Input:\ }  A path in $\mathcal P_{\w} \ni {\p}_{\w}=\hat{L}_{i-l}\to v_{r_1}\to \dots \to v_{r_m}\to \hat{L}_{i+m}$ for orientation $(l_i,l_{i+1})$.}
\BlankLine
{\bf Initialization:} extend $\w$ to $\w_0=\w s_{i_{\ell(w)+1}}\cdots s_{i_N}$;\\ complete $\pa(\w)$ to $\pa(\w_0)$ with orientation for $(l_i,l_{i+1})$;\\
set $p=q=0$ and $\hat\p_{\w}=\p_{\w}$.\\
\For{$p<m-1$}{follow $l_{i+m-p}$ with respect to the orientation to the next crossing with a line $l_{i+m-p-p'}$ with $p'\in[m-p-1]$,\\
    \If{$p'=m-p-1$}{
        {\bf Output:} $\hat\p_{\w}\to v_{(i+m-p,i+1)}\to L_{i+1}$.
        }
    \Else{replace $p$ by $p+p'$ and $\hat\p_{\w}$ by $\hat\p_{\w}\to v_{(i+m-p,i+m-p-p')}$ and start over.}
}
\For{$q<l$}{
follow $l_{i-l+q}$ against the orientation to the next crossing with a line $l_{i-l+q+q'}$ with $q'\in[l-q]$,\\
    \If{$q'=l-q$}{
        {\bf Output:} $L_i\to v_{(i,i-l+q+q')}\to\hat\p_{\w}$.
        }
    \Else{replace $q$ by $q+q'$ and $\hat\p_{\w}$ by $v_{(i-l+q+q',i-l+q)}\to\hat\p_{\w}$ and start over.}
}
\BlankLine
{\bf Output:\ } A path $\ind_{\w_0}(\p_{\w}):= L_{i}\to v_{(i,i-l+q)}\to \dots\to\p_{\w}\to \dots\to v_{(i+m-p,i+1)}\to L_{i+1}$.
\label{alg:induced path}
\caption{Constructing the induced path $\ind_{\w_0}(\p_{\w})$ from $\p_{\w}\in\mathcal P_{ \w}$.}
\end{algorithm}

We want to show the other implication, $\mathcal P_{\w}\subset \res_{\w}(\mathcal P_{\w_0})$. 
In Algorithm~\ref{alg:induced path} we give a construction to obtain a path in $\pa(\w_0)$ from a given path in $\mathcal P_{\w}$.
The following proposition shows that the algorithm always terminates and that the output is in fact a path in $\mathcal P_{\w_0}$.

\begin{proposition}\label{prop: algo induced path}
Algorithm~\ref{alg:induced path} terminates for all $\p_{\w}\in\mathcal P_{\w}$ and $\ind_{\w_0}(\p_{\w})\in\mathcal P_{\w_0}$.
\end{proposition}
\begin{proof}
By Proposition~\ref{prop:in stream} $\p_{\w}$ lies in the region of $\pa(\w)$ in between the lines $l_i$ and $l_{i+1}$. In particular, at some point there is a $p'$ with $l_{i+m-p-p'}=l_{i+1}$ terminating the first loop and a $q'$ with $l_{i-l+q+q'}=l_i$ terminating the second loop.

To see that $\ind_{\w_0}(\p_{\w})\in\mathcal P_{\w_0}$ observe that changing the lines as indicated by the algorithm avoids exactly the two situations from Figure~\ref{fig:rigorous} forbidden in rigorous paths. Hence, $\ind_{\w_0}(\p_{\w})$ is rigorous. 
\end{proof}

By Proposition~\ref{prop: algo induced path} we can define the following.

\begin{definition}\label{def: ind path}
Let $\w$ be a reduced expression of $w\in S_n$ and fix $\w_0=\w s_{i_{\ell(w)+1}}\cdots s_{i_N}$. 
For $\p_{\w}\in\mathcal P_{\w}$ we define the \emph{induced path} $\ind_{\w_0}(\p_{\w})\in\mathcal P_{\w_0}$ as the output of Algorithm~\ref{alg:induced path}.
\end{definition}

\begin{example}\label{exp: ind path}
Consider $\w=s_2s_1s_3$ and extend it to $\w_0=s_2s_1s_3s_2s_3s_1\in S_4$. We draw $\pa(\w)$ and endow it with the orientation for $(l_2,l_3)$. 
Figure~\ref{fig: ind path} shows a GP-path ${\p}_{\w}$ in blue to the left of the cut. The extension of $\p_{\w}$ in red to the right of the cut completes $\p_{\w}$ to the induced path $\ind_{\w_0}(\p_{\w})$ that is the output of Algorithm~\ref{alg:induced path}.
\end{example}

\begin{center}
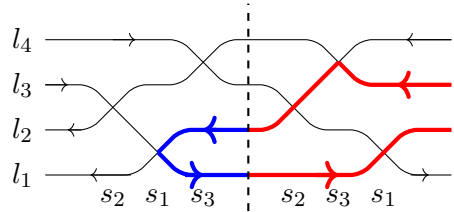
\begin{figure}
\centering
\begin{tikzpicture}[scale=.6]
\draw[rounded corners] (0,0) -- (2,0) -- (3,1) -- (5,1) -- (7,3) -- (9,3);
\draw[rounded corners] (0,1) -- (1,1) -- (2,2) -- (3,2) -- (4,3) -- (6,3) -- (7,2) -- (9,2);
\draw[rounded corners] (0,2) -- (1,2) -- (3,0) -- (7,0) -- (8,1) -- (9,1);
\draw[rounded corners] (0,3) -- (3,3) -- (4,2) -- (5,2) -- (6,1) -- (7,1) -- (8,0) -- (9,0);

\draw[->] (9,3) --  (8,3);
\draw[<-] (8.5,0) -- (8.1,0);
\draw[<-] (2,3) -- (1.5,3);
\draw[<-] (1,0) -- (1.9,0);
\draw[->] (0,2) -- (.5,2);
\draw[->] (.8,1) -- (.5,1);

\node at (1.5,-.5) {$s_2$};
\node at (2.5,-.5) {$s_1$};
\node at (3.5,-.5) {$s_3$};
\node at (5.5,-.5) {$s_2$};
\node at (6.5,-.5) {$s_3$};
\node at (7.5,-.5) {$s_1$};

\node at (-.5,0) {$l_1$};
\node at (-.5,1) {$l_2$};
\node at (-.5,2) {$l_3$};
\node at (-.5,3) {$l_4$};

%% p'_\w   c_4 
\draw[rounded corners, blue, ultra thick] (4.5,1) -- (3,1) -- (2.5,.5);
    \draw[rounded corners, blue, ultra thick] (2.5,0.5) -- (3,0) -- (4.5,0);
    \draw[blue, ultra thick, ->] (3.8,1) -- (3.5,1);
    \draw[blue, ultra thick, ->] (3.1,0) -- (3.5,0);
    \draw[rounded corners, ultra thick, red] (9,1) -- (8,1) -- (7,0) -- (4.5,0);   \draw[red, ultra thick, ->] (6,0) -- (6.5,0); 
    \draw[red, ultra thick, rounded corners] (4.5,1) -- (5,1)-- (6.5,2.5);
    \draw[red, ultra thick, rounded corners] (6.5,2.5) -- (7,2) -- (9,2);
    \draw[red, ultra thick, ->] (8,2) -- (7.8,2);
    
\draw[dashed, thick] (4.5,-.8) -- (4.5,3.8);    
\end{tikzpicture}
\caption{A path ${\p}_{\w}\in\mathcal P_{\w}$ for $\w=s_2s_1s_3$ and the induced path $\ind_{\w_0}(\p_{\w})\in\mathcal P_{\w_0}$ with $\w_0=\w s_2s_3s_1$.
}\label{fig: ind path}
\end{figure}
\end{center}

\begin{proposition}\label{prop:restrict}
Let $\w$ be a reduced expression of $w\in S_n$ and fix $\w_0=\w s_{i_{\ell(w)}+1}\cdots s_{i_N}$.
For every $\mathbf{p}_{\w}\in\mathcal P_{\w}$ there exists $\mathbf{p}_{\w_0}\in\mathcal P_{\w_0}$ such that, $\res_{\w}(\mathbf{p}_{\underline w_0})=\mathbf{p}_{\w}$.  
In particular, we have $\mathcal P_{\w}\subset\res_{\w}(\mathcal P_{\w_0})$.
\end{proposition}

\begin{proof}
By construction $\ind_{\w_0}({\p}_{\w})\in\mathcal P_{\w_0}$ satisfies $\res_{\w}(\ind_{\w_0}({\p}_{\w}))={\p}_{\w}$.
\end{proof}

Recall for $i\in[n-1]$ the definition $J(i)$ and $n_i$ from \eqref{eq:def J(i) and n_i}. To distinguish between the sets for $\w$ and $\w_0$, we use the notation $J(i)^{\w}$ (resp. $J(i)^{\w_0}$) and $n_i^{\w}$ (resp. $n_i^{\w_0}$). We define the following polyhedral objects from restricted paths and show they equal the (weighted) GP-cone, respectively polytope, in the subsequent key proposition for proving Theorem~\ref{thm:wt GP is wt string}. 

\begin{definition}\label{def: res GP cone}
Let $\w$ be a reduced expression of $w\in S_n$ and fix $\w_0=\w s_{i_{\ell(w)}+1}\cdots s_{i_N}$.
We define the  \emph{restricted weighted GP-cone} as
\begin{align}\label{eq:def res wt GPcone}
\res_{\w}(\mathcal C_{\w_0}):=\left\{ 
\mathbf x\in \mathbb R^{\ell(w)+n-1} \left\vert
\begin{matrix}
(c_{\res_{\w}(\p_{\w_0})})^t(\mathbf x)\ge 0, & \forall \ \p_{\w_0}\in\mathcal P_{\w_0},\\
(c_{[i:k]})^t(\mathbf x)\ge 0, & \forall i\in[n-1], 0\le k\le n_i^{\w}
\end{matrix}
\right.\right\}.
\end{align}
Similarly, we define $\res_{\w}(C_{\w_0}):=\{ \mathbf x\in \mathbb R^{\ell(w)}\mid  
(c_{\res_{\w}(\p_{\w_0})})^t(\mathbf x)\ge 0,  \forall \ \p_{\w_0}\in\mathcal P_{\w_0}\}$ the \emph{restricted GP-cone} and the polytope $\res_{\w}(\mathcal C_{\w_0}(\lambda)):=\res_{\w}(\mathcal C_{\w_0})\cap \pi^{-1}(\lambda)$ (see \eqref{eq: wt hyperplanes gp}) for $\lambda\in\mathbb R^{n-1}$.
\end{definition}

\vbox{
\begin{proposition}\label{prop:wt GP res}
For every $w \in S_{n}$ with reduced expression $\w$ and every extension $\w_0=\w s_{i_{\ell(w)}+1}\cdots s_{i_N}$ the following polyhedral objects coincide
\begin{enumerate}[(i)]
    \item $\mathcal C_{\w}=\res_{\w}(\mathcal C_{\w_0})$,
    \item $C_{\w}=\res_{\w}(C_{\w_0})$,
    \item $\mathcal{C}_{\w}(\lambda)=\res_{\w}(\mathcal C_{\w_0}(\lambda))$ for $\lambda\in\mathbb R^{n-1}$.
\end{enumerate}
\end{proposition}}

\begin{proof}
We start by showing (i), then (ii) and (iii) are direct implications.
Note that only the inequalities induced by GP-paths differ in the definition of $\mathcal C_{\w}$ \eqref{eq:def weighted GP cone}, resp. $\res_{\w}(\mathcal C_{\w_0})$ \eqref{eq:def res wt GPcone}.
By Proposition~\ref{prop:respath} we have $\mathcal C_{\w}\subseteq \res_{\w}(\mathcal C_{\w_0})$. By Proposition~\ref{prop:restrict} we deduce $\res_{\w}(\mathcal C_{\w_0})\subseteq\mathcal C_{\w}$ and hence, equality follows.
\end{proof}

We have now collected all ingredients necessary to provide the proof of Theorem~\ref{thm:wt GP is wt string}.

\begin{proof}[Proof of Theorem~\ref{thm:wt GP is wt string}]
We show $\mathcal Q_{\w}=\res_{\w}(\mathcal C_{\w_0})$ for every extension $\w_0=\w s_{i_{\ell(w)+1}}\cdots s_{i_N}$ and then apply Proposition~\ref{prop:wt GP res}.
By \cite[Lemma~3.3]{Cal02} (restated above) we know that
\[
\mathcal Q_{\w}=\mathcal Q_{{\w}_0}\cap\bigcap_{(i,k):\ w(\alpha_{i,k-1})>0}\{x_{(i,k)}=0\},
\]
as the $x_{(i,k)}$ appearing in the intersection of hyperplanes on the right correspond to the coordinates $x_{s_p}$ with $\ell(w)<p\le N$ in the extension of $\w$ to $\w_0$.
Further, we observe that if $c_{\p_{\w_0}}=\sum_{k}c_{(i_k,j_k)}$ then $c_{\res_{\w}(\p_{\w_0})}=\sum_{k: w(\alpha_{i_k,j_k-1})>0} c_{(i_k,j_k)}$.
Regarding the  normal vectors for weight inequalities $c_{[i:k]}$ (see \eqref{eq: def wt ineq GP}), observe that for $k>n_i^{\w}$ we obtain $c_i$ from $c_{[i:k]}$ when setting those $c_{(i_k,j_k)}$ to zero with $w(\alpha_{i_k,j_k-1})>0$.
Hence, $\mathcal Q_{\w}=\res_{\w}(\mathcal C_{\w_0})=\mathcal C_{\w}$ by Propositon~\ref{prop:wt GP res}.
Then $Q_{\w}=C_{\w}$ is a direct consequence and identifying $\Lambda^+$ with $\mathbb R^{n-1}$ using the fundamental weights, (iii) follows.
\end{proof}

\section{Double Bruhat cells and the superpotential}\label{subsec:super}

Recall the notation for cluster varieties from above.
In this section we explain the $\A$-cluster variety that can be associated to the quiver from a pseudoline arrangement.
Based on results of Berenstein-Fomin-Zelevinsky this variety is a double Bruhat cell (see Definition~\ref{def: double bruhat cell}).
We apply the construction of \cite{GHKK14} (see \S\ref{sec: prep cluster}) and recall results of Magee in \cite{Mag15} regarding the superpotential.

Recall that $SL_n$ has two cell decompositions (\emph{the Bruhat decompositions}) in terms of Bruhat cells indexed by elements of the symmetric group
\[
SL_n=\bigcup_{u\in S_n}BuB=\bigcup_{v\in S_n} B^-vB^-.
\]

\begin{definition}\label{def: double bruhat cell}
The \emph{double Bruhat cell} associated to $e$ and $w$ in $S_n$ is
\[
G^{e,w}:=B\cap B^-wB^-\subset SL_n.
\]
\end{definition}

The cluster structure of $G^{e,w}$ can be established as follows. Choose a reduced expression $\w$ and consider $\pa(\w)$.
Recall from Definition~\ref{def:quiver pa} that every face of $\pa(\w)$ corresponds to a vertex of $Q_{\w}$.
We therefore associate cluster variables to faces of $\pa(\w)$.
Let $F$ be such a face and assume the lines $l_{j_1},\dots,l_{j_k}$ pass below $F$.
In particular, $F$ is of level $k$.
Then associate the Pl\"ucker coordinate $\bar p_{j_1,\dots,j_k}\in\mathbb C[SL_n]$ to $F$, i.e. the minor of the columns $[k]$ and rows $\{{j_1},\dots,{j_k}\}$.
To remember it was associated with $F$, we set $A_F:=\bar p_{j_1,\dots,j_k}$.

\begin{definition}
Let $w\in S_n$ with reduced expression $\w$. Then the quiver $Q_{\w}$ together with the set of cluster variables $\mathbf A_{\w}:=\{A_F\mid F \text{ a face of }\pa(\w)\}$ forms the seed $s_{\w}:=(\mathbf A_{\w},Q_{\w})$.
\end{definition}

\begin{example}\label{exp: initial cluster var. S_5}
Recall from Example~\ref{exp:initial seed S_5} and Figure~\ref{fig:initial} the pseudoline arrangement $\pa(\hat \w_0)$ and the quiver $Q_{\hat \w_0}$ for $\hat\w_0\in S_5$.
To a face $F_{(i,j)}$ with $i\in[n-1]$ and $j\in[i+1,n]$ we associate following the above recipe the cluster variable $A_{(i,j)}:=\bar p_{i+1,\dots,j}$.
To the faces unbounded $F_i, i\in[4]$ to the left, we associate the variables $A_i:=\bar p_{5-i+1,\dots,5}$.
Note that the variables associated to the frozen vertices on the left (from bottom to top) are $\bar p_5,\bar p_{45},\bar p_{345},\bar p_{2345}$ and those associated to the frozen vertices on the right are $\bar p_{1},\bar p_{12},\bar p_{123},\bar p_{1234}$.
These Pl\"ucker coordinates are called \emph{consecutive minors}.
The collection of all cluster variables associated to this initial seed is
\[
\mathbf A_{\hat \w_0}=\{\bar p_{1},\bar p_{2},\bar p_{3},\bar p_4,\bar p_5,\bar p_{12},\bar p_{23},\bar p_{34},\bar p_{45},\bar p_{123},\bar  p_{234},\bar p_{345},\bar p_{1234},\bar
p_{2345}\}.
\]
\end{example}

\begin{example}\label{exp:initial seed w_0}
Consider $\hat\w_0\in S_n$ as in Examples~\ref{exp:pathGT} and \ref{exp:areaGT}. Then the collection $\mathbf A_{\hat \w_0}$ of associated cluster variables is
\[
\mathbf A_{\hat\w_0}=\{\bar p_{i,\dots, j}\mid i\in[n-1],j\in[i,n] \}, 
\]
where $\bar p_{i,\dots,j}$ is a frozen variable if either $i=1$ or $j=n$. Note that $\bar p_{1,\dots,n}=\det$, which is constant on $SL_n$, hence we disregard it. From now on we denote by $s_0$ the seed $s_0:=s_{\hat \w_0}=(\mathbf A_{\hat \w_0},Q_{\hat\w_0})$.
\end{example}

Berenstein-Fomin-Zelevinsky showed
\begin{theorem*}(\cite[Theorem~2.10]{BFZ05})
Let $w\in S_n$ with reduced expression $\w$. Then for the upper cluster algebra $\overline{\mathcal{Y}}(s_{\w})$ one has an isomorphism of algebras
\begin{align}\label{eq: cluster alg G^e,w}
\overline{\mathcal{Y}}(s_{\w})\otimes_{\mathbb Z} \mathbb C\cong \mathbb C[G^{e,w}].
\end{align}
\end{theorem*}

In particular, the Theorem implies the following: 
if $\w_1$ and $\w_2$ are two reduced expressions of $w\in S_n$ related by mutation in the sense of Definition~\ref{defn:mutation pa}, then the associated seeds $s_{\w_1}$ and $s_{\w_2}$ are related by cluster ($\A$-)mutation in the sense of \eqref{eq: A mutation}.
This explains our abuse of notation using the same letter $\mu$ for both types of mutation.

We now focus on the $\A$-cluster variety $G^{e,w_0}$ and the natural partial compactification using the frozen variables to study the superpotential as in \S\ref{sec: prep cluster}.
We partially compactify $G^{e,w_0}$ to $\bar G^{e,w_0}$ by allowing the frozen variables $\bar p_{[i]}$ and $\bar p_{[n-i,n]}$ for $i\in[n-1]$ to vanish.
Denote the resulting boundary divisor $D\subset \bar G^{e,w_0}$ and its irreducible components by
\[
D_{i}:=\{\bar p_{[i]}=0\}, \ \text{ resp. } \ D_{i,n}:=\{\bar p_{{[n-i,n]}}=0\}.
\]
There is an open embedding
$G^{e,w_0}\hookrightarrow SL_n/U$ given by  $g\mapsto g^tU$ and
up to codimension 2 the variety $\bar G^{e,w_0}$ agrees with $SL_n/U$ (this follows, for example, from \cite[Proposition~23]{Mag15}).
Hence, we have an isomorphism of rings $\mathbb C[\bar G^{e,w_0}]\cong \mathbb C[SL_n/U]$.
One of Magee's main results in \cite{Mag15} is the following.
\begin{theorem*}(\cite[Corollary~3]{Mag15})
The full Fock-Goncharov conjecture holds for $SL_n/U$.
\end{theorem*}

Moreover, Magee shows that there exists an optimized seed for every frozen vertex and therefore we can apply Algorithm~\ref{alg:superpot via opt seeds} stated in \S\ref{sec: prep cluster} to compute the superpotential. 
This is indeed what Magee did for the initial seed $s_0$ (see Example~\ref{exp:initial seed w_0}). 
Let $\X$ denote the Fock-Goncharov dual to the $\A$-cluster variety $G^{e,w_0}$ (see Definition~\ref{def:cluster variety}).
Recall that in the initial seed $s_0$ we have $N_{s_0}\cong \mathbb Z^{N+n-1}$ with basis $\{e_F\mid F\text{ face of }\pa(\w_0)\}$.
As before we set $e_{F_{(i,j)}}=:e_{(i,j)}$ and $e_{F_k}=:e_k$.
Further, recall that the superpotential $W:\X\to\mathbb C$ is given by the sum of $\vartheta$-functions associated to frozen variables. 
We denote by $\vartheta_{i}$ (resp. $\vartheta_{(i,n)}$) the $\vartheta$-functions associated to the frozen vertex $w_i$ (resp. $w_{(i,n)}$) in the initial quiver $Q_{s_0}$ (see Figure~\ref{fig:initial}) for $i\in[n-1]$.

\begin{proposition*}(\cite[Corollary~24]{Mag15})
Let $W:\X\to\mathbb C$ denote the superpotential. Then we have
$W\vert_{\X_{s_0}}=\sum_{i=1}^{n-1} \vartheta_{i}\vert_{\X_{s_0}}+\vartheta_{(n-i,n)}\vert_{\X_{s_0}}$, where 
\[
\vartheta_{i}\vert_{\X_{s_0}}=\sum_{k=0}^{n-1-i}z^{-e_i-\sum_{j=1}^k e_{(j,i+j)}}, \ \text{ and } \
\vartheta_{(i,n)}\vert_{\X_{s_0}}=\sum_{k=0}^{n-1-i}z^{-\sum_{j=0}^k e_{(i,n-j)}}, \ \text{ for } \ i\in[n-1].
\]
\end{proposition*}

\begin{example}
Consider $S_3$ and the initial seed with quiver $Q_{s_1s_2s_1}$. Then 
\begin{align*}
W\vert_{\X_{s_0}}&=\vartheta_{(1,3)}+\vartheta_{(2,3)}+\vartheta_{1}+\vartheta_{2}\\
&=z^{-e_{(1,3)}}+z^{-e_{(1,3)}-e_{(1,2)}}+z^{-e_{(2,3)}}+z^{-e_{1}}+z^{-e_{1}-e_{(1,2)}}+z^{-e_{2}}. \end{align*} 
\end{example}

\vbox{
\begin{definition}\label{def:super cone}
For $\w_0$ a reduced expression of $w_0\in S_n$ we define the following polyhedral objects by tropicalizing a sum of $\vartheta$-functions resp. the superpotential:
\begin{align*}
\Xi_{\w_0}&:=\{\mathbf x\in \mathbb R^{N+n-1}\mid W\vert_{\X_{\w_0}}^{\trop}(\mathbf x)\ge0\},\\
\mathsf \Xi_{\w_0}&:= \{\mathbf x\in \mathbb R^N\mid (\sum_{i=1}^{n-1}\vartheta_{(i,n)}\vert_{\X_{\w_0}})^{\trop}(\mathbf x)\ge 0\},\\
\Xi_{\w_0}(\lambda)&:= \Xi_{\w_0}\cap \tau^{-1}_{\w_0}(\lambda) \text{ for } \lambda\in \mathbb R^{n-1}.
\end{align*}
\end{definition}}

The $\A_{\text{prin}}$-construction in \cite{GHKK14} applied to our setting defines a flat family with base $\mathbb A^{N-2n+2}$ for every choice of seed, in particular for every ${\w}_0$.
The central fibre is by \cite[Theorem~8.39]{GHKK14} the toric variety associated to $\Xi_{\w_0}(\lambda)$ for $\lambda\in\mathbb Z_{>0}^{n-1}$.
One generic fibre is $SL_n/B$, hence we have a toric degeneration of the flag variety.
We do not go into the details on this construction but refer the reader to \cite[\S8]{GHKK14}.

\subsection{Relating to the area cones} Let $\w_0$ be an arbitrary reduced expression of $w_0\in S_n$. In what follows we show how to obtain an expression of the superpotential in any seed $s_{\w_0}$ associated to $\w_0$ by ``detropicalizing" the weighted cone $\mathcal S_{\w_0}$. 
We define it more generally for $\w$ a reduced expression of $w\in S_n$.
Denote by $\X_{\w}$ the cluster torus associated to the seed $s_{\w}$.

\begin{definition}\label{def:detrop area cone}
Let $\w$ be an arbitrary reduced expression of $w\in S_n$. Then the \emph{detropicalization} of the GP data, respectively the cone $\mathcal S_{\w}$ is defined as the function $W_{\mathcal S_{\w}}:\X_{\w}\to \mathbb C$ with
\begin{align}
W_{\mathcal S_{\w}}:=\sum_{\p\in\mathcal P_{\w}}z^{e_{\p}}+\sum_{i\in[n-1],0\le k\le n_i}z^{e_{[i:k]}}.
\end{align}
\end{definition}

The name is self-explanatory, observe that by definition we have
\[
\{\mathbf x\in \mathbb R^{\ell(w)+n-1}\mid W_{\mathcal S_{\w}}^{\trop}(\mathbf x)\ge 0\}=\mathcal S_{\w}.
\]

\begin{proposition}\label{prop: superpot area GT}
Let $\hat\w_0=s_1s_2s_1\cdots s_{n-1}s_{n-2}\cdots s_2s_1$ be the reduced expression associated to the initial seed $s_0$ as above. Then $W_{\mathcal S_{\hat\w_0}}=W\vert_{\X_{s_0}}$.
\end{proposition}

\begin{proof}
Recall from Example~\ref{exp:areaGT} the expressions $e_{\p_{i,j}}$ \eqref{eq: e_p in GT w_0} and $e_{[i:k]}$ \eqref{eq: f_i,k in GT w_0} for $i\in[n-1]$ and $j,k\in[i+1,n]$.
In comparison with \cite[Corollary~24]{Mag15} (restated above) we obtain
\begin{align*}
\vartheta_{(i,n)}\vert_{\X_{s_0}}= \sum_{j=i+1}^n z^{e_{\p_{i,j}}}, \ \text{ and } \
\vartheta_{i}\vert_{\X_{s_0}}= \sum_{k=0}^{n-1-i} z^{e_{[i:k]}}.
\end{align*}
As from Example~\ref{exp:pathGT} we know $\mathcal P_{\w_0}=\{\p_{i,j}\mid i\in[n-1],j\in[i+1,n]\}$, the claim follows.
\end{proof}

%%%%%%%%%%%%%%%%%%%%%%%%%%%%%%%%% MUTATION AREA

\subsection{Mutation of the area cone}\label{subsubsec: path and area mut}

Our aim is to generalize Proposition~\ref{prop: superpot area GT} for arbitrary reduced expressions $\w_0$. 
We achieve this by showing that the detropicalization of $\mathcal S_{\w_0}$ behaves as the superpotential does when applying $\mathcal X$-mutation. 
Further, we show that if $\mu(\w)$ and $\w$ are reduced expressions of $w\in S_n$, then $\mu^*(W_{\mathcal S_{\mu(\w)}})=W_{\mathcal S_{\w}}$, where $\mu^*: \mathbb C[\X_{\mu(\w)}]\to\mathbb C[\X_{\w}]$ is the pull-back of the cluster mutation as in \eqref{eq: def pullback X-mut}.
This follows from Lemma~\ref{lem: mutation paths} and Lemma~\ref{lem:mutation wt area}.
Recall from Definition~\ref{defn:mutation pa} the mutation of pseudoline arrangements.
The core of this subsection is the case-by-case analysis of how mutation effects GP-paths.

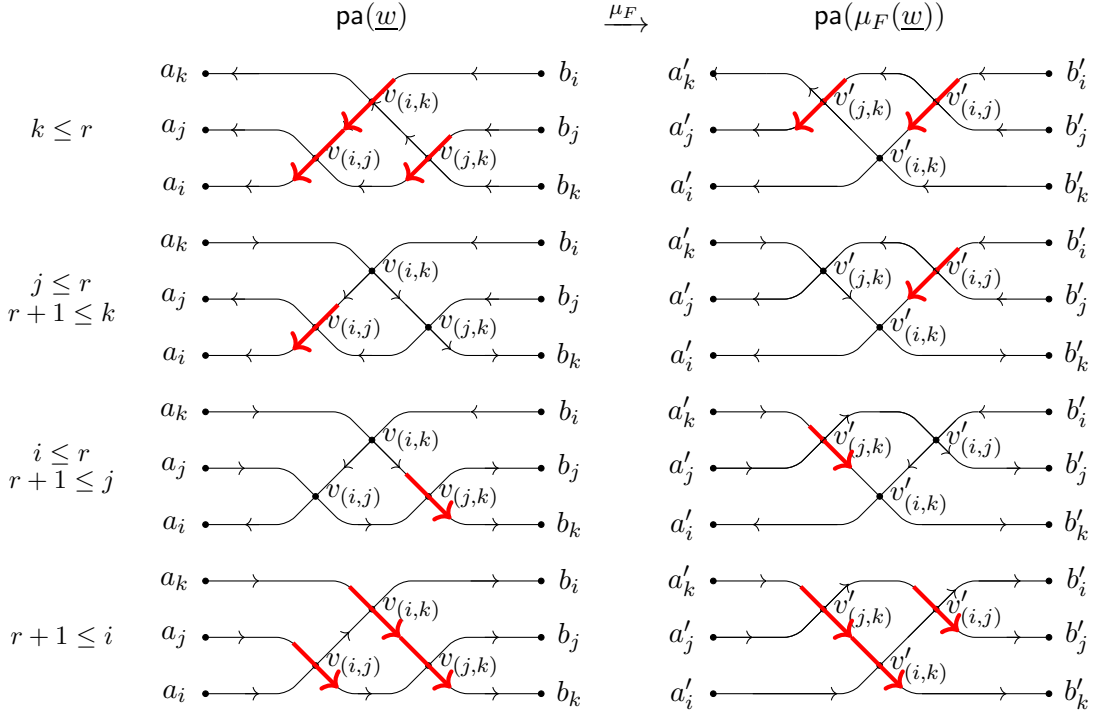
\begin{figure}
\centering
\begin{center}
\begin{tikzpicture}[scale=0.75]

%%%%%%%%%%%%%%%%%%%%%%%
%% LRR
\begin{scope}[yshift=-6cm]
\node at (-2,1.25) {\small $i\le r$};
\node at (-2,.75) {\small $r+1\le j$};

\draw [fill] (.55,2) circle [radius=0.05];
\draw [fill] (.55,1) circle [radius=0.05];
\draw [fill] (.55,0) circle [radius=0.05];
\draw [fill] (6.5,2) circle [radius=0.05];
\draw [fill] (6.5,1) circle [radius=0.05];
\draw [fill] (6.5,0) circle [radius=0.05];
\node at (0,2) {$a_{k}$};
\node at (0,1) {$a_{j}$};
\node at (0,0) {$a_{i}$};
\node at (7,2) {$b_{i}$};
\node at (7,1) {$b_{j}$};
\node at (7,0) {$b_{k}$};

\draw[fill] (3.5,1.5) circle [radius=.05];
\draw[fill] (2.5,.5) circle [radius=.05];
\draw[fill] (4.5,.5) circle [radius=.05];

\node[right] at (3.5,1.5) {$v_{(i,k)}$};
\node[right] at (2.5,.5) {$v_{(i,j)}$};
\node[right] at (4.5,.5) {$v_{(j,k)}$};

    \draw[rounded corners] (0.5,0) --(1.25,0) -- (2,0)-- (3,1) --(4,2) -- (5.25,2);
    \draw[->, rounded corners] (1.5,0) -- (1.25,0);
    \draw[->, rounded corners] (3.5,1.5) -- (3,1);
    \draw[->, rounded corners] (6.5,2) -- (5.25,2);
\draw[rounded corners] (1.25,1) -- (2,1) -- (3,0) -- (3.5,0)-- (4,0) -- (5,1) -- (5.57,1)-- (6.5,1);
    \draw[->] (0.5,1) -- (1.25,1);
    \draw[->] (3.25,0) -- (3.5,0);
    \draw[->] (5.5,1) -- (5.75,1);
\draw[rounded corners] (1.5,2) -- (3,2) -- (4,1) -- (5,0) -- (5.75,0) -- (6.5,0);
    \draw[->] (0.5,2) -- (1.5,2);
    \draw[->] (3.5,1.5) -- (4,1);
    \draw[->] (5.5,0) -- (5.75,0);
    \draw[red, ->, ultra thick] (4.1,0.9) -- (4.9,0.1);
    
\begin{scope}[xshift=9cm]
\draw [fill] (.55,2) circle [radius=0.05];
\draw [fill] (.55,1) circle [radius=0.05];
\draw [fill] (.55,0) circle [radius=0.05];
\draw [fill] (6.5,2) circle [radius=0.05];
\draw [fill] (6.5,1) circle [radius=0.05];
\draw [fill] (6.5,0) circle [radius=0.05];
\node at (0,2) {$a'_{k}$};
\node at (0,1) {$a'_{j}$};
\node at (0,0) {$a'_{i}$};
\node at (7,2) {$b'_{i}$};
\node at (7,1) {$b'_{j}$};
\node at (7,0) {$b'_{k}$};

\draw[fill] (3.5,.5) circle [radius=.05];
\draw[fill] (2.5,1.5) circle [radius=.05];
\draw[fill] (4.5,1.5) circle [radius=.05];
\node[right] at (3.5,.5) {$v'_{(i,k)}$};
\node[right] at (2.5,1.5) {$v'_{(j,k)}$};
\node[right] at (4.5,1.5) {$v'_{(i,j)}$};

\draw[rounded corners] (0.5,0) --(1.25,0) -- (3,0)-- (4,1) -- (5,2) -- (5.25,2);
\draw[rounded corners] (0.5,1) -- (2,1) -- (3,2) -- (4,2) -- (5,1) -- (6.5,1);
\draw[rounded corners] (1.5,2) -- (2,2) -- (3,1)-- (4,0) -- (5.75,0) -- (6.5,0);
 \draw[->] (0.5,1) -- (1.25,1);
\draw[->, red , ultra thick]  (2.25,1.75) -- (3,1);
    \draw[->] (2.75,0) -- (1.25,0);
   \draw[->] (4.85,1.85) -- (4,1);
    \draw[->] (6.5,2) -- (5.25,2);
\draw[->, rounded corners] (1.25,1) -- (2,1) -- (2.95,1.95);
\draw[->, rounded corners] (3.5,2) -- (4,2) -- (4.75,1.25);
\draw[->, rounded corners] (5.75,1)-- (6,1);
    \draw[->] (0.5,2) -- (1.5,2);
    \draw[->] (4.25,0) -- (5.75,0);      

\end{scope}
\end{scope}
%%%%%%%%%%%%%%%%%%%%%%%%%%%%%%%%
%%% RRR
\begin{scope}[yshift=-9cm]

\node at (-2,1) {\small $r+1\le i$};

\draw [fill] (.55,2) circle [radius=0.05];
\draw [fill] (.55,1) circle [radius=0.05];
\draw [fill] (.55,0) circle [radius=0.05];
\draw [fill] (6.5,2) circle [radius=0.05];
\draw [fill] (6.5,1) circle [radius=0.05];
\draw [fill] (6.5,0) circle [radius=0.05];
\node at (0,2) {$a_{k}$};
\node at (0,1) {$a_{j}$};
\node at (0,0) {$a_{i}$};
\node at (7,2) {$b_{i}$};
\node at (7,1) {$b_{j}$};
\node at (7,0) {$b_{k}$};

\draw[fill] (3.5,1.5) circle [radius=.05];
\draw[fill] (2.5,.5) circle [radius=.05];
\draw[fill] (4.5,.5) circle [radius=.05];

\node[right] at (3.5,1.5) {$v_{(i,k)}$};
\node[right] at (2.5,.5) {$v_{(i,j)}$};
\node[right] at (4.5,.5) {$v_{(j,k)}$};

    \draw[rounded corners] (0.5,0) --(1.25,0) -- (2,0)-- (3,1) --(4,2) -- (5.25,2) -- (6.5,2);
    \draw[<-, rounded corners] (1.5,0) -- (1.25,0);
    \draw[<-, rounded corners] (3.1,1.1) -- (3,1);
    \draw[<-, rounded corners] (5.75,2) -- (5.25,2);
\draw[rounded corners] (1.25,1) -- (2,1) -- (3,0) -- (3.5,0)-- (4,0) -- (5,1) -- (5.57,1)-- (6.5,1);
    \draw[->] (0.5,1) -- (1.25,1);
    \draw[->] (3.25,0) -- (3.5,0);
    \draw[->] (5.5,1) -- (5.75,1);
\draw[rounded corners] (1.5,2) -- (3,2) -- (4,1) -- (5,0) -- (5.75,0) -- (6.5,0);
    \draw[->] (0.5,2) -- (1.5,2);
%    \draw[->] (3.5,1.5) -- (4,1);
    \draw[->] (5.5,0) -- (5.75,0);    
\draw[->,red, ultra thick] (2.1,0.9) -- (2.9,0.1);  
\draw[->,red, ultra thick] (3.1,1.9) -- (4,1);
\draw[red, ultra thick, ->] (4,1) -- (4.9,0.1);

\begin{scope}[xshift=9cm]

\draw [fill] (.55,2) circle [radius=0.05];
\draw [fill] (.55,1) circle [radius=0.05];
\draw [fill] (.55,0) circle [radius=0.05];
\draw [fill] (6.5,2) circle [radius=0.05];
\draw [fill] (6.5,1) circle [radius=0.05];
\draw [fill] (6.5,0) circle [radius=0.05];
\node at (0,2) {$a'_{k}$};
\node at (0,1) {$a'_{j}$};
\node at (0,0) {$a'_{i}$};
\node at (7,2) {$b'_{i}$};
\node at (7,1) {$b'_{j}$};
\node at (7,0) {$b'_{k}$};

\draw[fill] (3.5,.5) circle [radius=.05];
\draw[fill] (2.5,1.5) circle [radius=.05];
\draw[fill] (4.5,1.5) circle [radius=.05];
\node[right] at (3.5,.5) {$v'_{(i,k)}$};
\node[right] at (2.5,1.5) {$v'_{(j,k)}$};
\node[right] at (4.5,1.5) {$v'_{(i,j)}$};

\draw[rounded corners] (0.5,0) --(1.25,0) -- (3,0)-- (4,1) -- (5,2) -- (5.25,2);
\draw[rounded corners] (0.5,1) -- (2,1) -- (3,2) -- (4,2) -- (5,1) -- (6.5,1);
\draw[rounded corners] (1.5,2) -- (2,2) -- (3,1)-- (4,0) -- (5.75,0) -- (6.5,0);
    \draw[<-] (2.75,0) -- (1.25,0);
    \draw[<-] (4.85,1.85) -- (4,1);
    \draw (6.5,2) -- (5.25,2);
    \draw[->] (5.25,2) -- (6,2);
\draw[->, rounded corners] (1.25,1) -- (2,1) -- (2.95,1.95);
\draw[->, rounded corners] (5.75,1)-- (6,1);
\draw[->] (4.25,0) -- (5.75,0);   
    \draw[->] (0.5,2) -- (1.5,2); 
     \draw[->] (0.5,1) -- (1.25,1);
    \draw[->, red, ultra thick] (4.1,1.9) -- (4.9,1.1);
    \draw[->, red, ultra thick] (2.1,1.9) -- (3,1);
    \draw[red, ultra thick, ->] (3,1) -- (3.9,0.1);

\end{scope}
\end{scope}

%%% LLL
\begin{scope}[yshift=0cm]
\node at (3.5,3) {$\pa(\w)$};
\node at (8,3) {$\xrightarrow{\mu_{F}}$};

\node at (-2,1) {\small$k\le r$};

\draw [fill] (.55,2) circle [radius=0.05];
\draw [fill] (.55,1) circle [radius=0.05];
\draw [fill] (.55,0) circle [radius=0.05];
\draw [fill] (6.5,2) circle [radius=0.05];
\draw [fill] (6.5,1) circle [radius=0.05];
\draw [fill] (6.5,0) circle [radius=0.05];
\node at (0,2) {$a_{k}$};
\node at (0,1) {$a_{j}$};
\node at (0,0) {$a_{i}$};
\node at (7,2) {$b_{i}$};
\node at (7,1) {$b_{j}$};
\node at (7,0) {$b_{k}$};

\draw[fill] (3.5,1.5) circle [radius=.05];
\draw[fill] (2.5,.5) circle [radius=.05];
\draw[fill] (4.5,.5) circle [radius=.05];

\node[right] at (3.5,1.5) {$v_{(i,k)}$};
\node[right] at (2.5,.5) {$v_{(i,j)}$};
\node[right] at (4.5,.5) {$v_{(j,k)}$};

\draw[rounded corners] (0.5,0) --(1.25,0) -- (2,0)-- (3,1) --(4,2) -- (5.25,2);
    \draw[->, rounded corners] (1.5,0) -- (1.25,0);
    \draw[->, rounded corners] (3.5,1.5) -- (3,1);
    \draw[->, rounded corners] (6.5,2) -- (5.25,2);
\draw[rounded corners] (.5,1) -- (2,1) -- (3,0) -- (3.5,0)-- (4,0) -- (5,1) -- (5.57,1)-- (6.5,1);
    \draw[<-] (1,1) -- (1.25,1);
    \draw[<-] (3.25,0) -- (3.5,0);
    \draw[<-] (5.5,1) -- (5.75,1);
\draw[rounded corners] (.5,2) -- (3,2) -- (4,1) -- (5,0) -- (5.75,0) -- (6.5,0);
    \draw[<-] (1,2) -- (1.5,2);
    \draw[<-] (3.5,1.5) -- (4,1);
    \draw[<-] (5.5,0) -- (5.75,0);
    \draw[<-] (4.1,0.9) -- (4.9,0.1);
   
\draw[->,red, ultra thick] (3,1) -- (2.1,0.1);
\draw[->,red, ultra thick] (3.9,1.9) -- (3,1);
\draw[<-,red, ultra thick] (4.1,0.1) -- (4.9,.9);
    
\begin{scope}[xshift=9cm]
\node at (3.5,3) {$\pa(\mu_F(\w))$};

\draw [fill] (.55,2) circle [radius=0.05];
\draw [fill] (.55,1) circle [radius=0.05];
\draw [fill] (.55,0) circle [radius=0.05];
\draw [fill] (6.5,2) circle [radius=0.05];
\draw [fill] (6.5,1) circle [radius=0.05];
\draw [fill] (6.5,0) circle [radius=0.05];
\node at (0,2) {$a'_{k}$};
\node at (0,1) {$a'_{j}$};
\node at (0,0) {$a'_{i}$};
\node at (7,2) {$b'_{i}$};
\node at (7,1) {$b'_{j}$};
\node at (7,0) {$b'_{k}$};

\draw[fill] (3.5,.5) circle [radius=.05];
\draw[fill] (2.5,1.5) circle [radius=.05];
\draw[fill] (4.5,1.5) circle [radius=.05];
\node[right] at (3.5,.5) {$v'_{(i,k)}$};
\node[right] at (2.5,1.5) {$v'_{(j,k)}$};
\node[right] at (4.5,1.5) {$v'_{(i,j)}$};

\draw[rounded corners] (0.5,0) --(1.25,0) -- (3,0)-- (4,1) -- (5,2) -- (5.25,2);
\draw[rounded corners] (0.5,1) -- (2,1) -- (3,2) -- (4,2) -- (5,1) -- (5.75,1);
\draw[rounded corners] (1.5,2) -- (2,2) -- (3,1)-- (4,0) -- (5.75,0) -- (6.5,0);
    \draw[->] (2.75,0) -- (1.25,0);
    \draw[->,red, ultra thick] (4.9,1.9) -- (4,1);
    \draw[->] (6.5,2) -- (5.25,2);
\draw[<-, rounded corners] (1.25,1) -- (2,1) -- (2.75,1.75);
\draw[<-,red, ultra thick] (2,1) -- (2.9,1.9);
\draw[<-, rounded corners] (3.5,2) -- (4,2) -- (4.75,1.25);
\draw[<-, rounded corners] (5.57,1)-- (6.5,1);
    \draw[<-] (0.5,2) -- (1.5,2);
    \draw[<-] (2.25,1.75) -- (3,1);
    \draw[<-] (4.25,0) -- (5.75,0);

\end{scope}
\end{scope}

%%% LLR
\begin{scope}[yshift=-3cm]
\node at (-2,1.25) {\small$j\le r$};
\node at (-2,.75) {\small$r+1\le k$};
\draw [fill] (.55,2) circle [radius=0.05];
\draw [fill] (.55,1) circle [radius=0.05];
\draw [fill] (.55,0) circle [radius=0.05];
\draw [fill] (6.5,2) circle [radius=0.05];
\draw [fill] (6.5,1) circle [radius=0.05];
\draw [fill] (6.5,0) circle [radius=0.05];
\node at (0,2) {$a_{k}$};
\node at (0,1) {$a_{j}$};
\node at (0,0) {$a_{i}$};
\node at (7,2) {$b_{i}$};
\node at (7,1) {$b_{j}$};
\node at (7,0) {$b_{k}$};

\draw[fill] (3.5,1.5) circle [radius=.05];
\draw[fill] (2.5,.5) circle [radius=.05];
\draw[fill] (4.5,.5) circle [radius=.05];

\node[right] at (3.5,1.5) {$v_{(i,k)}$};
\node[right] at (2.5,.5) {$v_{(i,j)}$};
\node[right] at (4.5,.5) {$v_{(j,k)}$};

\draw[rounded corners] (0.5,0) --(1.25,0) -- (2,0)-- (3,1) --(4,2) -- (5.25,2);
    \draw[->, rounded corners] (1.5,0) -- (1.25,0);
    \draw[->, rounded corners] (3.5,1.5) -- (3,1);
    \draw[->, rounded corners] (6.5,2) -- (5.25,2);
\draw[rounded corners] (.5,1) -- (2,1) -- (3,0) -- (3.5,0)-- (4,0) -- (5,1) -- (5.57,1)-- (6.5,1);
    \draw[<-] (1,1) -- (1.25,1);
    \draw[<-] (3.25,0) -- (3.5,0);
    \draw[<-] (5.5,1) -- (5.75,1);
\draw[rounded corners] (1.5,2) -- (3,2) -- (4,1) -- (5,0) -- (5.75,0) -- (6.5,0);
    \draw[->] (0.5,2) -- (1.5,2);
    \draw[->] (3.5,1.5) -- (4,1);
    \draw[->] (5.5,0) -- (5.75,0);
    \draw[->] (4.1,0.9) -- (4.9,0.1);
\draw[<-,red, ultra thick] (2.1,0.1)-- (2.9,0.9);

\begin{scope}[xshift=9cm]
\draw [fill] (.55,2) circle [radius=0.05];
\draw [fill] (.55,1) circle [radius=0.05];
\draw [fill] (.55,0) circle [radius=0.05];
\draw [fill] (6.5,2) circle [radius=0.05];
\draw [fill] (6.5,1) circle [radius=0.05];
\draw [fill] (6.5,0) circle [radius=0.05];

\node at (0,2) {$a'_{k}$};
\node at (0,1) {$a'_{j}$};
\node at (0,0) {$a'_{i}$};
\node at (7,2) {$b'_{i}$};
\node at (7,1) {$b'_{j}$};
\node at (7,0) {$b'_{k}$};

\draw[fill] (3.5,.5) circle [radius=.05];
\draw[fill] (2.5,1.5) circle [radius=.05];
\draw[fill] (4.5,1.5) circle [radius=.05];
\node[right] at (3.5,.5) {$v'_{(i,k)}$};
\node[right] at (2.5,1.5) {$v'_{(j,k)}$};
\node[right] at (4.5,1.5) {$v'_{(i,j)}$};

\draw[rounded corners] (0.5,0) --(1.25,0) -- (3,0)-- (4,1) -- (5,2) -- (5.25,2);
\draw[rounded corners] (0.5,1) -- (2,1) -- (3,2) -- (4,2) -- (5,1) -- (5.75,1);
\draw[rounded corners] (1.5,2) -- (2,2) -- (3,1)-- (4,0) -- (5.75,0) -- (6.5,0);
    \draw[->] (2.75,0) -- (1.25,0);
    \draw[->,red, ultra thick] (4.9,1.9) -- (4,1);
    \draw[->] (6.5,2) -- (5.25,2);
\draw[<-, rounded corners] (1.25,1) -- (2,1) -- (2.75,1.75);
\draw[<-, rounded corners] (3.5,2) -- (4,2) -- (4.75,1.25);
\draw[<-, rounded corners] (5.57,1)-- (6.5,1);
    \draw[->] (0.5,2) -- (1.5,2);
    \draw[->] (2.25,1.75) -- (3,1);
    \draw[->] (4.25,0) -- (5.75,0);
\end{scope}
\end{scope}

\end{tikzpicture}
\end{center}
\caption{The pseudoline arrangement $\pa(\w)$ (resp. $\pa(\mu_F(\w))$) locally around the face $F=F_{(i,j)}$ (resp. $F'=F'_{(j,k)}$) bounded by lines $l_i,l_j,l_k$ with $i<j<k$ and orientations $(l_r,l_{r+1})$ for all possible $r$.
The red arrows are those forbidden in GP-paths.}\label{fig:loc.orient}
\end{figure}

In Figure~\ref{fig:loc.orient} we display locally around the mutable face $F=F_{(i,j)}$ (resp. $F'$) the orientations of $\pa(\w)$ (resp. $\pa(\mu_F(\w))$). 
The red arrows indicate which passages are forbidden in GP-paths. 
In Tables~\ref{tab:(lll) cases of path mut} to \ref{tab:(rrr) cases of path mut} we list in the second column all possibilities how a GP-path ${\p}$ locally looks around the face $F$. 
In the third column of each table is a complete list of how GP-paths look locally around the face $F'$ obtained from $F$ by mutation $\mu_F$.

Recall the arrows for the quiver corresponding to $\pa(\w)$ and $\pa(\mu_F(\w))$ from Figure~\ref{fig:pseudo.mut}.
We call a face $E$ \emph{incoming} (resp. \emph{outgoing}) with respect to $F$ in $\pa(\w)$, if there is an arrow in the quiver $Q_{\w}$ from (resp. to) the vertex corresponding to $E$ to (resp. from) the vertex corresponding to $F$. 
We denote by $\text{In}_{F}$ the union of all incoming faces and by $\text{Out}_{F}$ the union of all outgoing faces.
See for example, Figure~\ref{fig:pseudo.mut}.

\vbox{
\begin{definition}\label{def:local type}
Let $\p\in\mathcal P_{\w}$ for $\w\in S_n$ and consider a mutable face $F$ of $\pa(\w)$. 
Set $\delta_{F\subset \area_{\p}}:=1$ if $F\in \area_{\p}$ and zero otherwise.
Then we define the \emph{$F$-local type} of $\p$ as the triple
\[
F(\p):=(i_{F,{\p}},x_{F,{\p}},o_{F,{\p}}):=(\#\{\text{In}_{F}\cap \area_{\p}\}, 
\delta_{F\in \area_{\p}},
\#\{\text{Out}_{F}\cap \area_{\p}\}).
\]

\end{definition}}
For example, if $\area_{\p}$ in Figure~\ref{fig:pseudo.mut} contains the faces $F,F_{\init_1}$ and $F_{\text{out}_2}$ but not $F_{\init_1}$ and $F_{\text{out}_1}$, then the $F$-local type of $\p$ is $(1,1,1)$.
The following lemma is a crucial observation on the $F$-local type of GP-paths.

    %%%%%%%%%%%%%%%%%%%%%%%%%%%
    
    \begin{table}
    \centering %\{.7}{
    \scalebox{.8}{
    \begin{tabular}{|c|c|c|c|} \hline
$F$-local type of $\p$ & ${\p}$ in $\pa(\w)$ & 
    $\p'=\mut_{F}({\p})$ in $\pa(\mu_F(\w))$ & $F'$-local type of $\p'$\\ \hline\hline
$(2,1,2)$ &
    $b_i\to v_{(i,k)}\to a_k$ &
    $b'_i\to v'_{(i,j)}\to v'_{(j,k)} \to a'_k$ &
     $(2,1,2)$ \\ \hline

$(1,1,2)$ &
    $b_{j}\to v_{(j,k)}\to v_{(i,k)}\to a_k$ &
    \begin{tabular}{c} 
        $b'_j\to v'_{(i,j)}\to v'_{(j,k)}\to a'_k$ \\ \hline
        $b'_j\to v'_{(i,j)}\to v'_{(i,k)}\to v'_{(j,k)}\to a'_k$ \end{tabular} & 
    \begin{tabular}{c}
        $(2,1,1)$ \\ \hline
        $(2,0,1)$
    \end{tabular}   \\ \hline 

 $(1,1,1)$&        
    $ b_j\to v_{(j,k)}\to v_{(i,k)}\to v_{(i,j)}\to a_j$ &
    $ b'_j\to v'_{(i,j)}\to v'_{(i,k)}\to v'_{(j,k)}\to a'_j$ &
     $(1,0,1)$ \\ \hline
    
 $(1,1,1)$ &
        $b_k \to v_{(j,k)}\to v_{(i,k)}\to a_k $   &
        $b'_k\to v'_{(i,k)}\to v'_{(j,k)}\to a'_k$     &
          $(1,0,1)$   \\ \hline
    
\begin{tabular}{c}
      $(1,1,0)$ \\ \hline
      $(1,0,0)$ \end{tabular} &    
    \begin{tabular}{c}
         $b_k\to v_{(j,k)}\to v_{(i,k)}\to v_{(i,j)}\to a_j$  \\ \hline
          $b_k\to v_{(j,k)}\to v_{(i,j)}\to a_j$
    \end{tabular} &
    $ b'_k\to v'_{(i,k)}\to v'_{(j,k)}\to a'_j$ &
          $(0,0,1)$ \\ \hline
  
 $(0,0,0)$ &
    $b_k\to v_{(j,k)}\to v_{(i,j)}\to a_i$ &
    $b'_k\to v'_{(i,k)}\to a'_i$ &
          $(0,0,0)$ \\ \hline 
    
        \end{tabular}}
    \caption{Shapes of paths locally around $F$ (resp. $F'$) in $\mathcal P_{\w}$ (resp. $\mathcal P_{\mu_F(\w)}$) for orientation $(l_r,l_{r+1})$ with $i<j<k\le r$ (see Figure~\protect{\ref{fig:loc.orient}}) and how they are mapped onto each other by $\mut_F$.}
    \label{tab:(lll) cases of path mut}
\end{table}

\vbox{
\begin{lemma}\label{lem:possible F-local types}
Let $\p\in\mathcal P_{\w}$ for $\w\in S_n$ and consider a mutable face $F$ of $\pa(\w)$. Then the following are all possible $F$-local types $\p$ can have:
\begin{itemize}
    \item[$i_{F,\p}=o_{F,\p}$:] then $F(\p) \in \{(0,0,0),(1,0,1),(1,1,1),(2,1,2)\}$;
    \item[$i_{F,\p}<o_{F,\p}$:] then $F(\p) \in \{(1,1,2),(0,0,1)\}$;
    \item[$i_{F,\p}>o_{F,\p}$:] then $F(\p) \in \{(1,0,0),(1,1,0),(2,0,1),(2,1,1)\}$.
\end{itemize}
Moreover, the $F$-local types of $\p$ with $i_{F,\p}>o_{F,\p}$ come in pairs as $((1,0,0),(1,1,0))$ or $((2,0,1),(2,1,1))$. Meaning that if a path of one type exists for a fixed orientation then so does a path of the corresponding other type for the same orientation.
\end{lemma}}

\begin{proof}
The lemma follows from case-by-case consideration of all possible shapes of $\p\in\mathcal P_{\w}$ around a mutable face $F$ of $\pa(\w)$. First observe, that $F$ can have two different shapes, depending on whether it is defined by simple reflections $s_ms_{m+1}s_m$ (as on the left in Figure~\ref{fig:pseudo.mut}) or by $s_{m+1}s_ms_{m+1}$ (as on the right in Figure~\ref{fig:pseudo.mut}).
We endow $\pa(\w)$ for either case of $F$ with all possible orientations $(l_r,l_{r+1})$. Then locally at $F$, there are four cases of orientation depending on $r$ and $r+1$ in relation to $i,j,k$ (see Figure~\ref{fig:loc.orient}).
We consider all possibilities for the path $\p$ to pass $F$ for each case of orientation and shape of $F$.
These are listed in Tables~\ref{tab:(lll) cases of path mut} to \ref{tab:(rrr) cases of path mut}, in the second column for $F$ as on the left of Figure~\ref{fig:loc.orient} and in the third for $F$ as on the right of Figure~\ref{fig:loc.orient}.
In the first and last columns of these tables we indicate the corresponding $F$-local type.
Observe, that the list in the claim of the lemma covers all occurring $F$-local types.

Regarding the second part of the claim, this also follows as an observation from Tables~\ref{tab:(lll) cases of path mut} to \ref{tab:(rrr) cases of path mut}.
\end{proof}

    \begin{table}
    \centering %\{.7}{
    \scalebox{.8}{
    \begin{tabular}{|c|c|c|c|} \hline
$F$-local type of $\p$ & ${\p}$ in $\pa(\w)$ & 
    $\p'=\mut_{F}({\p})$ in $\pa(\mu_F(\w))$ & $F'$-local type of $\p'$\\ \hline\hline

 $(1,0,1)$ &
    $b_i\to v_{(i,k)}\to v_{(j,k)}\to v_{(i,j)}\to a_i $     &
    $b'_i\to v'_{(i,j)} \to v'_{(j,k)} \to v'_{(i,k)}\to a'_i$ &
          $(1,1,1)$ \\ \hline

\begin{tabular}{c}
      $(2,1,1)$ \\ \hline
      $(2,0,1)$
    \end{tabular} &
    \begin{tabular}{c}
        $b_i\to v_{(i,k)}\to v_{(i,j)}\to a_j$   \\ \hline 
         $b_i\to v_{(i,k)}\to v_{(j,k)}\to v_{(i,j)}\to a_j$
    \end{tabular} &
    $ b'_i\to v'_{(i,j)}\to v'_{(j,k)}\to a'_j$ &
          $(1,1,2)$ \\ \hline

 $(1,0,1)$  &
    $b_i\to v_{(i,k)}\to v_{(j,k)}\to b_k$ &
    $ b'_i\to v'_{(i,j)}\to v'_{(j,k)}\to v'_{(i,k)}\to b'_k$ &
      $(1,1,1)$  \\ \hline
    
 $(0,0,1)$  &
    $b_j\to v_{(j,k)}\to v_{(i,j)}\to a_i$ &
    \begin{tabular}{c}
    $b'_j\to v'_{(i,j)}\to v'_{(j,k)}\to v'_{(i,k)}\to a'_i$\\
    $b'_i\to v'_{(i,j)}\to v'_{(i,k)}\to a'_i$
    \end{tabular}&
    \begin{tabular}{c}
          $(1,1,0)$  \\ \hline
          $(1,0,0)$
    \end{tabular} \\ \hline    

 $(1,0,1)$  &
    $b_j\to v_{(j,k)}\to v_{(i,j)}\to a_j $ &
    $b'_j\to v'_{(i,j)}\to v'_{(j,k)}\to a'_j $ &
     $(1,1,1)$  \\ \hline    
    
 $(0,0,1)$  &
    $b_{j}\to v_{(j,k)}\to b_k  $ &
    \begin{tabular}{c}
         $b'_j\to v'_{(i,j)} \to v'_{(j,k)}\to v'_{(i,k)}\to b'_k $  \\ \hline 
          $b'_j\to v'_{(i,j)}\to v'_{(i,k)}\to b'_k$
    \end{tabular} &
    \begin{tabular}{c}
         $(1,1,0)$ \\ \hline
          $(1,0,0)$ 
    \end{tabular} \\ \hline      

 $(1,0,1)$   &
    $a_k\to v_{(i,k)}\to v_{(j,k)}\to v_{(i,j)}\to a_i $ &
    $ a'_k\to v'_{(j,k)}\to v'_{(i,k)}\to a'_i $ &
     $(1,1,1)$ \\ \hline  

\begin{tabular}{c}
      $(2,1,1)$  \\ \hline
      $(2,0,1)$ 
    \end{tabular}&
    \begin{tabular}{c}
         $  a_k\to v_{(i,k)}\to v_{(i,j)}\to a_j $ \\ \hline 
          $a_k\to v_{(i,k)}\to v_{(j,k)}\to v_{(i,j)}\to a_j $
    \end{tabular} &
    $ a'_k\to v'_{(j,k)}\to a'_j $ &
      $(1,1,2)$ \\ \hline      
    
 $(1,0,1)$ &
    $a_k\to v_{(i,k)}\to v_{(j,k)}\to b_k$ &
    $ a'_k\to v'_{(j,k)}\to v'_{(i,k)}\to b'_k$ &
     $(1,1,1)$  \\ \hline   

        \end{tabular}}
    \caption{Shapes of paths locally around $F$ (resp. $F'$) in $\mathcal P_{\w}$ (resp. $\mathcal P_{\mu_F(\w)}$) for orientation $(l_r,l_{r+1})$ with $i<j\le r$ and $r+1\le k$ (see Figure~\protect{\ref{fig:loc.orient}}) and how they are mapped onto each other by $\mut_F$.}
    \label{tab:(llr) cases of path mut}
\end{table}

With notation as in the lemma, if $\p_1,\p_2$ are paths with $i_{F,\p_j}>o_{F,\p_j}, j=1,2$ such that $((i_{F,\p_1},x_{F,\p_1},o_{F,\p_1}),(i_{F,\p_2},x_{F,\p_2},o_{F,\p_2}))$ is one of the pairs, then we denote by $\p_1\oplus\p_2$ their formal sum.
If $\p_1$ and $\p_2$ are equal away from $F$, we denote this by $\p_1/F=\p_2/F$.
Observe, that this is the case here.
With this notation we define the following set of paths, respectively formal sums of paths.
\begin{align}\label{eq: def P_w,F}
\widetilde{\mathcal{P}}_{\w,F} := \left\{ \begin{matrix} \p,\\ \p_1\oplus\p_2\end{matrix} \right\vert 
\left.\begin{matrix}
\p\in\mathcal{P}_{\w} \text{ with }  i_{F,\p}=o_{F,\p} \text{ or } i_{F,\p}<o_{F,\p},\\
 \p_1,\p_2\in\mathcal P_{\w} \text{ with } i_{F,\p_j}>o_{F,\p_j}\text{ for } j=1,2
\end{matrix}\right\}.
\end{align}

Note that for every mutable face $F$ of $\pa(\w)$ every path in $\mathcal{P}_{\w}$ appears in $\widetilde{\mathcal{P}}_{\w,F}$ either on its own or as a formal summand. This additional structure on $\mathcal{P}_{\w}$ allows us to define \emph{mutation} on it.

\begin{definition}\label{def:mut path}
Let $w\in S_n$ with reduced expressions $\w$ and $\mu_F(\w)$, where $F$ is a mutable face in $\pa(\w)$. Denote by $F'$ the corresponding face in $\pa(\mu_F(\w))$. We define  $\mut_F:\widetilde{\mathcal{P}}_{\w,F}\to \widetilde{\mathcal{P}}_{\mu_F(\w),F'}$ depending on the $F$-local type by 
\begin{itemize}
    \item[$i_{F,\p}=o_{F,\p}$:] $\mut_F(\p)=\p'$ with $\p/ F=\p'/ F'$, where for $F(\p)\in\{(0,0,0),(2,1,2)\}$ we have $F(\p)=F'(\p')$, and for $F(\p)\in\{(1,0,1),(1,1,1)\}$ we have $F'(\p')=(i_{F,\p},\vert x_{F,\p}-1\vert,o_{F,\p})$;
    \item[$i_{F,\p}<o_{F,\p}$:] $\mut_F(\p)=\p_1'\oplus\p_2'$ with $\p/F=\p'_1/F'=\p'_2/F'$, for $F(\p)\in\{(0,0,1),(1,1,2)\}$ with $F'(\p_1')=(o_{F,\p},x_{F,\p},i_{F,\p})$ and $F'(\p'_2)=(o_{F,\p},\vert x_{F,\p}-1\vert,i_{F,\p})$;
    \item[$i_{F,\p}>o_{F,\p}$:] $\mut_F(\p_1\oplus \p_2)=\p'$ with $\p_1/F=\p_2/F=\p'/F'$, for $(F(\p_1),F(\p_2))$ either being $((1,0,0),(1,1,0))$ or $((2,1,1),(2,0,1))$ with $F'(\p')=(o_{F,\p_1},x_{F,\p_1},i_{F,\p_1})$.
\end{itemize}
\end{definition}

Consider the torus $\mathcal X_{\w}$ corresponding to the seed (associated with) $\pa(\w)$. 
For the lattice $N_{\w}$ we have the basis $\{e_E\}_{E \text{ face of }\pa(\w)}$. 
Then $e_{\p}\in N$ for ${\p}\in \mathcal{P}_{\w,F}$ is an expression in this basis and  $z^{e_{\p}}$ a function on $\mathcal X_{\w}$. 
To extend our definition of $e_{\p}$ in \eqref{eq:def area ineq} for $\p\in \mathcal P_{\w}$ to $\p\in\widetilde{\mathcal{P}}_{\w,F}$, we set $z^{e_{{\p}_1\oplus {\p}_2}}:=z^{e_{{\p}_1}}+z^{e_{{\p}_2}}$.
Then for every mutable face $F$ of $\pa(\w)$ we have
\[
\left\{\mathbf x\in \mathbb R^{\ell(w)}\left\vert (\sum_{\p\in\widetilde{\mathcal{P}}_{\w,F}}z^{e_{\p}})^{\trop}(\mathbf x)\ge 0\right.\right\}=S_{\w}.
\]
The following is the key lemma of this section.

\begin{table}
    \centering
    \scalebox{.8}{
    \begin{tabular}{|c|c|c|c|} \hline
$F$-local type of $\p$ & ${\p}$ in $\pa(\w)$
    & $\p'=\mut_{F}({\p})$ in $\pa(\mu_F(\w))$ & $F'$-local type of $\p'$\\ \hline\hline
 $(1,1,1)$ &
    $b_i\to v_{(i,k)}\to v_{(i,j)}\to a_i $ &
    $ b'_i\to v'_{(i,j)}\to v'_{(i,k)}\to a'_i $ &
     $(1,0,1)$ \\ \hline
    
\begin{tabular}{c} 
     $(1,0,0)$ \\ \hline
     $(1,1,0)$ 
    \end{tabular} &
    \begin{tabular}{c}
         $b_i\to v_{(i,k)}\to v_{(j,k)}\to b_{j} $  \\ \hline
         $ b_i\to v_{(i,k)}\to v_{(i,j)}\to v_{(j,k)}\to b_{j} $
    \end{tabular} &
    $b'_i\to v'_{(i,j)}\to b_{k} $ &
     $(0,0,1)$ \\ \hline
    
 $(1,1,1)$ &
    $b_i\to v_{(i,k)}\to v_{(i,j)}\to v_{(j,k)}\to b_k $ &
    $b'_i\to v'_{(i,j)}\to v'_{(i,k)}\to b'_k $ &
     $(1,0,1)$ \\ \hline    
    
 $(1,1,2)$ &
    $ a_j\to v_{(i,j)}\to  a_i $ &
    \begin{tabular}{c}
         $ a'_j\to v'_{(j,k)}\to v'_{(i,k)}\to a'_i $  \\
         $a'_j\to v'_{(j,k)}\to v'_{(i,j)}\to v'_{(i,k)}\to a'_i $ 
    \end{tabular} &
    \begin{tabular}{c} 
     $(2,1,1)$ \\ \hline
     $(2,0,1)$ 
    \end{tabular} \\ \hline

 $(1,1,1)$ &
    $ a_j\to v_{(i,j)}\to v_{(j,k)}\to b_j  $ &
    $ a'_j\to v'_{(j,k)}\to v'_{(i,j)}\to b'_j $ &
     $(1,0,1)$ \\ \hline 

 $(1,1,2)$ &
    $ a_j\to v_{(i,j)}\to v_{(j,k)}\to b_k $ &
    \begin{tabular}{c}
         $ a'_j\to v'_{(j,k)}\to v'_{(i,k)}\to b'_k $  \\
         $ a'_j\to v'_{(j,k)}\to v'_{(i,j)}\to v'_{(i,k)}\to b'_k $ 
    \end{tabular}  &
    \begin{tabular}{c} 
     $(2,1,1)$ \\ \hline
     $(2,0,1)$ 
    \end{tabular}\\ \hline 

 $(1,1,1)$  &
     $ a_k\to v_{(i,k)}\to v_{(i,j)}\to a_i $   &
    $  a'_k\to v'_{(j,k)}\to v'_{(i,j)}\to v'_{(i,k)}\to a'_i $ &
     $(1,0,1)$ \\ \hline  
    
    \begin{tabular}{c} 
     $(1,0,0)$ \\ \hline
     $(1,1,0)$ 
    \end{tabular} &
    \begin{tabular}{c}
         $a_k\to v_{(i,k)}\to v_{(j,k)}\to b_j $  \\
         $a_k\to v_{(i,k)}\to v_{(i,j)}\to v_{(j,k)}\to b_j $ 
    \end{tabular}  &
    $ a'_k\to v'_{(j,k)}\to v'_{(i,j)}\to b'_j $ &
     $(0,0,1)$ \\ \hline  

 $(1,1,1)$  &
    $ a_k\to v_{(i,k)}\to v_{(i,j)}\to  v_{(j,k)}\to b_k $ &
    $  a'_k\to v'_{(j,k)}\to v'_{(i,j)}\to v'_{(i,k)}\to b'_k $ &
     $(1,0,1)$ \\ \hline  
    \end{tabular}}
    \caption{Shapes of paths locally around $F$ (resp. $F'$) in $\mathcal P_{\w}$ (resp. $\mathcal P_{\mu_F(\w)}$) for orientation $(l_r,l_{r+1})$ with $i\le r$ and $r+1\le j<k$ (see Figure~\protect{\ref{fig:loc.orient}}) and how they are mapped onto each other by $\mut_F$.}
    \label{tab:(lrr) cases of path mut}
\end{table}

\vbox{
\begin{lemma}\label{lem: mutation paths}
Let $w\in S_n$ with reduced expressions $\w$ and $\mu_F(\w)$, where $F$ is a mutable face of $\pa(\w)$ and $F'$ the corresponding face of $\pa(\mu_F(\w))$ (i.e. $\mu_{F'}(\mu_F(\w))=\w$).
Let $\{e_E\}_E$ denote the basis for $N_{\w}$ and $\{e'_E\}_E$ the basis for $N_{\mu_F(\w)}$. Then for ${\p}\in \widetilde{\mathcal P}_{\w,F}$ we have
\[
\mu_{F'}^*(z^{e_{\p}})=z^{e'_{\mut_F(\p)}}.
\]
\end{lemma}}

\begin{proof}
We prove the claim case-by-case depending on the $F$-local type of $\p$ as in Lemma~\ref{lem:possible F-local types}.
As notation we use $n\in N_{\w}$ (resp. $n'\in N_{\mu_F(\w)}$) referring to an expression of $n$ is the basis $\{e_E\}_{E\text{ face of }\pa(\w)}$ (resp. $\{e'_E\}_{E\text{ face of }\pa(\mu_F(\w))}$).
Consider ${\p}\in \mathcal P_{\w}$, then 
\[
e_{\p}=-\sum_{E\subset \area_{\p}}e_E=-\sum_{E\subset ( \text{In}_{F}\cup \text{Out}_{F}) \cap \area_{\p}} e_E - \sum_{E\not \subset ( \text{In}_{F}\cup \text{Out}_{F}) \cap \area_{\p}} e_E=: n_{\p}+m_{\p}.
\]
As by definition $\mut_F$ effects a path only locally around $F$, we have $m_{\p}=m_{\mut_F(\p)}$ (resp. $m_{\p}=m_{\p'_1}=m_{\p'_2}$ if $\mut_F(\p)=\p'_1\oplus\p'_2\in\widetilde{\mathcal P}_{\mu_F(\w),F'}$).
Both have the same expressions in bases $\{e_E\}_E$ and $\{e'_E\}_E$ as the corresponding basis elements are not effected by mutation: 
only basis elements corresponding to vertices (i.e. faces of $\pa(\w)$) adjacent to $F$ (i.e. in $\text{In}_F\cup \text{Out}_F$) are changed by mutation in \eqref{eq: def mutation lattice basis}.
We use this fact throughout the proof.
Denote basis elements associated with faces $F_{\init},F_{\init_1},F_{\init_2}\in \text{In}_F$ by $e_{\init},e_{\init_1},e_{\init_2}$ and similarly for $e_{\text{out}}$. After mutation, $e'_{\init}$ is associated with the face $F'_{\init}\in \text{Out}_{F'}$ in $\pa(\mu_F(\w))$. 

\begin{table}
    \centering
    \scalebox{.8}{
    \begin{tabular}{|c|c|c|c|} \hline
$F$-local type of $\p$ & ${\p}$ in $\pa(\w)$
    & $\p'=\mut_{F}({\p})$ in $\pa(\mu_F(\w))$ & $F'$-local type of $\p'$\\ \hline\hline
 $(1,0,1)$ &
    $ a_i\to v_{(i,j)}\to v_{(i,k)}\to b_i $ &
    $ a'_i\to v'_{(i,k)}\to v'_{(i,j)}\to b'_i $ &
     $(1,1,1)$  \\ \hline      

    \begin{tabular}{c} 
     $(2,0,1)$   \\ \hline
     $(2,1,1)$ \end{tabular}&
    \begin{tabular}{c}
         $ a_i\to v_{(i,j)}\to v_{(i,k)}\to v_{(j,k)}\to b_j  $  \\
         $ a_i\to v_{(i,j)}\to v_{(j,k)}\to b_j $ 
    \end{tabular}   &
    $ a'_i\to v'_{(i,k)}\to v'_{(i,j)}\to b'_j $ &
     $(1,1,2)$ \\ \hline 

 $(2,1,2)$  &
    $ a_i\to v_{(i,j)}\to v_{(j,k)}\to b_k $ &
    $ a'_i\to v'_{(i,k)}\to b'_k $ &
     $(2,1,2)$ \\ \hline 

 $(1,0,1)$ &
    $ a_j\to v_{(i,j)}\to v_{(i,k)}\to v_{(j,k)} \to b_j  $ &
    $ a'_j\to v'_{(j,k)}\to v'_{(i,k)}\to v'_{(i,j)}\to b'_j $ &
     $(1,1,1)$ \\ \hline  
    
 $(0,0,0)$  &
    $ a_k\to v_{(i,k)}\to b_i $ &
    $ a'_k\to v'_{(j,k)}\to v'_{(i,j)}\to b'_i $ &
     $(0,0,0)$ \\ \hline  

 $(0,0,1)$  &
    $ a_j\to v_{(i,j)}\to v_{(i,k)}\to b_i $ &
    \begin{tabular}{c}
         $ a'_j\to v'_{(j,k)}\to v'_{(i,j)}\to b'_i $  \\
         $ a'_j\to v'_{(j,k)}\to v'_{(i,k)}\to v'_{(i,j)}\to b'_i $ 
    \end{tabular}   &
    \begin{tabular}{c} 
     $(1,0,0)$ \\ \hline
     $(1,1,0)$ \end{tabular} \\ \hline    

    \end{tabular}}
    \caption{Shapes of paths locally around $F$ (resp. $F'$) in $\mathcal P_{\w}$ (resp. $\mathcal P_{\mu_F(\w)}$) for orientation $(l_r,l_{r+1})$ with $r+1\le i<j<k$ (see Figure~\protect{\ref{fig:loc.orient}}) and how they are mapped onto each other by $\mut_F$.}
    \label{tab:(rrr) cases of path mut}
\end{table}

We distinguish the cases as in Lemma~\ref{lem:possible F-local types}.
\begin{itemize}
    \item[$i_{F,\p}<o_{F,\p}$] From Lemma~\ref{lem:possible F-local types} we know that in this case $n_{\p}=-e_{\init}-e_F-e_{\text{out}_1}-e_{\text{out}_2}$ (resp. $n_{\p}=-e_{\text{out}}$) and $\mut_F({\p})={\p'}_1\oplus {\p'}_2$ with $\p_1',\p_2'$ as in Definition~\ref{def:mut path}.
    Then $n'_{\p'_1}=-e'_{\init}-e'_F-e'_{\text{out}_1}-e'_{\text{out}_2}$ (resp. $n'_{\p'_1}=-e'_{\text{out}}-e_F$) and $n'_{\p'_2}=-e'_{\init}-e'_{\text{out}_1}-e'_{\text{out}_2}$ (resp. $n'_{\p'_1}=-e'_{\text{out}}$). 
    We compute using formulas \eqref{eq: def mutation lattice basis}, \eqref{eq: def pullback X-mut} and the observation that $m'_{\p'_1}=m'_{\p'_2}$:
    \begin{align*}
    \mu_{F'}^*(z^{n_{\p}+m_{\p}}) &= z^{-e_{\init}-e_F-e_{\text{out}_1}-e_{\text{out}_2}+m_{\p}}(1+z^{e_F}) \\
    &= z^{-e'_{\init}-e'_{\text{out}_1}-e'_{\text{out}_2}+m'_{\p'}}(1+z^{-e'_F})\\
    &= z^{n'_{\p'_1}+m'_{\p'_1}}+z^{n'_{\p'_2}+m'_{\p'_2}} 
    =z^{e'_{{\p'}_1}}+z^{e'_{{\p'}_2}}
    \stackrel{\text{(by def.)}}{=} z^{e'_{{\p'}_1\oplus {\p'}_2}}=
    z^{e'_{\mut_F(\p)}}\\
    (\text{resp. }  \mu_{F'}^*(z^{n_{\p}+m_{\p}}) &= z^{-e_{\text{out}}+m_{\p}}(1+z^{e_F})
    = z^{-e'_{\text{out}}+m'_{\p'}}(1+z^{-e'_F})\\
    &= z^{n'_{\p'_1}+m'_{\p'_1}}+z^{n'_{\p'_2}+m'_{\p'_2}} = z^{e'_{\mut_F(\p)}}).
    \end{align*}
    \item[${i_{F,{\p}}=o_{F,{\p}}}$] In this case $\mut_F({\p})={\p'}\in\widetilde{\mathcal P}_{\mu_F(\w),F'}$ as in Definition~\ref{def:mut path}.
    We divide into three cases: $i_{F,{\p}}\in\{0,1,2\}$. 
    If $i_{F,{\p}}=0$, consider $\area_{\p}=F_1\cup\dots \cup F_r$ then $\area_{{\p}'}=F'_1\cup \dots \cup F'_r$. Further, 
    \begin{align*}
    \mu_{F'}^*(z^{e_{\p}})=\mu_{F'}^*(z^{m_{\p}})=z^{m_{\p}}=z^{m_{{\p}'}}=z^{e'_{{\p}'}}=z^{e'_{\mut_F({\p})}}.
    \end{align*}
    If $i_{F,{\p}}=1$ we have $n_{\p}= -e_F-e_{\init}-e_{\text{out}}$ (resp. $n_{\p}= -e_{\init}-e_{\text{out}}$).
    We have $n'_{\p'}= -e'_{\init}-e'_{\text{out}}$ (resp. $n'_{\p'}= -e'_F-e'_{\init}-e'_{\text{out}}$) and compute
    \begin{align*}
    \mu_{F'}^*(z^{n_{\p}})&= z^{-e_F-e_{\init}-e_{\text{out}}} = z^{e'_F-(e'_{\init}+e'_F)-e'_{\text{out}}} = z^{-e'_{\init}-e'_{\text{out}}}=z^{n'_{\p'}}\\
    (\text{resp. }    \mu_{F'}^*(z^{n_{\p}})&= z^{-e_{\init}-e_{\text{out}}} = z^{-(e'_{\init}+e'_F)-e'_{\text{out}}} = z^{-e'_F-e'_{\init}-e'_{\text{out}}}=z^{n'_{\p'}}).
    \end{align*}
    If $i_{F,{\p}}=2$ we have $n_{\p}=-e_{\init_1}-e_{\init_2}-e_F-e_{\text{out}_1}-e_{\text{out}_2}$. Now $n'_{{\p}'}=-e'_{\init_1}-e'_{\init_2}-e'_F-e'_{\text{out}_1}-e'_{\text{out}_2}$ and we compute
    \begin{align*}
    \mu_{F'} ^*(z^{n_{\p}})&= z^{-e_{\init_1}-e_{\init_2}-e_F-e_{\text{out}_1}-e_{\text{out}_2}} = z^{-(e'_{\init_1}+e'_F) - (e'_{\init_2}+e'_F)-(-e'_F)-e'_{\text{out}_1}-e'_{\text{out}_2}}\\
    &= z^{-e'_{\init_1}-e'_{\init_2}-e'_F-e'_{\text{out}_1}-e'_{\text{out}_2}} = z^{n'_{\p'}}.
    \end{align*}
    In all three cases the claim follows from the computation.
    \item[$i_{F,{\p}}>o_{F,{\p}}$] In this case by Lemma~\ref{lem:possible F-local types} there are paths $\p_1,\p_2\in\mathcal P_{\w}$ with $\p_1\oplus \p_2\in\widetilde{\mathcal{P}}_{\w,F}$ and $\mut_{F}(\p_1\oplus\p_2)=\p'\in\widetilde{\mathcal P}_{\mu_F({\w}),F'}$ as in Definition~\ref{def:mut path}. 
    We have $n_{\p_1}=-e_{\init}$ and $n_{\p_2}=-e_{\init}-e_F$ (resp. $n_{\p_1}=-e_{\init_1}-e_{\init_2}-e_{\text{out}}$ and $n_{\p_2}=-e_{\init_1}-e_{\init_2}-e_F-e_{\text{out}}$). For $\p'$ we have $n'_{\p'}=-e'_{\init}$ (resp. $n'_{\p'}=-e'_{\init_1}-e'_{\init_2}-e'_F-e'_{\text{out}}$). We compute
    \begin{align*}
    \mu_{F'}^*(z^{n_{\p_1}}+z^{n_{\p_2}}) &=
    z^{-e_{\init}}(1+z^{e_F})^{-1}+z^{-e_{\init}-e_F}(1+z^{e_F})^{-1} \\
    &= (z^{-e'_{\init}-e'_F}+z^{-e'_{\init}})(1+z^{-e'_F})^{-1} = z^{-e'_{\init}} = z^{n'_{\p'}}\\
    (\text{resp. } \mu_{F'}^*(z^{n_{\p_1}}+z^{n_{\p_2}}) &=
    z^{-e_{\init_1}-e_{\init_2}-e_{\text{out}}}(1+z^{e_F})^{-1}+z^{-e_{\init_1}-e_{\init_2}-e_F-e_{\text{out}}} (1+z^{e_F})^{-1} \\
    &= (z^{-e'_{\init_1}-e'_{\init_2}-2e'_F-e'_{\text{out}}} + z^{-e'_{\init_1}-e'_{\init_2}-e'_F-e'_{\text{out}}})(1+z^{-e'_F})^{-1} \\
    &= z^{-e'_{\init_1}-e'_{\init_2}-e'_F-e'_{\text{out}}} = z^{n'_{\p'}}).
    \end{align*}
    In both cases the claim follows.
\end{itemize}
\end{proof}

Before proving a generalization of Proposition~\ref{prop: superpot area GT} we have to show that the normal vectors associated to the weight inequalities $e_{[i:k]}$ \eqref{eq:def area wt ineq} mutate as expected.
We use the notation as in Lemma~\ref{lem: mutation paths} and its proof.
Recall the normal vectors of the weight inequalities for $\mathcal S_{\w}$ from \eqref{eq:def area wt ineq}.
For $i\in[n-1]$ let $e_{[i:0]},\dots,e_{[i:n_i]}$ be those for $\w$ as expressions in $\{e_E\}_E$ and $e'_{[i:0]},\dots,e'_{[i,n'_i]}$ those for $\mu_F(\w)$ as expressions in $\{e'_E\}_E$.

\begin{lemma}\label{lem:mutation wt area}
With notation as above we have for every $i\in[n-1]$
\[
\mu_{F'}^*\left(\sum_{k=0}^{n_i}z^{e_{[i:k]}}\right)=\sum_{k'=0}^{n'_i}z^{e'_{[i:k']}}.
\]
\end{lemma}
\begin{proof}
We treat the case where $F$ is of level $l$ and $F'$ of level $l+1,$ with $ l\in[n-2]$ (the proof of the other case is similar).
Recall that $e_{[i:k]}=-e_{F_i}-e_{F_{j_1}}-\dots-e_{F_{j_k}}$, where $k\in[0,n_i]$, $s_{j_1},\dots,s_{j_{n_i}}=s_i$ in $\w$, and $F_{j_k}$ is bounded to the left by the crossing in $\pa(\w)$ induced by $s_{j_k}$.
Let $F'_i,F'_{j_1},\dots,F'_{j_{n'_i}}$ be the corresponding faces in $\pa(\mu_F(\w))$.
In particular, if $i\not\in\{l,l+1\}$ we have 
\[
\mu_{F'}^*(z^{-e_{F_i}-e_{F_{j_1}}-\dots-e_{F_{j_k}}})=z^{-e'_{F'_i}-e'_{F'_{j_1}}-\dots-e'_{F'_{j_k}}}.
\]
We therefore focus on the cases $i\in\{l,l+1\}$.
\begin{itemize}
    \item[$i=l$] As $F$ is of level $l$ we have $F=F_{j_k}$ for one $k\in[n_l],s_{j_k}=s_l$. 
    By \eqref{eq: def mutation lattice basis} we have $\mu_F(e_{F_l})=e'_{F'_l}$ and $\mu_{F}(e_{F_{j_r}})=e'_{F'_{j_r}}$ for $r\in[k-1]$, hence
    \[
    \mu_{F'}^*(z^{-e_{F_l}-e_{F_{j_1}}-\dots-e_{F_{j_r}}})=z^{-e'_{F'_l}-e'_{F'_{j_1}}-\dots-e'_{F'_{j_r}}}.
    \]
    Still by \eqref{eq: def mutation lattice basis} we have $\mu_F(e_{F_{k-1}})=e'_{F'_{j_{k-1}}}+e'_{F'_{j_k}}, \mu_F(e_{F_{j_k}})=-e'_{F'_{j_k}}$ and $\mu_F(e_{F_{j_s}})=e'_{F'_{j_s}}$ for $s\in[k+1,n_l]$. Plugging in to \eqref{eq: def pullback X-mut} we obtain
    \begin{align*}
    \mu_{F'}^*(z^{e_{[l:k-1]}}+z^{e_{[l:k]}}+z^{e_{[l:k+1]}})&=  z^{-e_{F_l}-\dots-e_{F_{j_{k-1}}}}(1+z^{e_{F_{j_k}}})^{-1}\\ &\quad\quad + z^{-e_{F_l}-\dots-e_{F_{j_k}}}(1+z^{e_{F_{j_k}}})^{-1} + z^{-e_{F_l}-\dots-e_{F_{j_{k+1}}}}\\
    &= z^{-e'_{F'_l}-\dots-e'_{F'_{j_{k-1}}}} + z^{-e'_{F'_l}-\dots-e'_{F'_{j_{k+1}}}}\\
    &= z^{e'_{[l:k-1]}}+z^{e_{[l:k]}}.
    \end{align*}
    Note that the index shift in the last equality comes from the fact that $\pa(\mu_F(\w))$ has one less face of level $l$ than $\pa(\w)$ as $F'$ is of level $l+1$. So the claim follows for level $l$.
    \item[$i=l+1$] Let $F_{j_r}$ be the face of level $l+1$ in $\text{Out}_F$ and $F_{j_{r+1}}$ the one in $\text{In}_F$. Then we compute with notation as above
    \begin{align*}
    \mu_{F'}^*(z^{e_{[l:r]}}+z^{e_{[l:r+1]}})&= z^{-e_{F_{l+1}}-\dots-e_{F_{j_r}}}(1+z^{e_F})+z^{-e_{F_{l+1}}-\dots-e_{F_{j_{r+1}}}}\\
    &= z^{-e'_{F_{l+1}}-\dots-e'_{F'_{j_r}}} + z^{-e'_{F'_{l+1}}-\dots-e'_{F'_{j_r}}-e'_{F'}} + z^{-e'_{F'_{l+1}}-\dots-e'_{F'_{j_r}}-e'_{F'}-e'_{F'_{j_{r+1}}}}\\
    &= z^{e'_{[l+1:r]}}+z^{e'_{[l+1:r+1]}}+z^{e_{[l+1:r+2]}}.
    \end{align*}
    As before the index shift occurs because $\pa(\mu_F(\w))$ has additionally the face $F'$ of level $l+1$ in comparison to $\pa(\w)$.
\end{itemize}
\end{proof}

We can now prove the following theorem.

\vbox{
\begin{theorem}\label{thm: superpot and area potential}
Let $\w_0$ be an arbitrary reduced expression of $w_0\in S_n$. Then the superpotential expressed in the seed given by $\pa(\w_0)$ satisfies $W\vert_{\X_{\w_0}}=W_{\mathcal S_{\w_0}}$. In particular,
\[
W\vert_{\X_{\w_0}}=\sum_{\p\in\mathcal P_{\w_0}}z^{e_{\p}}+\sum_{i\in[n-1],0\le k\le n_i}z^{e_{[i:k]}}.
\]
\end{theorem}}

\begin{proof}
By Proposition~\ref{prop: superpot area GT} the claim is true for the seed $s_0$ with $\w_0=s_1s_2s_1\cdots s_{n-1}\cdots s_2s_1$. 
Now Lemmas~\ref{lem: mutation paths} and \ref{lem:mutation wt area} imply that the claim holds for all seeds that are related to $s_0$ by a finite sequence of mutations.
As there are only finitely many reduced expressions for $w_0$ and they are all related by mutation as defined in Definition~\ref{defn:mutation pa} the claim is true for all $\w_0$.
\end{proof}

\begin{remark}
The theorem has the following consequence relating GP-paths with $\mathcal X$-cluster variables and $\vartheta$-functions. 
We have encoded GP-paths in two different ways: in terms of their crossing points in the pseudoline arrangement $c_{(i,j)}$'s and in terms of the faces of the pseudoline arrangement they enclose $e_{F_i}$'s.
For a fixed orientation of $\pa(\w_0)$, say $(l_i,l_{i+1})$, the detropicalization of the GP-data is a sum of corresponding $\mathcal X$-variables.
In fact, by Theorem~\ref{thm: superpot and area potential} this sum of $\mathcal X$-variables is the $\vartheta$-function associated to the frozen vertex $w_{(i,n)}$ expressed in the seed $s_{\w_0}$:
every GP-path $\p$ yields one summand $z^{e_\p}$ of $\vartheta_{(i,n)}\vert_{\mathcal X_{\w_0}}$.
\end{remark}

\vbox{
\begin{corollary}\label{cor:area is trop super}
For every reduced expression $\w_0\in S_n$ the following polyhedral objects coincide
\begin{enumerate}[(i)]
    \item $\mathcal S_{\w_0}=\Xi_{\w_0}$,
    \item $S_{\w_0}=\mathsf \Xi_{\w_0}$,
    \item $\mathcal S_{\w_0}(\lambda)=\Xi_{\w_0}(\lambda)$ for $\lambda\in\mathbb R^{n-1}$.
\end{enumerate}
\end{corollary}}

\begin{proof}
The claim in (i) follows immediately from Theorem~\ref{thm: superpot and area potential} by tropicalizing. 
Then (iii) follows by definition as we intersect both cones with the same collection of hyperplanes.
To see (ii), recall from the proof of Proposition~\ref{prop: superpot area GT} that for the initial seed $s_0$ the $\vartheta$-functions $\vartheta_{(i,n)}$ correspond to GP-paths.
Then the claim follows by Lemma~\ref{lem: mutation paths} and the proof of Theorem~\ref{thm: superpot and area potential}.
\end{proof}

%%%%%%%%%%%%%%%%%%%%%%%%%%%%%%%% MUTATION AREA END

\section{Applications of Theorem~\ref{thm:unimod}}\label{subsec:apply}

We have seen in the last two subsections how the cones and polytopes defined in \S\ref{sec:pa and gp} arise from a representation theoretic point of view and in the context of cluster varieties.
The following theorem is the main combinatorial result of this section. We obtain it as an application of the unimodular equivalences in Theorem~\ref{thm:unimod}.

\vbox{
\begin{theorem}\label{thm:application unimod w_0}
Let $\w_0$ be an arbitrary reduced expression of $w_0\in S_n$. Then the following polyhedral objects are unimodularly equivalent
\begin{enumerate}[(i)]
    \item $\mathcal Q_{\w_0} \cong \Xi_{\w_0}$ via $\Psi_{\w_0}$,
    \item $Q_{\w_0}\cong \mathsf \Xi_{\w_0}$ via $\Psi_{\w_0}\vert_{\mathbb R^N}$,
    \item $\mathcal Q_{\w_0}(\lambda)\cong \Xi_{\w_0}(\lambda)$ for $\lambda\in\mathbb R^{n-1}$ via $\Psi_{\w_0}$.
\end{enumerate}
\end{theorem}
}

\begin{proof}
Combine Theorem~\ref{thm:unimod} with Theorem~\ref{thm:wt GP is wt string} and Corollary~\ref{cor:area is trop super}.
\end{proof}

\begin{remark}
For the special case of the initial seed $s_0$ the theorem can also be proved by combining results of Magee and Littelmann. 
In \cite{Lit98} Littelmann shows that the string polytope $\mathcal Q_{\w_0}(\lambda)$ for $\w_0=s_1s_2s_1\cdots s_{n-1}s_{n-2}\cdots s_2s_1$ is unimodularly equivalent to the \emph{Gelfand-Tsetlin polytope} defined in \cite{GT50}.
Magee shows in \cite{Mag15} that $\Xi_{s_0}$ (resp. $\Xi_{s_0}(\lambda)$) is unimodularly equivalent to the Gelfand-Tsetlin cone (resp. polytope). 
Combining both, one obtains Theorem~\ref{thm:application unimod w_0} for $s_0$.
In fact, to understand Magee's result was driving motivation behind this project.
\end{remark}

By the construction of toric varieties associated to polytopes as in \cite[\S2.1 and \S2.3]{CLS11} and the toric degenerations of Caldero \cite{Cal02} and Gross-Hacking-Keel-Kontsevich \cite{GHKK14} we obtain the following corollary from Theorem~\ref{thm:application unimod w_0} relating these toric varieties.
It is the main result regarding toric degenerations of flag varieties in this section and an answer to Question~\ref{Q:cluster} in the introduction.

\begin{corollary}
Let $\w_0$ be an arbitrary reduced expression of $w_0\in S_n$ and $\lambda\in\mathbb Z_{>0}^{n-1}$. We have an induced isomorphism of the following toric varieties that are degenerations (resp. normalizations of such) of $SL_n/B$
\[
X_{\mathcal Q_{\w_0}}(\lambda)\cong X_{\Xi_{\w_0}}(\lambda).
\]
\end{corollary}

In order to achieve a similar result for Schubert varieties, we study the restriction of the superpotential in the following subsection.

\subsection{Restricted Superpotential and Schubert varieties}
Caldero's degeneration works more generally for Schubert varieties. 
As we have seen above, he uses the degeneration for the flag variety and by a quotient construction on the level of rings he obtains a family for the Schubert variety.
For the cones, taking this quotient corresponds to setting certain variables to zero, or equivalently, restricting the defining GP-paths as in Definition~\ref{def: res path}.
In a similar fashion we want to proceed with the superpotential.
We show how the polytopes defining toric degenerations of Schubert varieties arise in the setting of \cite{GHKK14}.

Consider $w\in S_n$ with reduced expression $\w$ and extension $\w_0=\w s_{i_{\ell(w)+1}}\cdots s_{i_N}$. 
Recall that for a seed corresponding to $\w_0$ we have a basis $\{e_F\mid F \text{ face of }\pa(\w_0)\}$ for $N_{\w_0}$ and further $\mathbb C[\X_{\w_0}]=\mathbb C[z^{\pm e_F}\mid F\text{ face of }\pa(\w_0)]$.
Then $\{e_F\mid F\text{ face of }\pa(\w)\}$ generates a sublattice in $N_{\w_0}$, which we denote by $N_{\w}$ with dual lattice $M_{\w}$ a quotient of $M_{\w_0}$.
We have the torus $\X_{\w}=T_{M_{\w}}=\Spec(\mathbb C[N_{\w}])$ associated with $M_{\w}$ as in \eqref{def: seed tori}.
In particular, $\mathbb C[\X_{\w}]=\mathbb C[z^{\pm e_{F}}\mid F\text{ face of }\pa(\w)]$ and we have a restriction morphism between the Laurent polynomial rings
\[
\res_{\w}:\mathbb C[\X_{\w_0}]\to\mathbb C[\X_{\w}], \ \ f\mapsto f\vert_{\X_{\w}}.
\]

We are interested in the restrictions to $\X_{\w}$ of the superpotential $W\vert_{\X_{\w_0}}$ and the detropicalization $W_{\mathcal S_{\w_0}}$ of $\mathcal S_{\w_0}$ (they are equal by Theorem~\ref{thm: superpot and area potential}). We want to show that they coincide with the detropicalization of $\mathcal S_{\w}$.
In analogy with Definiton~\ref{def:super cone} for $w_0$ we consider for arbitrary $w$ the following polyhedral objects.

\begin{definition}\label{def: res super cone}
For $w\in S_n$ with reduced expression $\w$ and an extension $\w_0=\w s_{i_{\ell(w)+1}}\cdots s_{i_N}$  polyhedral objects by tropicalizing the restriction of a sum of $\vartheta$-functions resp. the superpotential: 
\begin{align*}
\res_{\w}(\Xi_{\w_0}) & :=  \{ \mathbf x\in \mathbb R^{\ell(w)+n-1}\mid \res_{\w}(W\vert_{\X_{\w_0}})^{\trop}(\mathbf x)\ge 0 \},\\
\res_{\w}(\mathsf{\Xi}_{\w_0}) & :=  \{ \mathbf x\in\mathbb R^{\ell(w)}\mid \res_{\w}(\sum_{i=1}^{n-1}\vartheta_{(i,n)}\vert_{\X_{\w_0}})^{\trop}(\mathbf x)\ge 0 \},\\
\res_{\w}(\Xi_{\w_0}(\lambda)) & :=   \res_{\w}(\Xi_{\w_0})  \cap\tau_{\w}^{-1}(\lambda) \text{ for } \lambda\in\mathbb R^{n-1}.
\end{align*}    
\end{definition}

\begin{center}
\begin{figure}
\centering
\begin{tikzpicture}[scale=.7]
\draw[rounded corners] (0,0) -- (1,0) -- (4,3) -- (9,3);
\draw[rounded corners]  (0,1) -- (1,1) -- (2,0) -- (5,0) -- (6,1) -- (7,1) -- (8,2) -- (9,2);
\draw[rounded corners]  (0,2) -- (2,2) -- (3,1) -- (4,1) -- (5,2) -- (7,2) -- (8,1) -- (9,1);
\draw[rounded corners]  (0,3) -- (3,3) -- (6,0) -- (9,0);

\node at (-.5,0) {$l_1$};
\node at (-.5,1) {$l_2$};
\node at (-.5,2) {$l_3$};
\node at (-.5,3) {$l_4$};

\node at (.5,.5) {$F_{1}$};
\node at (1,1.5) {$F_{2}$};
\node at (2,2.5) {$F_{3}$};
\node at (3.5,0.5) {$F_{(1,2)}$};
\node at (3.5,1.5) {$F_{(1,3)}$};
\node at (5.5,1.5) {$F_{(3,4)}$};
\node at (6.5,.5) {$F_{(2,4)}$};
\node at (4.5,2.5) {$F_{(1,4)}$};
\node[red] at (8.5,1.5) {$F_{(2,3)}$};

\draw[red, dashed, thick] (7.25,-.8) -- (7.25,3.5);

\node at (1.5,-.5) {$s_1$};
\node at (2.5,-.5) {$s_2$};
\node at (3.5,-.5) {$s_3$};
\node at (4.5,-.5) {$s_2$};
\node at (5.5,-.5) {$s_1$};
\node at (7.75,-.5) {$s_2$};
\end{tikzpicture}
\caption{Restriction/Extension of a pseudoline arrangement.}\label{fig: restr.pseudo}
\end{figure}
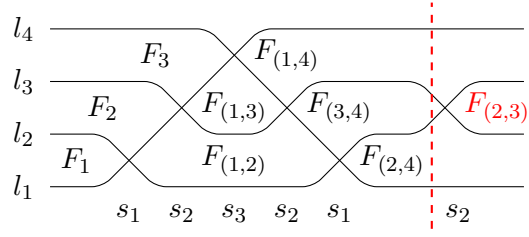
\end{center}

\begin{example}\label{exp: res superpot}
Consider $\w=s_1s_2s_3s_2s_1\in S_4$ with extension $\w_0=\w s_2$. 
We compute the superpotential in $W{\vert}_{\X_{\w_0}}\in\mathbb C[\X_{\w_0}]$.
\begin{align*}
W\vert_{\X_{\w_0}}&= (z^{-e_3}+z^{-e_3-e_{(1,4)}})+ (z^{-e_2} + z^{-e_2-e_{(1,3)}} + z^{-e_2-e_{(1,3)}-e_{(3,4)}})\\
&+ (z^{-e_1}+z^{-e_1-e_{(1,2)}}+z^{-e_1-e_{(1,2)}-e_{(2,4)}}) + (z^{-e_{(2,4)}}+z^{-e_{(2,4)}-e_{(3,4)}}) + (z^{-e_{2,3}})\\
&+ (z^{-e_{(1,4)}}+z^{-e_{(1,4)}-e_{(1,3)}}+z^{-e_{(1,4)}-e_{(1,3)}-e_{(3,4)}}+z^{-e_{(1,4)}-e_{(1,3)}-e_{(3,4)}-e_{(1,2)}}).
\end{align*}
From Figure~\ref{fig: restr.pseudo} we see that $F_{(2,3)}$ is a face of $\pa(\w_0)$, but not of $\pa(\w)$. Hence,
\begin{align*}
{\res}_{\w}(W\vert_{\X_{\w_0}}) &=  (z^{-e_3}+z^{-e_3-e_{(1,4)}})+ (z^{-e_2} + z^{-e_2-e_{(1,3)}} + z^{-e_2-e_{(1,3)}-e_{(3,4)}})\\
&+ (z^{-e_1}+z^{-e_1-e_{(1,2)}}+z^{-e_1-e_{(1,2)}-e_{(2,4)}}) + (z^{-e_{(2,4)}}+z^{-e_{(2,4)}-e_{(3,4)}}) \\
&+ (z^{-e_{(1,4)}}+z^{-e_{(1,4)}-e_{(1,3)}}+z^{-e_{(1,4)}-e_{(1,3)}-e_{(3,4)}}+z^{-e_{(1,4)}-e_{(1,3)}-e_{(3,4)}-e_{(1,2)}}).
\end{align*}
\end{example}

\begin{proposition}\label{prop: res superpt}
Let $w\in S_n$ and consider a reduced expression $\w$ with an extension to $\w_0=\w s_{i_{\ell(w)+1}}\cdots s_{i_N}$. Then
\[
\res_{\w}(W\vert_{\X_{\w_0}})=W_{\mathcal S_{\w}}.
\]
\end{proposition}

\begin{proof}
Recall the restriction of GP-paths defined in Definition~\ref{def: res path}. 
By Propositions~\ref{prop:respath} and \ref{prop:restrict} we have seen $\res_{\w}(\mathcal P_{\w_0})=\mathcal P_{\w}$. 
To avoid confusion we denote as before for $i\in[n-1]$ by $n_i^{\w}:=\#\{j\mid s_{i_j}=s_i \text{ in }\w\}$ and $n_i^{\w_0}:=\#\{j\mid s_{i_j}=s_i\text{ in }\w_0\}$.
Using Theorem~\ref{thm: superpot and area potential} we compute
\begin{align*}
\res_{\w}(W\vert_{\X_{\w_0}}) &= \sum_{\p\in\mathcal P_{\w_0}}z^{e_{\p}}\vert_{\X_{\w}}+\sum_{i\in[n-1],0\le k\le n_i^{\w_0}}z^{e_{[i:k]}}\vert_{\X_{\w}}\\
&= \sum_{\p\in\res_{\w}(\mathcal P_{\w_0})}z^{e_{\p}}+\sum_{i\in[n-1],0\le k\le n_i^{\w}}z^{e_{[i:k]}} = W_{\mathcal S_{\w}}.
\end{align*}
\end{proof}

\vbox{
\begin{theorem}\label{thm: application unimod w}
Let $w\in S_n$ and consider a reduced expression $\w$ with an extension to $\w_0=\w s_{i_{\ell(w)+1}}\cdots s_{i_N}$. Then the following polyhedral objects are unimodularly equivalent
\begin{enumerate}[(i)]
    \item $\mathcal Q_{\w}\cong \res_{\w}(\Xi_{\w_0})$ via $\Psi_{\w}$,
    \item $Q_{\w}\cong \res_{\w}(\mathsf\Xi_{\w_0})$ via $\Psi_{\w}\vert_{\mathbb R^{\ell(w)}}$,
    \item $\mathcal Q_{\w}(\lambda)\cong  \res_{\w}(\Xi_{\w_0}(\lambda))$ for $\lambda\in\mathbb R^{n-1}$ via $\Psi_{\w}$. 
\end{enumerate}
\end{theorem}}

\begin{proof}
For (i) combine Proposition~\ref{prop: res superpt} with Theorem~\ref{thm:unimod} and Theorem~\ref{thm:wt GP is wt string}, which directly implies (iii).
To see (ii), recall that by Lemma~\ref{lem: mutation paths} and the proof of Proposition~\ref{prop: superpot area GT} we have
\[
\sum_{\p\in\mathcal P_{\w_0}} z^{e_{\p}} = \sum_{i\in[n-1]}\vartheta_{(i,n)}\vert_{\X_{\w_0}}.
\]
By the proof of Proposition~\ref{prop: res superpt} the same equality when replacing $\w_0$ by $\w$. Then the claim follows by Theorem~\ref{thm:unimod} and Theorem~\ref{thm:wt GP is wt string}.
\end{proof}

Before stating the following corollary relating the toric degenerations of Schubert varieties by Caldero \cite{Cal02} to the toric degenerations of flag varieties by Gross-Hacking-Keel-Kontsevich \cite{GHKK14} we briefly remind you about the \emph{Orbit-Cone-Correspondence} for toric varieties (see \cite[\S3.2]{CLS11}).
For a (full-dimensional) polytope $P\subset \mathbb R^n$ denote by $\Sigma_P\subset\mathbb R^n$ its \emph{normal fan} (see \cite[Remark~2.3.3]{CLS11}).
Every cone $\sigma\in\Sigma_P$ corresponds to a torus orbits in $X_{\Sigma_{P}}$ of dimension $n-\dim \sigma$ (\cite[Theorem~3.2.6]{CLS11}).
The closure of each torus orbit is a toric variety.
For a face $Q$ of $P$ let $\sigma_Q\in\Sigma_P$ be the cone in $\Sigma_P$ spanned by the normal vectors of all facets of $P$ containing $Q$.
Then by \cite[Proposition~3.2.9]{CLS11} the toric variety $X_Q$ is isomorphic to the closure of the torus orbit corresponding to the cone $\sigma_Q\in\Sigma_P$.

Consider an arbitrary $w\in S_n$ with a reduced expression $\w$ and $\w_0=\w s_{i_{\ell(w)+1}}\cdots s_{i_N}$ an extension. 
For every $\lambda\in\Lambda^{++}$ recall that the toric variety $X_{\mathcal Q_{\w}(\lambda)}$ is (the normalization of) a toric degeneration of $X_w$ by \cite{Cal02}. Similarly, $X_{\Xi_{\w_0}(\lambda)}$ is a flat degeneration of $SL_n/B$ by \cite{GHKK14}. 
We can now formulate the geometric version of our main result on toric degenerations of Schubert varieties.

\begin{corollary}
The toric variety $X_{\mathcal Q_{\w}(\lambda)}$ is isomorphic to a subvariety of $X_{\Xi_{\w_0}(\lambda)}$.
More precisely, we have
\[
X_{\mathcal Q_{\w}(\lambda)}\cong X_{\res_{\w}(\Xi_{\w_0}(\lambda))},
\]
where $ X_{\res_{\w}(\Xi_{\w_0}(\lambda))}$ is the closure of the torus orbit corresponding to the cone $\sigma_{\res_{\w}(\Xi_{\w_0}(\lambda))}\in \Sigma_{\Xi_{\w_0}(\lambda)}$.
\end{corollary}

\begin{proof}
By definition $\res_{\w}(\Xi_{\w_0}(\lambda))$ is a union of faces of $\Xi_{\w_0}(\lambda)$. Theorem~\ref{thm: application unimod w}(iii) implies in particular, that $\res_{\w}(\Xi_{\w_0}(\lambda))$ is a polytope itself, hence a face of $\Xi_{\w_0}(\lambda)$.
Further, the unimodular equivalence $Q_{\w}(\lambda)\cong \res_{\w}(\Xi_{\w_0}(\lambda))$ induces an isomorphism of toric varieties $X_{Q_{\w}(\lambda)}\cong X_{\res_{\w}(\Xi_{\w_0}(\lambda))}$.
Then the Corollary follows by \cite[Proposition~3.2.9]{CLS11}.
\end{proof}

\subsection{Restriction vs. superpotential for double Bruhat cells} 
We conclude with an example that shows how $\res_{\w}(W\vert_{\X_{\w_0}})$ is essentially different from a function one would obtain from applying Algorithm~\ref{alg:superpot via opt seeds} to the quiver $Q_{\w}$

\begin{example}
Let $s = s_{\w}$ be the seed of the reduced expression $\w=s_1s_2s_3s_2s_1 \in S_4$ as in Figure~\ref{fig: restr.pseudo}. The corresponding quiver is pictured in Figure~\ref{restr.quiver}. We apply Algorithm~\ref{alg:superpot via opt seeds} and compute optimized seeds for all frozen vertices in $Q_{\w}$.
As $w_{3}$ and $w_{(2,4)}$ are sinks in $Q_{\w}$ we set $\vartheta_{3}\vert_{\X_s}=z^{-e_{3}}$ and $\vartheta_{(2,4)}\vert_{\X_s}=z^{-e_{(2,4)}}$, where $\{e_1,e_2,e_3,e_{(1,2)},e_{(1,3)},e_{(1,4)},e_{(2,4)},e_{(3,4)}\}$ is the lattice basis associated to $s$.

\begin{center}
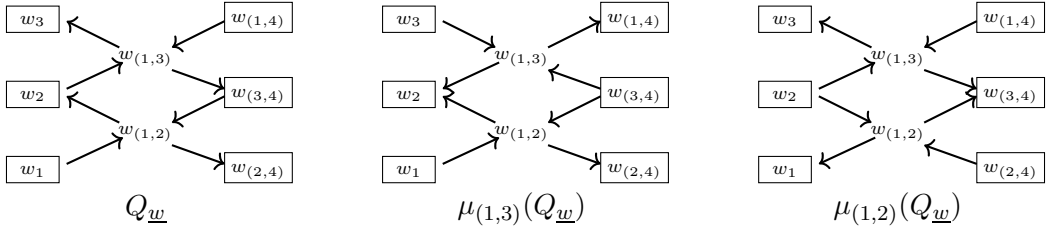
\begin{figure}[ht]
\centering
\begin{tikzpicture}

\node at (0,0) {\tiny$w_{1}$};
    \draw (-.35,-.15) rectangle (.35,0.2);
\node at (3,0) {\tiny $w_{(2,4)}$};
    \draw (2.55,-.15) rectangle (3.45,0.25);
\node at (1.5,0.5) {\tiny $w_{(1,2)}$};
\node at (0,1) {\tiny $w_{2}$};
    \draw (-.35,.85) rectangle (.35,1.2);
\node at (3,1) {\tiny $w_{(3,4)}$};
        \draw (2.55,.85) rectangle (3.45,1.25);
\node at (1.5,1.5) {\tiny $w_{(1,3)}$};
\node at (0,2) {\tiny $w_{3}$};
    \draw (-.35,1.85) rectangle (.35,2.2);
\node at (3,2) { \tiny $w_{(1,4)}$};
        \draw (2.55,1.85) rectangle (3.45,2.25);

%% von und zu 23        
\draw[->, thick] (1.15,1.65) -- (0.45,2);
\draw[->, thick] (2.55,2) -- (1.85,1.65);
\draw[->, thick] (0.45,1.1) -- (1.2,1.45);
\draw[->, thick] (1.85,1.35) -- (2.55,1.1);

%% von und zu 13
\draw[->, thick] (1.15,0.65) -- (0.45,1);
\draw[->, thick] (2.55,1) -- (1.85,.65);
\draw[->, thick] (0.45,.1) -- (1.2,.45);
\draw[->, thick] (1.85,.35) -- (2.55,.1);
\node at (1.5,-.5) {$Q_{\w}$};
\begin{scope}[xshift=5cm]

\node at (0,0) {\tiny $w_{1}$};
    \draw (-.35,-.15) rectangle (.35,0.2);
\node at (3,0) {\tiny $w_{(2,4)}$};
    \draw (2.55,-.15) rectangle (3.45,0.25);
\node at (1.5,0.5) {\tiny $w_{(1,2)}$};
\node at (0,1) {\tiny $w_{2}$};
    \draw (-.35,.85) rectangle (.35,1.2);
\node at (3,1) {\tiny $w_{(3,4)}$};
        \draw (2.55,.85) rectangle (3.45,1.25);
\node at (1.5,1.5) {\tiny $w_{(1,3)}$};
\node at (0,2) {\tiny $w_{3}$};
    \draw (-.35,1.85) rectangle (.35,2.2);
\node at (3,2) {\tiny $w_{(1,4)}$};
        \draw (2.55,1.85) rectangle (3.45,2.25);

%% von und zu 23        
\draw[<-, thick] (1.15,1.65) -- (0.45,2);
\draw[<-, thick] (2.55,2) -- (1.85,1.65);
\draw[<-, thick] (0.45,1.1) -- (1.2,1.45);
\draw[<-, thick] (1.85,1.35) -- (2.55,1.1);

%% von und zu 13
\draw[->, thick] (1.15,0.65) -- (0.45,1);
\draw[->, thick] (2.55,1) -- (1.85,.65);
\draw[->, thick] (0.45,.1) -- (1.2,.45);
\draw[->, thick] (1.85,.35) -- (2.55,.1);

\node at (1.5,-.5) {$\mu_{(1,3)}(Q_{\w})$};

\begin{scope}[xshift=5cm]
\node at (0,0) {\tiny $w_{1}$};
    \draw (-.35,-.15) rectangle (.35,0.2);
\node at (3,0) {\tiny $w_{(2,4)}$};
    \draw (2.55,-.15) rectangle (3.45,0.25);
\node at (1.5,0.5) {\tiny $w_{(1,2)}$};
\node at (0,1) {\tiny $w_{2}$};
    \draw (-.35,.85) rectangle (.35,1.2);
\node at (3,1) {\tiny $w_{(3,4)}$};
        \draw (2.55,.85) rectangle (3.45,1.25);
\node at (1.5,1.5) {\tiny $w_{(1,3)}$};
\node at (0,2) {\tiny $w_{3}$};
    \draw (-.35,1.85) rectangle (.35,2.2);
\node at (3,2) {\tiny $w_{(1,4)}$};
        \draw (2.55,1.85) rectangle (3.45,2.25);
        
%% von und zu 23        
\draw[->, thick] (1.15,1.65) -- (0.45,2);
\draw[->, thick] (2.55,2) -- (1.85,1.65);
\draw[->, thick] (0.45,1.1) -- (1.2,1.45);
\draw[->, thick] (1.85,1.35) -- (2.55,1.1);

%% von und zu 13
\draw[<-, thick] (1.15,0.65) -- (0.45,1);
\draw[<-, thick] (2.55,1) -- (1.85,.65);
\draw[<-, thick] (0.45,.1) -- (1.2,.45);
\draw[<-, thick] (1.85,.35) -- (2.55,.1);

\node at (1.5,-.5) {$\mu_{(1,2)}(Q_{\w})$};
\end{scope}
\end{scope}

\end{tikzpicture}
\caption{The quivers $Q_{\w}$, $\mu_{(1,3)}(Q_{\w})$ and $\mu_{(1,2)}(Q_{\w})$ for $\w=s_1s_2s_3s_2s_1$. The boxes denote frozen variables.}\label{restr.quiver}
\end{figure}
\end{center}
\vspace{-.5cm}

For the other variables we have to find a mutation sequence to an optimized seed. Mutation at $w_{(1,3)}$ (resp. $w_{(1,2)}$) yields the quiver $\mu_{(1,3)}(Q_{\w})$ (resp. $\mu_{(1,2)}(Q_{\w})$) in Figure~\ref{restr.quiver}. 
The seed $\mu_{(1,3)}(s)$ is optimized for $w_{(1,4)}$ and $w_{2}$, so $\vartheta_{(1,4)}\vert_{\X_{\mu_{(1,3)}(\w)}}=z^{-e'_{(1,4)}}$ and $\vartheta_{2}\vert_{\X_{\mu_{(1,3)}(\w)}}=z^{-e'_{2}}$. 
In $\X_{\w}$ we obtain $\vartheta_{(1,4)}\vert_{\X_s}=z^{-e_{(1,4)}}+z^{-e_{(1,4)}-e_{(1,3)}}$ and $\vartheta_{2}\vert_{\X_{s}}=z^{-e_{2}}+z^{-e_{2}-e_{(1,3)}}$. 
Proceeding analogously with $\mu_{(1,2)}(s)$, optimized for $w_{(3,4)}$ and $w_{1}$, we obtain a function on $\X_{\w}$
\begin{align*}
F &:= (z^{-e_3}) + (z^{-e_{2}}+z^{-e_2-e_{(1,3)}})
+ (z^{-e_1}+ z^{-e_1-e_{(1,2)}}) + (z^{-e_{(2,4)}})\\
&+ (z^{-e_{(3,4)}}+z^{-e_{(3,4)}-e_{(1,2)}}) + (z^{-e_{(1,4)}}+z^{-e_{(1,4)}-e_{(1,3)}}).
\end{align*}

Comparing to Example~\ref{exp: res superpot} where $\w_0=\w s_2$ we observe that $F\not = \res_{\w}(W\vert_{\X_{\w_0}})$. 
Tropicalizing $\res_{\w}(W\vert_{\X_{\w_0}})$ we get the following set of inequalities defining the cone $\mathcal S_{\w}\subset \mathbb R^8$
\begin{align*}
-x_3 &\ge 0, &-x_3-x_{(1,4)}   &\ge 0, & & \\
-x_2 &\ge 0, &-x_{2}-x_{(1,3)} &\ge 0, &-x_2-x_{(1,3)}-x_{(3,4)} &\ge 0\\
-x_1 &\ge 0, &-x_{1}-x_{(1,2)} &\ge 0, &-x_{1}-x_{(1,2)}-x_{(2,4)} &\ge 0, \\
-x_{(2,4)} &\ge 0, &-x_{(2,4)}-x_{(3,4)} &\ge 0, &-x_{(1,4)}-x_{(1,3)}-x_{(3,4)} &\ge 0,\\
-x_{(1,4)} &\ge 0, &-x_{(1,4)}-x_{(1,3)} &\ge 0, &-x_{(1,4)}-x_{(1,3)}-x_{(3,4)}-x_{(1,2)} &\ge0. 
\end{align*}
From $F^{\trop}$ we get inequalities defining a cone $\mathcal D_{F}\subset\mathbb R^{8}$:
\begin{align*}
-x_3 &\ge 0, & &\\
-x_2 &\ge 0, &-x_2-x_{(1,3)} &\ge 0,\\
-x_1 &\ge 0, &-x_1-x_{(1,2)} &\ge 0,\\
-x_{(2,4)} &\ge 0, & & \\
-x_{(3,4)}& \ge 0, &-x_{(3,4)}-x_{(1,2)} &\ge 0,\\
-x_{(1,4)} &\ge 0, &-x_{(1,4)}-x_{(1,3)} &\ge 0.
\end{align*}

Observe that $\mathcal D_{F}\subset \mathcal S_{\w}$. We compute the polytopes $\mathcal S_{\w}(\lambda)$ and $\mathcal D_{F}\cap \tau_{\w}^{-1}(\lambda)$ for $\lambda=(1,1,1)$ and their lattice points using \emph{polymake}\cite{GJ00}. The outcome is
\[
\vert \mathcal S_{\w}(\lambda)\cap \mathbb Z^8\vert = 49 = \dim_{\mathbb C}H^0(X_w,L_\lambda) > \vert \mathcal D_{F}\cap \tau_{\w}^{-1}(\lambda)\cap\mathbb Z^{8}\vert =30.
\]
In particular, the toric variety $X_{\mathcal D_{F}\cap \tau_{\w}^{-1}(\lambda)}$ can not be a flat degeneration of the Schubert variety $X_w$.
However, this observation is not too surprising from a geometric point of view, as the restricted superpotential and the function $F$ correspond to different partial compactifications of the $\A$-cluster variety $G^{e,w}$ associated with $\mathcal Y(s_{\w})$.
When considering the restricted superpotential, the cluster variety we are dealing with is $G^{e,w_0}$ and its compactification $\bar G^{e,w_0}$ with boundary divisors
\[
\{\bar p_1=0\},\{\bar p_{12}=0\},\{\bar p_{123}=0\},\{\bar p_{4}=0\},\{\bar p_{34}=0\},\{\bar p_{234}=0\}.
\]
Recall that $G^{e,w_0}$ is $SL_4/U$ up to codimension 2.
The Schubert variety of our interest is $X_{w}$ with $s_1s_2s_3s_2s_1=w$. 
It is given by $\{\bar p_{34}=0\}$ as a subvariety $SL_4/B$. 
Note that in fact, whenever we have a reduced expression $\w$ and an extension $\w_0=\w s_{i_{\ell(w)+1}}\cdots s_{i_N}$, then the Pl\"ucker coordinates that appear as $\A$-cluster variables for faces of $\pa(\w_0)$ that are not faces of $\pa(\w)$ vanish identically on $X_w$.
When restricting the superpotential, we consider the divisor of $\bar G^{e,w_0}$ (resp. $SL_4/U$) given by $\{\bar p_{34}=0\}$, which is closely related to $X_w$.

The function $F$ on the other hand corresponds to the $\A-$cluster variety $G^{e,w}$ and its partial compactification $\bar G^{e,w}$ with boundary divisors
\[
\{\bar p_{1}=0\},\{\bar p_{12}=0\},\{\bar p_{123}=0\},\{\bar p_4=0\},\{\bar p_{24}=0\},\{\bar p_{234}=0\}.
\]
In this case, the defining equation for $X_w$ in $SL_4/B$ is not part of the boundary, so there is no reason to expect information for the Schubert variety from the potential $F$ encoding this boundary.
\end{example}

\bibliographystyle{alpha} 
\bibliography{Trop.bib}

\end{document}